\documentclass[fleqn,10pt]{article}
\usepackage{latexsym, graphicx, epsfig, amsmath, amssymb,amsfonts}
\usepackage{natbib,amsthm,version}
\usepackage{amsbsy,bm,multirow,enumerate}
\usepackage[titletoc,page]{appendix}
\usepackage[mathscr]{eucal}
\usepackage{mathtools}
\usepackage{color}
\usepackage[utf8]{inputenc}
\usepackage[english]{babel}
\usepackage{amsthm}
\usepackage{enumerate}
\usepackage[hidelinks]{hyperref}
\usepackage{url}
\usepackage{subfigure}
\usepackage[lined,boxed,linesnumbered,ruled]{algorithm2e}
\usepackage{float}
\usepackage{arydshln}
\usepackage{dsfont,bbm}

\DeclareMathOperator*{\argmin}{arg\,min}

\newcommand{\mbs}[1]{\mathbf{#1}}

\newtheorem{remark}{Remark}[section]

\theoremstyle{definition}

\title{
  Numerical Approximation of Partial Differential Equations
  by a Variable Projection Method with
  Artificial Neural Networks
} 
\author{
  Suchuan Dong\thanks{Author of correspondence. Email: sdong@purdue.edu}, \ Jielin Yang
  \\
  Center for Computational and Applied Mathematics \\
  Department of Mathematics \\
  Purdue University, USA 
 } 

\date{(January 24, 2022)}

\begin{document}
\maketitle


\begin{abstract}

  We present a method for solving linear and nonlinear partial
  differential equations (PDE) based on the variable projection framework
  and artificial neural networks.
  For linear PDEs, enforcing the boundary/initial value
  problem on the collocation points
  gives rise to a separable nonlinear
  least squares problem about the network coefficients. We reformulate this problem
  by the variable projection approach
  to eliminate the linear output-layer coefficients, leading to
  a reduced problem about the hidden-layer coefficients only.
  The reduced problem is solved first by the nonlinear least squares method
  to determine the hidden-layer coefficients, and then the output-layer
  coefficients are computed by the linear least squares method.
  For nonlinear PDEs, enforcing the boundary/initial value problem
  on the collocation points gives rise to a nonlinear least squares
  problem that is not separable, which precludes 
  the variable projection strategy for such problems.
  To enable the variable projection approach for nonlinear PDEs,
  we first linearize the problem with a Newton iteration, using
  a particular linearization formulated in terms of the updated approximation
  field. The linearized system is solved by the variable projection
  framework together with artificial neural networks.
  Upon convergence of the Newton iteration, the neural-network coefficients
  provide the representation of the solution field to the original
  nonlinear problem. We present ample numerical examples
  with linear and nonlinear PDEs to demonstrate the 
  performance of the method developed herein. For smooth field solutions,
  the errors of the current
  method decrease exponentially as the number of collocation points
  or the number of output-layer coefficients increases. 
  We compare extensively the current method with the extreme learning
  machine (ELM) method from a previous work.
  Under identical conditions and network configurations, 
  the current method  exhibits an accuracy significantly superior to the ELM method.

\end{abstract}


\vspace{0.05cm}



\section{Introduction}
\label{sec:intro}

This work concerns the numerical approximation of partial
differential equations (PDE) with artificial neural networks (ANN),
and we explore the use of the variable projection (VarPro)
approach~\cite{GolubP1973,GolubP2003} together with ANNs
for solving linear and nonlinear PDEs.
Neural network-based PDE solvers, especially those based on
deep neural networks (DNN) and deep learning~\cite{GoodfellowBC2016},
have flourished in the past few years;
see e.g.~\cite{SirignanoS2018,RaissiPK2019,EY2018,HeX2019,LuoY2020,ZangBYZ2020,DongN2020,Samaniegoetal2020,WangYP2020,MaoJK2020,LuMMK2021,LiangJHY2021}, and
the recent review~\cite{Karniadakisetal2021} and
the references therein.
%
%
The DNN-based PDE solvers are fairly straightforward
to implement, by encoding the PDEs,
the boundary and initial conditions into a cost function and then using
some flavor of gradient descent (or back propagation) type 
optimization algorithms to minimize this cost.
Their weakness lies in the limited accuracy and the high
computational cost (long network-training time).
Another promising class of neural network-based methods 
for computational PDEs has recently
appeared~\cite{DongL2021,DongL2021bip,DwivediS2020,FabianiCRS2021,DongY2021},
which are based on a type of randomized
neural networks called extreme learning machines (ELM)~\cite{HuangZS2006,HuangCS2006}.
With these methods the weight/bias coefficients in the hidden layers of
the neural network are set to random values and are fixed.
Only the coefficients of the linear output layer are trainable, and
they are trained by a linear least squares method for linear PDEs
and by a nonlinear least squares method for nonlinear PDEs~\cite{DongL2021}.
It has been shown in~\cite{DongL2021} that the accuracy and the computational
cost (network training time) of the ELM-based method are
considerably superior to those of the aforementioned DNN-based
PDE solvers. In addition, the computational performance  of the ELM-type
method from~\cite{DongL2021} is observed to be comparable to
or exceed that of the classical finite element method (FEM).
Further extensions and improvements to the ELM method of~\cite{DongL2021}
have been documented in~\cite{DongY2021} recently,
which compares systematically the improved method with
the classical and high-order FEM for solving a number of
linear and nonlinear PDEs.
The improved ELM method outperforms the classical second-order FEM
by a considerable margin, and it outcompetes the high-order FEM
when the problem size is not very small~\cite{DongY2021}.


Variable projection (VarPro) is a classical approach for solving separable
nonlinear least squares (SNLLS) problems~\cite{GolubP1973,GolubP2003}.
These problems are separable in the sense that
the unknown  parameters or variables can be separated into two sets: the linear
parameters and the nonlinear parameters.
Problems of this kind often involve a model function that is a linear combination
of parameterized nonlinear basis functions.
The basic idea of VarPro is to treat the linear parameters as dependent on
the nonlinear parameters, and then eliminate the linear parameters from the
problem by using the linear least squares method.
This gives rise to a reduced, but generally more complicated, nonlinear least squares
problem that involves only the nonlinear parameters~\cite{GolubP2003}.
One can then solve the reduced problem for the nonlinear parameters
by a nonlinear least squares method,
typically involving a Gauss-Newton type algorithm coupled with
trust region or backtracking line search strategies~\cite{DennisS1996,Bjorck2015}.
Upon attaining the nonlinear parameters, one then computes the linear parameters
by the linear least squares method.
Although the reduced problem is in general more complicated, the
benefits of variable projection are typically very significant.
These include the reduced dimension of the parameter space,
better conditioning, and faster convergence with the reduced
problem~\cite{RuheW1980,SjobergV1997,GolubP2003}.
In some sense the idea of variable projection to least squares problems can be
analogized to the Schur complement in linear algebra or
the static condensation in computational mechanics (see e.g.~\cite{KarniadakisS2005}).


The VarPro algorithm was originally developed in~\cite{GolubP1973}, and
has been improved and generalized by a number of researchers and applied
to many areas in the past few
decades~\cite{Kaufman1975,RuheW1980,GolubP2003,ChungHN2006,Osborne2007,MullenS2009,OlearyR2013,AskhamK2018,ChenGCL2019,SongXHZ2020,ErichsonZMBKA2020,LeeuwenA2021,NewmanCCR2021}.
In~\cite{GolubP1973} the authors have proved the equivalence between the solution
of the VarPro reduced formulation and that of the original problem, and
developed differentiation formulas for the orthogonal projectors and
the Moore-Penrose pseudoinverses, which are critical to the computation of
the Jacobian matrix in the nonlinear least squares solution of
the reduced problem.
An important simplification to the VarPro algorithm is suggested
in~\cite{Kaufman1975}, which involves computing an approximate Jacobian
rather than the true Jacobian.
This significantly reduces the
per-iteration cost of VarPro, with generally insignificant or negligible
sacrifice to the accuracy for many problems~\cite{GanCCC2018}.
The variable projection algorithms for problems with constraints
on the linear or nonlinear parameters are investigated
in e.g.~\cite{KaufmanP1978,SimaH2007,OlearyR2013,CornelioPN2014}, among
others. The implementations of the VarPro method have been
discussed in~\cite{Krogh1974,OlearyR2013}.
In~\cite{GolubP2003} the original developers of VarPro have reviewed
the developments of this method up to the early 2000s and compiled an extensive list of 
areas for its ongoing and potential applications.
A generalization of the variable projection approach has been considered
in~\cite{RuheW1980}, which deals with two separate classes of
variables without requiring one class to be linear; see also
more recent contributions on the generalization of VarPro
in e.g.~\cite{AravkinL2012,ShearerG2013,HerringNR2018,LeeuwenA2021}.
We would also like to mention the simplification of the Jacobian
matrix in~\cite{RuanoJF1991}, and the algorithm of~\cite{McLooneBI1998},
which resembles the variable projection approach in some sense;
see a comparison of these algorithms with the variable projection
method in~\cite{GanCCC2018}.
An approach related to variable projection is the so-called
block coordinate descent~\cite{NocedalW1999}, which alternates
between the minimization of two separate sets of variables
involved in the problem~\cite{RuheW1980,ChungHN2006,CyrGPPT2020}.


The VarPro algorithm or its variants for training neural networks
have been the subject of several studies in
the literature~\cite{WeiglB1993,WeiglGB1993,WeiglB1994,SjobergV1997,PereyraSW2006,KimL2008,NewmanRHW2020,NewmanCCR2021}.
The projection learning algorithm developed in~\cite{WeiglB1993,WeiglGB1993,WeiglB1994}
is in the same spirit as variable projection, and it computes
the linear parameters by the linear least squares method and
the nonlinear parameters by a gradient descent scheme.
In~\cite{SjobergV1997} the authors have proved that the reduced nonlinear
functional of the variable projection approach,
while seemingly more complicated, leads to a better-conditioned
problem and always converges faster than the original problem; see
also~\cite{RuheW1980}. The VarPro method together with the Levenberg-Marquardt
algorithm is employed for the training of two-layered neural networks
in~\cite{PereyraSW2006,KimL2008} and compared with other related approaches.
In the recent works~\cite{NewmanRHW2020,NewmanCCR2021} the authors extend
the variable projection approach to deal with non-quadratic objective functions,
such as the cross-entropy function in classification tasks, and also
present a stochastic optimization method (termed ``slimTrain'')
based on variable projection
for training deep neural networks with attractive properties.


In the current work we focus on the variable projection approach 
for solving partial differential equations. We numerically approximate
the solution fields to linear and nonlinear PDEs by exploiting 
variable projection together with artificial neural networks.
For computational PDEs, the issues one would encounter with VarPro
are a little different from those for data fitting problems or
function approximations, which account for the majority of
applications the VarPro algorithm is developed for in the literature.
%
For solving PDEs, we do not have the data for the field function to be solved for,
unlike in data fitting problems.
What we do have are the conditions (or constraints) the solution field needs
to satisfy, namely, the PDEs, the boundary conditions, and also the initial conditions
if the problem is time-dependent.
In order to deal with this type of problems, the variable projection method
needs to be adapted accordingly.

The general approach with VarPro and artificial neural networks for
solving PDEs is as follows.
We employ a feed-forward neural network with one or more hidden layers to
represent the field solution to the PDE, requiring that the output layer
be linear (i.e.~applying no activation function) and with zero bias.
We enforce the PDEs on a set of collocation points in the domain, and
enforce the boundary/initial conditions on a set of collocation
points on the appropriate boundaries of the spatial (or spatial-temporal) domain.
This gives rise to a set of discrete equations about the field function
to be solved for, which depends on the weight/bias coefficients in
the output/hidden layers of the neural network.
In turn, this set of equations leads to a nonlinear least squares problem about
the neural-network coefficients, providing an opportunity for the variable
projection method if this nonlinear least squares problem is separable.

It is necessary to distinguish two types of linearities (or nonlinearities)
before the variable projection approach can be used to solve 
the above nonlinear least squares problem.
The first type concerns whether the neural-network coefficients
are linear (or nonlinear) with respect to the output field of the network.
Since no activation function is applied to the output layer,
the output-layer coefficients are linear and the hidden-layer coefficients
are nonlinear with respect to the network output.
The second type concerns whether the boundary/initial value problem
with the given PDE is linear (or nonlinear)
with respect to the field function to be solved for.

If the boundary/initial problem is linear,
i.e.~both the PDE and the boundary/initial conditions are linear with respect to
the solution field, then the aforementioned nonlinear least squares problem
is separable. The output-layer coefficients of the neural network
are the linear parameters and the hidden-layer coefficients are
the nonlinear parameters in this separable nonlinear least squares problem.
In this case, employing the VarPro approach for training the neural network
to solve the given boundary/initial value problem
would be conceptually straightforward.

On the other hand, if the boundary/initial value problem is nonlinear,
i.e.~either the PDE itself or the associated boundary/initial conditions are
nonlinear with respect to the solution field, the aforementioned nonlinear
least squares problem is not separable. In this case all the weight/bias
coefficients in the neural network become nonlinear parameters in
the aforementioned nonlinear least squares problem.
Therefore, the variable projection approach cannot be directly used
for solving nonlinear PDEs (or problems with nonlinear boundary/initial conditions).
How to enable the variable projection method to solve nonlinear PDEs
is the focus of the current work.

In this paper we present a Newton-variable projection method together with
artificial neural networks for solving
nonlinear boundary/initial value problems
(nonlinear PDEs or nonlinear boundary/initial conditions).
Given a nonlinear boundary/initial value problem, we first linearize
the problem for the Newton iteration, with a particular linearized
form. More specifically, the linearization is formulated in terms of
the updated approximation field, not the increment field.
This linearization form is critical to the accuracy of the
current Newton-VarPro method. The linearized system (PDE and boundary/initial conditions)
is linear with respect to the updated approximation field,
and it is solved by the variable projection approach together with
the neural networks.
Therefore, to solve nonlinear PDEs, the current method
involves an overall Newton iteration.
Within each iteration, we use the VarPro method together
with ANNs to solve the linearized
system to attain the updated field approximation.
Upon convergence of the Newton iteration, the weight/bias coefficients of
the neural network contain the representation of the solution
field to the original nonlinear problem.

The VarPro method together with ANNs
for solving linear PDEs has been considered first in this paper.
We discuss in some detail how to
implement the Jacobian matrix together with neural networks, and how to introduce
perturbations in VarPro when solving the reduced problem in order to
prevent the solution from being trapped to the local minima in
the nonlinear least squares computation.
We have presented a number of numerical examples, involving both linear and
nonlinear PDEs, to test the performance of the VarPro method developed here.
We observe that, for smooth field solutions,
the VarPro errors decrease exponentially as the number of
collocation points or the number of output-layer coefficients
in the neural network increases, which is reminiscent of the spectral convergence of
traditional high-order methods~\cite{KarniadakisS2005,SzaboB1991,ZhengD2011,DongS2012,Dong2018,Dong2015clesobc,LinYD2019,YangD2019,YangD2020}.
We also compare extensively the performance of the current VarPro method with
that of the ELM method from~\cite{DongL2021,DongY2021}.
The numerical results show that, under identical conditions and
network configurations, the VarPro method is considerably more accurate
than the ELM method, especially when the size of the neural network is small. 
On the other hand, the computational cost (i.e.~network training time)
of the VarPro method is usually much higher than that of the ELM method.

In the current work the VarPro method and the neural networks are implemented
based on the Tensorflow (www.tensorflow.org) and Keras (keras.io) libraries
in Python. The scipy and numpy libraries in Python are used for the linear and nonlinear
least squares computations. All the numerical tests are carried out on
a MAC computer in the authors' institution.


The main contribution of this paper lies in the Newton-VarPro method together with
artificial neural networks for solving nonlinear partial differential
equations. To the best of the authors' knowledge, this work seems also to be the first
time when the variable projection approach (with ANNs) is
extended and adapted  to solving
linear partial differential equations.


The rest of this paper is structured as follows.
In Section \ref{sec:method} we first outline how to solve linear PDEs with 
the variable projection approach together with ANNs.
The computations for the reduced residual function and
the Jacobian matrix of the reduced problem, and the VarPro algorithm with perturbations
are discussed in detail. Then we introduce the Newton-VarPro method
together with ANNs for solving nonlinear PDEs.
In Section \ref{sec:tests} we present several numerical examples with
linear and nonlinear PDEs to demonstrate the accuracy of the VarPro method
developed herein. The performance of the current VarPro method
is compared extensively with that of the ELM method from~\cite{DongL2021,DongY2021}.
Section~\ref{sec:summary} then concludes the presentation with
some closing remarks and comments on the presented method.

\section{Variable Projection with Artificial Neural Networks
for Computational PDEs}
\label{sec:method}

We develop an algorithm combining the variable projection (VarPro) framework
with artificial neural networks (ANN) for numerically approximating PDEs.
For linear PDEs, the ANN representation of the solution field leads to
a separable nonlinear least squares problem, which can be
solved by the variable projection approach.
For nonlinear PDEs, on the other hand, the 
ANN representation of the solution field
leads to a nonlinear least squares (NLLSQ) problem that is
not separable, preventing the use of the variable projection strategy.
We overcome this issue by a combined Newton-variable projection method, which
enables the variable projection approach in solving nonlinear PDEs.
In the following subsections we first illustrate the
VarPro/ANN algorithm for solving linear PDEs, and then introduce the Newton-VarPro/ANN
algorithm for solving nonlinear PDEs.

\subsection{Variable Projection Method for Solving Linear PDEs}
\label{sec:lin}

Consider a domain $\Omega\subset \mbs R^d$ ($d=1$ to $3$)
and the following linear boundary-value problem on $\Omega$,
\begin{subequations}\label{eq_1}
  \begin{align}
  &
  Lu = f(\mbs x),  \label{eq_1a} \\
  &
  Bu = g(\mbs x), \quad\text{on}\ \partial\Omega. \label{eq_1b}
  \end{align}
\end{subequations}
In these equations $\mbs x=(x_1,\dots,x_d)$ denotes the coordinate,
$u(\mbs x)$ is the field solution to be solved for,
$L$ denotes a linear differential
operator, $B$ denotes a linear algebraic or differential operator 
on the boundary $\partial\Omega$ representing the boundary conditions,
and $f(\mbs x)$ and $g(\mbs x)$
are prescribed non-homogeneous terms in the domain or on
the boundary.
We assume that $L$ may include linear differential operators
with respect to the time $t$ (e.g.~$\frac{\partial}{\partial t}$,
$\frac{\partial^2}{\partial t^2}$). In such a case, this becomes
an initial boundary-value problem, and we treat the time $t$
in the same way as the spatial coordinates. We designate
the last coordinate $x_d$ as $t$, and $\Omega$ becomes a
spatial-temporal domain. 
Accordingly, we assume that the
boundary condition~\eqref{eq_1b} in this case should include appropriate
initial condition(s) with respect to $t$, which will be imposed
only on the portion of $\partial\Omega$ corresponding to the initial condition(s). 
The point here is that the equations~\eqref{eq_1a}--\eqref{eq_1b} may denote
a time-dependent problem, and we will not distinguish the
stationary and time-dependent cases in the following discussions.
We assume that the problem~\eqref{eq_1} is well-posed.


We  approximate the solution field $u(\mbs x)$
by a feed-forward neural network~\cite{GoodfellowBC2016} with
$(L+1)$ layers, where $L$ is an integer satisfying $L\geqslant 2$.
The input layer (layer $0$) of the neural network contains $d$ nodes which
represent the coordinate $\mbs x$, and the output layer (layer $L$)
contains $1$ node which
represents the solution $u$.
The $(L-1)$ layers in between are the hidden layers.
From layer to layer the network logic
represents an affine transform followed by a node-wise
function composition with an activation function
$\sigma(\cdot)$~\cite{GoodfellowBC2016}.
The coefficients of the affine transforms
are referred to as the weight and bias coefficients of the neural network.
For the convenience of presentation, we
use the vector $[M_0, M_1,\dots, M_L]$ to denote the architecture
of the neural network, where $M_i$ ($0\leqslant i\leqslant L$)
denotes the number of nodes in layer $i$, with $M_0=d$
and $M_L=1$. We also use $M=M_{L-1}$ to denote the number of nodes
in the last hidden layer in what follows.
The weight/bias coefficients in all the hidden layers and
in the output layer are the trainable parameters of the neural network.

In the current work, we make the assumption that
the output layer contains no bias (or zero bias), and no
activation function (or equivalently it uses the identity activation function
$\sigma(x)=x$). So the output layer of the neural network
is linear in this paper.

Let $\Phi_j(\bm\theta,\mbs x)$ ($1\leqslant j\leqslant M$)
denote the output fields of the last hidden layer,
where $\bm\theta=(\theta_1,\dots,\theta_{N_h})^T$ 
denotes the vector of weight/bias coefficients
in all  the hidden layers of the network, with $N_h=\sum_{i=1}^{L-1}M_i(M_{i-1}+1)$.
Then we have the following expansion relation,
\begin{equation}\label{eq_2}
  u(\mbs x) = \sum_{j=1}^M \beta_j\Phi_j(\bm\theta, \mbs x) = \bm\Phi(\bm\theta,\mbs x)\bm\beta
\end{equation}
where
$\bm\Phi(\bm\theta,\mbs x)=[\Phi_1(\bm\theta,\mbs x), \dots,\Phi_M(\bm\theta,\mbs x)]$
denotes the set of output fields of the last hidden layer, and
$\bm\beta=[\beta_1,\dots,\beta_M]^T$ is the vector of
weight coefficients of the output layer.
Note that $(\bm\theta,\bm\beta)$ are the trainable parameters of
the neural network. Note also that $M$ represents the number of nodes
in the last hidden layer, as well as the number of output-layer coefficients.


We choose a set of $N$ ($N\geqslant 1$) collocation points
on $\Omega$, which can be chosen according to a certain distribution
(e.g.~random, uniform). Among them $N_b$ ($1\leqslant N_b\leqslant N-1$)
collocation points reside on the boundary $\partial\Omega$,
and the rest of the points are from the interior of $\Omega$.
We use $\mathbb{X}$ to denote the set of all the collocation points and
$\mathbb{X}_b$ to denote the set of collocation points on $\partial\Omega$.
In the current paper, for simplicity we assume that
$\Omega$ is a rectangular domain,
given by the interval $[a_i,b_i]$ ($1\leqslant i\leqslant d$)
in the $i$-th direction. We employ a uniform set of grid
points (including the boundary end points) in each direction
as the collocation points 
for the numerical tests in Section \ref{sec:tests}.


The input training data to the neural network consist of
the coordinates of the all the $N$ collocation points on $\Omega$.
We use the $N\times d$ matrix $\mbs X$ to denote
the input data. Each row of $\mbs X$ denotes the coordinates
a collocation point. Let the $N\times 1$ matrix $\mbs U$ (column vector)
denote the output data of the neural network, which represents
the solution field $u(\mbs x)$ evaluated on all the $N$
collocation points. We use the $N\times M$ matrix $\bm\Psi$
to denote the output data of the last hidden layer of the neural
network, which represents the output fields $\bm\Phi(\bm\theta,\mbs x)$
of the last hidden layer evaluated on all
the $N$ collocation points.


Inserting the expansion \eqref{eq_2} into \eqref{eq_1}, and
enforcing the equation \eqref{eq_1a} on all the collocation
points from $\mathbb{X}$ and the equation \eqref{eq_1b}
on all the boundary collocation points from $\mathbb{X}_b$,
we arrive at the following system,
\begin{subequations}\label{eq_3}
  \begin{align}
    &
    \sum_{j=1}^M\left[L\Phi_j(\bm\theta,\mbs x_p) \right]\beta_j = f(\mbs x_p),
    \quad 1\leqslant p\leqslant N, \ \text{where}\ \mbs x_p\in \mathbb{X}, 
    \label{eq_3a} \\
    &
    \sum_{j=1}^M \left[B\Phi_j(\bm\theta,\mbs x_q) \right]\beta_j = g(\mbs x_q),
    \quad 1\leqslant q\leqslant N_b,\ \text{where}\ \mbs x_q\in \mathbb{X}_b. \label{eq_3b}
  \end{align}
\end{subequations}
This is a system of $(N+N_b)$ algebraic equations about the trainable
parameters $(\bm\theta,\bm\beta)$, with $(N_h+M)$ unknowns. Note that
for a given $\bm\theta$ the terms
$L\Phi_j(\bm\theta,\mbs x_p)$ and $B\Phi_j(\bm\theta,\mbs x_q)$
in the above equations can be computed by forward evaluations of
the neural network and auto-differentiations.

We seek a least squares solution for $(\bm\theta,\bm\beta)$ to the system \eqref{eq_3}.
This system is linear with respect to $\bm\beta$, and nonlinear with
respect to $\bm\theta$. This leads to a separable nonlinear least squares problem.
Therefore we adopt the variable projection approach~\cite{GolubP1973}
for the least squares solution of the system~\eqref{eq_3}.

To make the formulation more compact, we re-write the system~\eqref{eq_3}
into a matrix form,
\begin{equation}\label{eq_4}
  \mbs H(\bm\theta)\bm\beta = \mbs S,
  \ \ 
\text{where}\ \
  \mbs H(\bm\theta)=\begin{bmatrix}
  \vdots \\
  L\bm\Phi(\bm\theta,\mbs x_p)\\
  \vdots  \\ \hdashline[1pt/1pt]
  \vdots \\
  B\bm\Phi(\bm\theta,\mbs x_q)\\
  \vdots
  \end{bmatrix}_{(N+N_b)\times M},
  \ \ 
  \mbs S = \begin{bmatrix}
  \vdots \\
  f(\mbs x_p)\\
  \vdots  \\ \hdashline[1pt/1pt]
  \vdots \\
  g(\mbs x_q)\\
  \vdots
  \end{bmatrix}_{(N+N_b)\times 1}.
\end{equation}
For any given $\bm\theta$, the least squares solution for
the linear parameters $\bm\beta$ to this system is given by
\begin{equation}\label{eq_7}
  \bm\beta = \left[\mbs H(\bm\theta)\right]^+\mbs S,
\end{equation}
where the superscript in $\mbs H^{+}$ denotes the
Moore-Penrose pseudo-inverse of $\mbs H$.
Define the residual function of the system~\eqref{eq_4} by
\begin{equation}\label{eq_8}
  \mbs r(\bm\theta) = \mbs H(\bm\theta)\bm\beta - \mbs S
  = \mbs H(\bm\theta) \mbs H^+(\bm\theta)\mbs S - \mbs S,
\end{equation}
where the linear parameter $\bm\beta$ has been eliminated
by using equation~\eqref{eq_7}.
We compute the optimal nonlinear parameters $\bm\theta_{opt}$ by minimizing
the Euclidean norm of the residual function $\mbs r$,
\begin{equation}\label{eq_9}
  \bm\theta_{opt}= \argmin_{\bm\theta} \frac12\|\mbs r(\bm\theta) \|^2
  =\argmin_{\bm\theta} \frac12\| \mbs H(\bm\theta) \mbs H^+(\bm\theta)\mbs S - \mbs S  \|^2,
\end{equation}
where $\|\cdot\|$ denotes the Euclidean norm.
After $\bm\theta_{opt}$ is obtained, we can compute the optimal linear
parameters $\bm\beta_{opt}$ based on equation \eqref{eq_7} or
by solving equation~\eqref{eq_4} using the linear least squares method.
Outlined above is the essence of the variable projection
approach for solving the system~\eqref{eq_3} for $(\bm\theta,\bm\beta)$.


The problem represented by~\eqref{eq_9} is a nonlinear
least squares problem about $\bm\theta$ only, where the linear parameter
$\bm\beta$ has been eliminated.
We solve this problem by a Gauss-Newton algorithm combined with
a trust region strategy. Specifically, in the current paper
we solve this problem by
employing the nonlinear least squares library routine
``scipy.optimize.least\_squares'' from the scipy package in Python,
which implements the Gauss-Newton method together with
a trust region reflective algorithm~\cite{BranchCL1999,ByrdSS1988}.

\begin{algorithm}[tb]\label{alg_1}
  \DontPrintSemicolon
  \SetKwInOut{Input}{input}\SetKwInOut{Output}{output}
  \Input{$\bm\theta$; input data $\mbs X$ to neural network;
    source data $\mbs S$.
  }
  \Output{$\mbs r(\bm\theta)$.}
  \BlankLine\BlankLine
  update the hidden-layer coefficients of the neural network by $\bm\theta$\;
  \eIf{$\bm\theta=\bm\theta_{s}$}{
    retrieve $\mbs H(\bm\theta_s)$, and set $\mbs H(\bm\theta)=\mbs H(\bm\theta_s)$\;
    retrieve $\bm\beta^{LS}(\bm\theta_s)$, and
    set $\bm \beta^{LS}(\bm\theta)=\bm \beta^{LS}(\bm\theta_s)$\;
  }{
    compute $\mbs H(\bm\theta)$ using the input data $\mbs X$\;
    solve equation~\eqref{eq_4} by the linear least squares method to get $\bm\beta^{LS}(\bm\theta)$\;
    set $\bm\theta_s=\bm\theta$, and save $\mbs H(\bm\theta)$ and $\bm\beta^{LS}(\bm\theta)$\;
  }
  \BlankLine
  
  compute $\mbs r(\bm\theta)$ by equation~\eqref{eq_10}\;
  
  \caption{Computing the residual $\mbs r(\bm\theta)$}
\end{algorithm}

The scipy routine ``least\_squares()'' requires
two functions as input, which are needed by the Gauss-Newton algorithm.
These are, for any given $\bm\theta$,
\begin{itemize}
\item a function for computing the residual $\mbs r(\bm\theta)$, and
\item a function for computing the Jacobian
  matrix $\frac{\partial\mbs r}{\partial\bm\theta}$.
\end{itemize}
The computation for $\mbs r(\bm\theta)$ is straightforward.
For a given $\bm\theta$, we first solve equation~\eqref{eq_4}
by the linear least squares method for the minimum-norm least
squares solution $\bm\beta^{LS}$. Then we compute the residual
according to equation~\eqref{eq_8} as follows,
\begin{equation}\label{eq_10}
  \mbs r(\bm\theta) = \mbs H(\bm\theta)\bm\beta^{LS} - \mbs S.
\end{equation}
Note that the Moore-Penrose inverse $\mbs H^+(\bm\theta)$ is
not explicitly computed in the implementation.
In the current paper we employ the linear least squares routine
``scipy.linalg.lstsq'' from scipy to solve \eqref{eq_4}
for $\bm\beta^{LS}$.
The computation for $\mbs r(\bm\theta)$ is summarized in
the Algorithm~\ref{alg_1}.

\begin{remark}\label{rem_1}
  Let us elaborate on, for a given $\bm\theta$ and the input data $\mbs X$,
  how to compute the matrix $\mbs H(\bm\theta)$ on line $6$
  of Algorithm~\ref{alg_1}.
  As defined in~\eqref{eq_4}, $H(\bm\theta)$ consists of the terms
  $L\bm\Phi(\bm\theta,\mbs x_p)$ ($\mbs x_p\in\mathbb{X}$) and
  $B\bm\Phi(\bm\theta,\mbs x_q)$ ($\mbs x_q\in\mathbb{X}_b$).
  These terms involve the output fields of the last hidden layer
  $\bm\Phi(\bm\theta,\mbs x)$, and their derivatives up to a certain
  order, evaluated on all the collocation points.
  All these terms can be computed by evaluating the neural network
  on the input data $\mbs X$ and by auto-differentiations.
  Specifically, in our implementation we have
  created a sub-model to the neural network in Keras,
  with the neural network's input as its input and
  with the output of the neural network's last hidden layer
  as the sub-model's output. Let us refer to this sub-model
  as the last-hidden-layer-model. Let $m$ denote the order
  of the PDE~\eqref{eq_1a}, and we assume that the
  hidden-layer coefficients have been updated
  by the given $\bm\theta$. Then computing $\mbs H(\bm\theta)$
  involves the the following procedure:
  \begin{enumerate}[(i)]
  \item evaluate the last-hidden-layer-model on the input $\mbs X$
    to get $\bm\Phi(\bm\theta,\mbs x)$ on all the collocation points;
  \item compute the derivatives of $\bm\Phi(\bm\theta,\mbs x)$ with respect to
    $\mbs x$, up to the order $m$, on all
    the collocation points by a forward-mode auto-differentiation;
  \item compute $L\bm\Phi(\bm\theta,\mbs x)$ on all the collocation points
    based on the data for $\bm\Phi(\bm\theta,\mbs x)$ and its derivatives;
  \item extract the boundary data (i.e.~on the boundary collocation points)
    for $\bm\Phi(\bm\theta,\mbs x)$ and its derivatives from those data
    attained from steps (i) and (ii);
  \item compute $B\bm\Phi(\bm\theta,\mbs x)$ based on the boundary data
    for $\bm\Phi(\bm\theta,\mbs x)$ and its derivatives;
  \item assemble $L\bm\Phi(\bm\theta,\mbs x)$ (on all the collocation points)
    and $B\bm\Phi(\bm\theta,\mbs x)$ (on the boundary collocation points)
    to form $\mbs H(\bm\theta)$ based on equation~\eqref{eq_4}.
  \end{enumerate}
  Note that we employ the forward-mode auto-differentiations to compute the derivatives
  of $\bm\Phi(\bm\theta,\mbs x)$ in step (ii) above, because the number of nodes
  in the last hidden layer ($M$) is typically much larger than that
  in the input layer ($d$). In this case the forward-mode auto-differentiation
  is significantly faster than the reverse-mode auto-differentiation.
  In our implementation we have used
  the ``ForwardAccumulator'' from the Tensorflow library for the forward-mode
  auto-differentiations.

\end{remark}


For computing the Jacobian matrix $\frac{\partial\mbs r}{\partial\bm\theta}$
we consider the following formula, which is due to~\cite{GolubP1973},
\begin{equation}\label{eq_11}
  \begin{split}
  \frac{\partial\mbs r}{\partial\bm\theta}=&
  \left[\mbs I - \mbs H(\bm\theta)\mbs H^+(\bm\theta) \right]
  \frac{\partial\mbs H}{\partial\bm\theta}\mbs H^+(\bm\theta)\mbs S \
  + \
  \left[\mbs H^+(\bm\theta) \right]^T\frac{\partial\mbs H^T}{\partial\bm\theta}
  \left[\mbs I - \mbs H(\bm\theta)\mbs H^+(\bm\theta) \right]\mbs S \\
  \approx &
  \left[\mbs I - \mbs H(\bm\theta)\mbs H^+(\bm\theta) \right]
  \frac{\partial\mbs H}{\partial\bm\theta}\mbs H^+(\bm\theta)\mbs S,
  \end{split}
\end{equation}
where $\mbs I$ denotes the identity matrix and equation~\eqref{eq_8} has been
used. Note that here
we have adopted the simplification suggested by~\cite{Kaufman1975} to keep only
the first term for an approximation of $\frac{\partial\mbs r}{\partial\bm\theta}$.
So the Jacobian matrix is computed only approximately. This
greatly simplifies the computation, and as observed in~\cite{Kaufman1975,GolubP2003}
only slightly or  moderately increases the number of Gauss-Newton iterations.

In light of~\eqref{eq_11}, we compute the approximate Jacobian matrix
as follows. For any given $\bm\theta$, note that
\begin{equation}\label{eq_12}
  \mbs J_0(\bm\theta) \equiv
  \frac{\partial\mbs H}{\partial\bm\theta}\mbs H^+(\bm\theta)\mbs S
  = \frac{\partial\mbs H}{\partial\bm\theta}\bm\beta^{LS}
  = \frac{\partial\mbs V}{\partial\bm\theta},
\end{equation}
where $\bm\beta^{LS}$ is the least squares solution of~\eqref{eq_4},
and
\begin{equation}\label{eq_12a}
  \mbs V=\mbs H(\bm\theta)\bm\beta_c^{LS}, \quad \text{with}\
  \bm\beta_c^{LS} = \left.\bm\beta^{LS}\right|_{\bm\theta}.
\end{equation}
Here $\bm\beta_c^{LS}$ is a {\em constant}
vector that equals $\bm\beta^{LS}$ at the given $\bm\theta$.
The vector $\mbs V(\bm\theta)$ of length $(N+N_b)$
represents the field
$\begin{bmatrix}
  Lu(\mbs x) \\
  Bu(\mbs x)
  \end{bmatrix}
$
evaluated on the collocation points (and the boundary collocation points),
with $\bm\theta$ as the hidden-layer coefficients and
$\bm\beta_c^{LS}$ as the output-layer coefficients in the neural network.
We would like to emphasize that $\bm\beta_c^{LS}$ is considered to be constant
and does not depend on $\bm\theta$ when computing
$\mbs J_0(\bm\theta)=\frac{\partial\mbs V}{\partial\bm\theta}$.
For a given $\bm\theta$, $\mbs J_0(\bm\theta)$ can be
computed by an auto-differentiation of the neural network.

In light of~\eqref{eq_12} we transform~\eqref{eq_11} into
\begin{equation}\label{eq_13}
  \frac{\partial\mbs r}{\partial\bm\theta} = \mbs J_0(\bm\theta)
  - \mbs H(\bm\theta)\mbs H^+(\bm\theta)\mbs J_0(\bm\theta)
  = \mbs J_0(\bm\theta) - \mbs J_1(\bm\theta).
\end{equation}
The term $\mbs J_1(\bm\theta) = \mbs H(\bm\theta)\mbs H^+(\bm\theta)\mbs J_0(\bm\theta)$
can be computed as follows.
For any given $\bm\theta$, we first solve the following system for the
matrix $\mbs K(\bm\theta)$ by the linear least squares method,
\begin{equation}\label{eq_14}
  \mbs H(\bm\theta)\mbs K(\bm\theta) = \mbs J_0(\bm\theta).
\end{equation}
Then we compute $\mbs J_1(\bm\theta)$ by
\begin{equation}\label{eq_15}
  \mbs J_1(\bm\theta) = \mbs H(\bm\theta)\mbs K(\bm\theta).
\end{equation}

\begin{algorithm}[tb]\label{alg_2}
  \DontPrintSemicolon
  \SetKwInOut{Input}{input}\SetKwInOut{Output}{output}
  \Input{$\bm\theta$; input data $\mbs X$ to neural network;
    source data $\mbs S$.
  }
  \Output{$\frac{\partial\mbs r}{\partial\bm\theta}$.}
  \BlankLine\BlankLine
  update the hidden-layer coefficients of the neural network by $\bm\theta$\;
  \eIf{$\bm\theta=\bm\theta_{s}$}{
    retrieve $\mbs H(\bm\theta_s)$, and set $\mbs H(\bm\theta)=\mbs H(\bm\theta_s)$\;
    retrieve $\bm\beta^{LS}(\bm\theta_s)$, and
    set $\bm \beta^{LS}(\bm\theta)=\bm \beta^{LS}(\bm\theta_s)$\;
  }{
    compute $\mbs H(\bm\theta)$ using the input data $\mbs X$\;
    solve equation~\eqref{eq_4} by the linear least squares method to get $\bm\beta^{LS}(\bm\theta)$\;
    set $\bm\theta_s=\bm\theta$, and save $\mbs H(\bm\theta)$ and $\bm\beta^{LS}(\bm\theta)$\;
  }
  \BlankLine
  
  compute $\mbs V(\bm\theta)$ by equation \eqref{eq_12a}\;
  compute $\mbs J_0(\bm\theta)$ based on equation~\eqref{eq_12} by auto-differentiations\;
  solve equation~\eqref{eq_14} for $\mbs K(\bm\theta)$ by the linear least squares method\;
  compute $\mbs J_1(\bm\theta)$ by equation~\eqref{eq_15}\;
  compute $\frac{\partial\mbs r}{\partial\bm\theta}$ by equation~\eqref{eq_13}\;
  
  \caption{Computing the Jacobian matrix $\frac{\partial\mbs r}{\partial\bm\theta}$}
\end{algorithm}

Therefore, in order to compute the Jacobian matrix we first solve equation~\eqref{eq_4}
for $\bm\beta^{LS}$ by the linear least squares method, and then use~\eqref{eq_12}
to compute $\mbs J_0(\bm\theta)$. We then compute $\mbs J_1(\bm\theta)$ by
equations~\eqref{eq_14} and~\eqref{eq_15}. Finally
the Jacobian matrix $\frac{\partial\mbs r}{\partial\bm\theta}$ is computed
by equation \eqref{eq_13}.
These computations involve only the linear least squares method and
the auto-differentiations of the neural network.
The computation for the Jacobian matrix is summarized
in Algorithm~\ref{alg_2}.


\begin{remark}\label{rem_2}
  Let us elaborate on how to compute the matrix $\mbs J_0(\bm\theta)$,
  which has a dimension $(N+N_b)\times N_h$ ($N_h$ denoting the
  total number of hidden-layer coefficients), on the lines $10$ and $11$
  in Algorithm~\ref{alg_2}.
  This is for a given $\bm\theta$, $\bm\beta$ ($\bm\beta=\bm\beta_c^{LS}$),
  and the input data $\mbs X$ to the neural network.
  Based on equation~\eqref{eq_12a}, the column vector $\mbs V(\bm\theta)$
  consists of the terms $Lu(\mbs x_p)$ ($\mbs x_p\in\mathbb{X}$) and
  $Bu(\mbs x_q)$ ($\mbs x_q\in\mathbb{X}_b$), where $u(\mbs x)$
  is the output field of the neural network obtained with
  the given $(\bm\theta,\bm\beta_c^{LS})$ as the hidden-layer coefficients
  and the output-layer coefficients, respectively.
  It should be noted that $Lu$ and $Bu$ involve the derivatives of $u(\mbs x)$
  with respect to $\mbs x$ (not $\bm\theta$).
  Based on equation~\eqref{eq_12}, the matrix $\mbs J_0(\bm\theta)$
  consists of the terms
  $\left.\frac{\partial (Lu)}{\partial\bm\theta}\right|_{(\bm\theta,\mbs x_p)}$
  ($\mbs x_p\in\mathbb{X}$) and
  $\left.\frac{\partial (Bu)}{\partial\bm\theta}\right|_{(\bm\theta,\mbs x_q)}$
  ($\mbs x_q\in\mathbb{X}_b$).
  These terms can be computed by evaluating the neural network on
  the input data $\mbs X$ and by auto-differentiations with respect to
  $\mbs x$ and $\bm\theta$. We assume again that the PDE~\eqref{eq_1a} is
  of the $m$-th order.
  Given ($\bm\theta,\bm\beta,\mbs X$), we compute $\mbs J_0(\bm\theta)$
  specifically by the following procedure:
  \begin{enumerate}[(i)]
  \item update the hidden-layer coefficients of the neural network by $\bm\theta$,
    and update the output-layer coefficients by $\bm\beta$;
  \item evaluate the neural network on the input $\mbs X$ to obtain
    the output field $u(\mbs x)$ on all the collocation points;
  \item compute the derivatives of $u(\mbs x)$ with respect to $\mbs x$,
    up to the order $m$, by a reverse-mode auto-differentiation;
  \item compute the derivative, with respect to the hidden-layer coefficients,
    for $u(\mbs x)$ and for its derivatives with respect to $\mbs x$
    from steps (ii) and (iii), on all the collocation points
    by a reverse-mode auto-differentiation;
  \item compute $\frac{\partial (Lu)}{\partial\bm\theta}$ on all the collocation
    points based on the data for $u(\mbs x)$ and its derivatives
    from the previous step;
  \item extract the boundary data (i.e.~on the boundary collocation points)
    for $u(\mbs x)$ and its derivatives from the data  obtained from step (iv);
  \item compute $\frac{\partial (Bu)}{\partial\bm\theta}$ based on
    the boundary data for $u(\mbs x)$ and its derivatives from the previous step;
  \item assemble $\frac{\partial (Lu)}{\partial\bm\theta}$ (on all the collocation
    points) and $\frac{\partial (Bu)}{\partial\bm\theta}$ (on the boundary collocation
    points) to form $\mbs J_0(\bm\theta)$.
  \end{enumerate}
  When computing the derivatives of $u(\mbs x)$
  with respect to $\mbs x$ and with respect to the hidden-layer
  coefficients in the steps (iii) and (iv) above,
  in our implementation we have employed a vectorized map (tf.vectorized\_map) together with
  the gradient tape (tf.GradientTape) in the Tensorflow library to vectorize
  the gradient computations.
  
\end{remark}

\begin{remark}\label{rem_3}
  In Algorithms~\ref{alg_1} and~\ref{alg_2}, we have
  saved the matrix $\mbs H(\bm\theta)$ and the vector $\bm\beta^{LS}$
  when they are computed for a new $\bm\theta$; see the lines $2$ to $9$
  in both algorithms. The goal of this extra storage is to save computations. 
  During the Gauss-Newton iterations, the Algorithm~\ref{alg_2}
  is typically invoked to compute the Jacobian matrix
  for the same $\bm\theta$, following the call to the Algorithm~\ref{alg_1}
  for computing the residual $\mbs r(\bm\theta)$.
  In this case one avoids the re-computation
  of the matrix $\mbs H(\bm\theta)$ and the vector $\bm\beta^{LS}$
  for the same $\bm\theta$.
  
\end{remark}

\begin{algorithm}[tb]\label{alg_3}
  \DontPrintSemicolon
  \SetKwInOut{Input}{input}\SetKwInOut{Output}{output}
  \Input{input data $\mbs X$ to neural network;
    source data $\mbs S$; initial guess $\bm\theta_0$;
    maximum perturbation magnitude $\delta>0$;
    preference probability $p$ ($p\in[0,1]$), with default value $p=0.5$.
  }
  \Output{$\bm\theta$ and $\bm\beta$.}
  \BlankLine\BlankLine
  call scipy.optimize.least\_squares routine to solve~\eqref{eq_9},
  using $\bm\theta_0$ as the
  initial guess, with Algorithms~\ref{alg_1} and~\ref{alg_2} as input arguments\;
  set $\bm\theta\leftarrow$ returned solution, and $c\leftarrow$ returned cost\;
  \BlankLine
  \If{$c$ is above a threshold}{
    set $\delta_{\text{pref}}=\text{None}$\;
    \For{$i\leftarrow 1$ \KwTo maximum-number-of-sub-iterations}{
      generate a uniform random number $\xi\in[0,1]$\;
      \eIf{($\delta_{\text{pref}}$ is not $\text{None}$) and ($\xi<p$) }{
        generate a uniform random number $\delta_1 \in[0, \min(1.1\delta_{\text{pref}},\delta)]$\;
      }{
        generate a uniform random number $\delta_1\in[0,\delta]$\;
      }
      \BlankLine
      generate a uniform random vector $\Delta\bm\theta$ of
      the same shape as $\bm\theta$
      on the interval $[-\delta_1,\delta_1]$\;
      set $\bm\vartheta_0\leftarrow \bm\theta + \Delta\bm\theta$\;
      \BlankLine
      call scipy.optimize.least\_squares routine to solve~\eqref{eq_9},
      using $\bm\vartheta_0$ as the
      initial guess, with Algorithms~\ref{alg_1} and~\ref{alg_2} as input arguments\;
      \If{the returned cost is less than $c$}{
        set $\bm\theta\leftarrow$ returned solution, and $c\leftarrow$ returned cost\;
        set $\delta_{\text{pref}}=\delta_1$\;
      }
      \BlankLine
      \If{$c$ is not above a threshold}{ break\;}
    }
  }
  \BlankLine
  solve equation~\eqref{eq_4} for $\bm\beta$ by the linear least squares method\;
  
  \caption{Variable projection algorithm with perturbations}
\end{algorithm}


To solve the system~\eqref{eq_3a}--\eqref{eq_3b} with
the variable projection approach, we first solve the reduced problem~\eqref{eq_9}
for $\bm\theta$ by the the nonlinear least squares method, and
then we solve the equation~\eqref{eq_4} for $\bm\beta$ by
the linear least squares method.
To make the nonlinear least squares computation for~\eqref{eq_9} more robust
(from being trapped to local minima),
in our implementation we have incorporated a perturbation to
the initial guess and a sub-iteration procedure,
in a way analogous to the NLLSQ-perturb method from~\cite{DongL2021}.
The sub-iteration procedure will be triggered 
if the nonlinear least squares computation fails to converge or
the converged cost value is not small enough.
The overall variable projection algorithm with perturbations for solving
the system~\eqref{eq_3} is summarized in Algorithm~\ref{alg_3}.
The perturbations to the initial guess of the nonlinear least squares
computation are generated on the lines $6$ to $13$ in Algorithm~\ref{alg_3}.

\begin{remark}\label{rem_4a}
  In Algorithm~\ref{alg_3}, when generating the perturbation magnitude $\delta_1$,
  we have incorporated a preferred perturbation magnitude $\delta_{\text{pref}}$
  and a preference probability $p$.
  Here $\delta_{\text{pref}}$ keeps the last perturbation magnitude $\delta_1$
  that has resulted in a reduction in the converged cost.
  The lines $6$ to $11$ of Algorithm~\ref{alg_3} basically means that,
  with a probability $p$,
  we will generate the next perturbation magnitude $\delta_1$
  based on the preferred magnitude $\delta_{\text{pref}}$. Otherwise, we will generate
  the next perturbation magnitude based on the original maximum magnitude $\delta$.
  After the algorithm hits upon a favorable perturbation magnitude,
  the employment of $\delta_{\text{pref}}$ and the probability $p$ tends to promote
  the use this value. When the algorithm is close to convergence,
  this also tends to reduce the amount of the perturbation,
  which is conducive to achieving convergence.
  In the current paper we employ a preference probability $p=0.5$ with the
  variable projection algorithm for all
  the numerical tests in Section~\ref{sec:tests}.

\end{remark}

\begin{remark}\label{rem_4}
  If the problem consisting of equations~\eqref{eq_1a}--\eqref{eq_1b}
  is time dependent, for longer-time or long-time simulations
  we employ the block time marching scheme from~\cite{DongL2021}
  together with the variable projection algorithm developed here.
  The basic idea is as follows. If the domain $\Omega$ has a large dimension in time,
  we first divide the temporal dimension
  into a number of windows (referred to as time blocks),
  so that each time block has a moderate size in time.
  We solve the problem using
  the variable projection algorithm on the spatial-temporal domain of
  each time block individually and successively.
  After one time block is computed, the field solution (and
  also possibly its derivatives) evaluated at the last time instant of this block
  is used as the initial condition(s) for the subsequent time block.
  We refer the reader to~\cite{DongL2021} for more detailed discussions
  of the block time marching scheme.

\end{remark}


\begin{remark}\label{rem_4b}
  It would be interesting to compare the current VarPro method with the extreme
  learning machine (ELM) method from~\cite{DongL2021,DongY2021} for solving PDEs.
  With ELM, the weight/bias coefficients in all the hidden layers
  of the neural network are pre-set to random values and are fixed, while
  the output-layer coefficients are computed by the linear least squares method
  for solving linear PDEs and by the nonlinear least squares method for solving nonlinear
  PDEs~\cite{DongL2021}. In~\cite{DongL2021,DongY2021} the hidden-layer
  coefficients are set and fixed to uniform random values generated on
  the interval $[-R_m,R_m]$, where $R_m$ is a user-provided constant (hyperparameter).
  The constant $R_m$ has an influence on the accuracy of ELM, and the optimal
  $R_m$ value (denoted by $R_{m0}$ in~\cite{DongY2021})
  can be computed by the method from~\cite{DongY2021} based on
  the differential evolution algorithm. It is crucial to note that in ELM
  all the hidden-layer coefficients are fixed (not trained) once they are set.

  With the VarPro method, the output-layer coefficients
  are always computed by the linear least squares method, once the hidden-layer
  coefficients are determined. The weight/bias coefficients in
  the hidden layers are determined by considering the reduced
  problem, which eliminates the linear output-layer
  coefficients. The hidden-layer coefficients are computed by solving the
  reduced problem using the nonlinear least squares method.
  With the VarPro approach, the hidden-layer coefficients of the neural network
  are trained/computed first by solving the reduced problem,
  and then the output-layer coefficients are computed by
  the linear least squares method afterwards. If the maximum number of iterations
  is set to zero in the nonlinear least squares solution of
  the reduced problem, the VarPro algorithm will be reduced to essentially
  the ELM method.

  With the same neural network architecture and under the same settings,
  the VarPro method is in general significantly more accurate than
  the ELM method. In particular, VarPro can produce highly accurate solutions
  when the number of nodes in the last hidden layer is not large.
  In contrast, the result produced by ELM in this case is usually much less accurate
  or utterly inaccurate. VarPro achieves the higher accuracy at the price of
  the computational cost. Because VarPro needs to solve the reduced problem
  for the hidden-layer coefficients by a nonlinear least squares computation,
  its computational cost is usually much higher than that of the ELM method, which only
  computes the output-layer coefficients by the linear least squares method
  (for linear PDEs). In numerical simulations with VarPro, we initialize the hidden-layer
  coefficients (i.e.~the initial guess $\bm\theta_0$ in Algorithm~\ref{alg_3})
  to uniform random values generated on the
  interval $[-R_m,R_m]$, with $R_m=1$ in general
  (or with $R_m$ set to a user-provided value).
  We observe from numerical experiments that the VarPro method is less sensitive or
  insensitive to the random coefficient initializations (the $R_m$ constant)
  than ELM. We provide numerical experiments in Section~\ref{sec:tests}
  for comparisons between the VarPro and the ELM methods.

\end{remark}


\subsection{Newton-Variable Projection Method for Solving Nonlinear PDEs}

We next develop a method based on variable projection
for solving nonlinear PDEs. The notations and settings here follow those
of Section~\ref{sec:lin}.

Consider the following nonlinear boundary value problem on
the domain $\Omega$ in $d$ dimensions,
\begin{subequations}\label{eq_15t}
  \begin{align}
    &
    Lu + F(u) = f(\mbs x),  \label{eq_15a} \\
    &
    Bu + G(u) = g(\mbs x), \quad \text{on}\ \partial\Omega, \label{eq_15b}
  \end{align}
\end{subequations}
where $F(u)$ and $G(u)$ are nonlinear operators on the solution field $u(\mbs x)$
and also possibly on its derivatives, and $L$, $B$, $f$ and $g$ have
the same meanings as in the equations~\eqref{eq_1a}--\eqref{eq_1b}.
We assume that the highest-order term occurs in the linear differential
operator $L$, and that the nonlinear terms $F(u)$ and $G(u)$ involve
only the lower-order derivatives (if any).
We again assume that the $L$ operator may involve time derivatives.
In such a case we treat the time $t$ in the same way as
the spatial coordinate $\mbs x$, as discussed in Section~\ref{sec:lin}.
We assume that this problem is well-posed.

We approximate the field solution $u(\mbs x)$ to the system~\eqref{eq_15t}
by a feed-forward neural network with $(L+1)$ layers,
following the same configurations and settings as
discussed in Section~\ref{sec:lin}.
Substituting the expansion relation~\eqref{eq_2} for $u(\mbs x)$
into equations~\eqref{eq_15a} and~\eqref{eq_15b},
and enforcing these two equations on all the collocation
points from $\mathbb{X}$ and on all the boundary collocation
points from $\mathbb{X}_b$ respectively, we arrive at an
algebraic system of $(N+N_b)$ equations
about the $(N_h+M)$ unknown
neural-network coefficients $(\bm\theta,\bm\beta)$.
We seek a least squares solution to this system, thus
leading to a nonlinear least squares problem.
This algebraic system, however, is nonlinear with respect to both
$\bm\theta$ and $\bm\beta$, because of the nonlinear terms
$F(u)$ and $G(u)$ in~\eqref{eq_15a}--\eqref{eq_15b}.
This is not a separable nonlinear least squares problem.
The variable projection approach apparently cannot be used
for solving this system, at least with the above straightforward formulation.


To circumvent the above issue and enable the use
of the variable projection strategy, we consider the linearization
of the system~\eqref{eq_15a}--\eqref{eq_15b} with the Newton's
method. Let $u^{k}$ denote the approximation of the
solution at the $k$-th Newton iteration. We linearize
this system as follows,
\begin{subequations}\label{eq_16}
  \begin{align}
    &
    Lu^{k+1} + F(u^k) + F'(u^k) \left(u^{k+1} - u^k \right) = f(\mbs x), \label{eq_16a} \\
    &
    Bu^{k+1} + G(u^k)+ G'(u^k)\left(u^{k+1} - u^k \right) = g(\mbs x),
    \quad \text{on}\ \partial\Omega, \label{eq_16b}
  \end{align}
\end{subequations}
where $F'(u)$ and $G'(u)$ denote the derivatives with respect to $u$.
We further re-write the linearized system into,
\begin{subequations}\label{eq_17}
  \begin{align}
    &
    Lu^{k+1} + F'(u^k) u^{k+1} = f(\mbs x) - F(u^k) + F'(u^k)u^k, \label{eq_17a} \\
    &
    Bu^{k+1} + G'(u^k)u^{k+1} = g(\mbs x) - G(u^k) + G'(u^k)u^k,
    \quad \text{on}\ \partial\Omega. \label{eq_17b}
  \end{align}
\end{subequations}
Given $u^k$, this system represents a linear boundary value problem
about the updated approximation field $u^{k+1}$.
Therefore, the VarPro/ANN algorithm developed in Section~\ref{sec:lin}
can be used to solve this linearized system~\eqref{eq_17a}--\eqref{eq_17b} for $u^{k+1}$.
Upon convergence of the Newton iteration, the solution to
the original nonlinear system~\eqref{eq_15a}--\eqref{eq_15b} will be obtained
and represented by the neural-network coefficients.


\begin{remark}\label{rem_5}
  It is important to notice that the above formulation leads to
  a linearized system about the updated approximation field  $u^{k+1}$ directly.
  This linearization form is crucial to the
  high  accuracy for solving nonlinear PDEs with
  the variable projection approach and
  artificial neural networks.

  An alternative and perhaps more commonly-used form of
  linearization for the Newton's method
  is often formulated in terms of the increment field.
  Let
  \begin{equation}\label{eq_18}
  u^{k+1} = u^k + v,
  \end{equation}
  where $v$ is the increment field at the step $k$. Then
  the increment is given by the following linearized system,
  \begin{subequations}\label{eq_19}
    \begin{align}
      &
      Lv + F'(u^k)v = f(\mbs x) - \left[Lu^k + F(u^k)  \right], \label{eq_19a} \\
      &
      Bv + G'(u^k)v = g(\mbs x) - \left[Bu^k + G(u^k)  \right],
      \quad \text{on}\ \partial\Omega. \label{eq_19b}
    \end{align}
  \end{subequations}
  So the increment field $v$ can be computed by the
  variable projection approach from the above system.
  The updated approximation $u^{k+1}$ is given by equation~\eqref{eq_18}.

  There are two issues with the form of linearization given by~\eqref{eq_19}
  when using variable projection and artificial neural networks.
  First, with this form $u^{k+1}$ can only be computed in the physical space (i.e.~on
  the collocation points), and it is not represented in terms of
  the neural network (i.e.~given by the network coefficients).
  Note that with the system~\eqref{eq_19} the increment field $v$
  is computed by the VarPro/ANN algorithm and is represented
  by the hidden-layer and output-layer coefficients of the neural network.
  But $u^{k+1}$ is computed by equation~\eqref{eq_18}. This can
  only be performed in the physical space, not in terms of the
  neural-network coefficients, due to the nonlinearity of the network output
  with respect to the hidden-layer coefficients.
  Second, upon convergence of the Newton iteration, the solution
  to the nonlinear system (i.e.~the converged $u^{k+1}$) is given in the physical
  space (on the collocation points), not represented by the neural
  network. Therefore, one needs to additionally
  convert this solution from physical
  space to the neural network representation, by solving a function
  approximation problem using the neural network and variable projection.
  This extra step is necessary in order to evaluate the solution field on the
  points other than the training collocation points in the domain.

  The form of linearization given by~\eqref{eq_16}, on the other hand, does not suffer
  from these issues. The updated approximation field $u^{k+1}$
  computed by the variable projection method
  is directly represented by the neural-network coefficients,
  as well as the solution to the original nonlinear system upon convergence.
  We observe that the solution obtained based on the formulation~\eqref{eq_16}
  is considerably more accurate, typically by two orders of magnitude or more,
  than that obtained based on the formulation~\eqref{eq_19}.
  It should be noted that the system~\eqref{eq_19} can be
  transformed into the system~\eqref{eq_16} by the substitution
  $v=u^{k+1}-u^k$.

\end{remark}


Within each Newton iteration we solve the linear boundary
value problem~\eqref{eq_17} using the variable projection method.
In order to make the following discussions more concise,
we introduce the following notation to drop the superscripts,
\begin{equation}\label{eq_20}
  \left\{
  \begin{split}
    &
    u(\mbs x) = u^{k+1}(\mbs x), \quad
    w(\mbs x) = u^k(\mbs x), \\
    &
    f_a(\mbs x) = f(\mbs x) - F(u^k) + F'(u^k)u^k, \quad
    g_a(\mbs x) = g(\mbs x) - G(u^k) + G'(u^k)u^k.
  \end{split}
  \right.
\end{equation}
Then the system~\eqref{eq_17} is re-written into,
\begin{subequations}\label{eq_21}
  \begin{align}
    &
    Lu + F'(w)u = f_a(\mbs x), \label{eq_21a}  \\
    &
    Bu + G'(w)u = g_a(\mbs x), \quad\text{on}\ \partial\Omega. \label{eq_21b}
  \end{align}
\end{subequations}

Let us next consider the solution of~\eqref{eq_21} with
the variable projection approach.
This system is similar to~\eqref{eq_1}. The solution
procedure mirrors that of Section~\ref{sec:lin}.
So we only summarize the most important steps below.
We use a feed-forward neural network to represent
the solution $u(\mbs x)$ to the system~\eqref{eq_21},
with the same settings and configurations for
the neural network and the collocation points
as given in Section~\ref{sec:lin}.
Substituting the expansion~\eqref{eq_2} into~\eqref{eq_21},
and enforcing the equation~\eqref{eq_21a} on all
the collocation points and equation~\eqref{eq_21b}
on all the boundary collocation points, we get
the following system in matrix form,
\begin{equation}\label{eq_22}
    \mbs H(\bm\theta)\bm\beta = \mbs S,
  \ 
  \text{where}\
  \mbs H(\bm\theta)=\begin{bmatrix}
  \vdots \\
  L\bm\Phi(\bm\theta,\mbs x_p) + F'(w)\bm\Phi(\bm\theta,\mbs x_p)  \\
  \vdots  \\ \hdashline[1pt/1pt]
  \vdots \\
  B\bm\Phi(\bm\theta,\mbs x_q) + G'(w)\bm\Phi(\bm\theta,\mbs x_q) \\
  \vdots
  \end{bmatrix}_{(N+N_b)\times M},
  \ 
  \mbs S = \begin{bmatrix}
  \vdots \\
  f_a(\mbs x_p)\\
  \vdots  \\ \hdashline[1pt/1pt]
  \vdots \\
  g_a(\mbs x_q)\\
  \vdots
  \end{bmatrix}_{(N+N_b)\times 1},
\end{equation}
where $\mbs x_p\in\mathbb{X}$ and $\mbs x_q\in\mathbb{X}_b$.
Following the same developments as given
by equations~\eqref{eq_7} and~\eqref{eq_9},
we arrive at the reduced nonlinear least squares problem~\eqref{eq_9}
about $\bm\theta$, with the understanding that
the terms $\mbs H(\bm\theta)$ and $\mbs S$
in all those equations are now defined by~\eqref{eq_22}.
We then invoke the Algorithm~\ref{alg_3} to
compute $(\bm\theta,\bm\beta)$, with the understanding that
on line $24$ of that algorithm the ``equation~\eqref{eq_4}''
is now replaced by equation~\eqref{eq_22} when computing $\bm\beta$.


\begin{remark}\label{rem_6}
  It should be noted that, depending on the form of the nonlinear operator $F(u)$,
  the terms $F'(w)\bm\Phi$ in the matrix $\mbs H(\bm\theta)$ may involve
  the derivatives of $\bm\Phi$. For example,
  with $F(u)=u\frac{\partial u}{\partial x}$,
  we have
  $
  F'(w)\bm\Phi = \frac{\partial w}{\partial x}\bm\Phi
  + w\frac{\partial\bm\Phi}{\partial x}.
  $
  The extra terms $F'(w)\bm\Phi$ and $G'(w)\bm\Phi$ in $\mbs H(\bm\theta)$ 
  do not add to the difficulty in computing the matrices $\mbs H(\bm\theta)$
  and $\mbs J_0(\bm\theta)$.
  Computing $\mbs H(\bm\theta)$ and $\mbs J_0(\bm\theta)$
  follows the same procedures as outlined
  in the Remarks~\ref{rem_1} and~\ref{rem_2}.
  The only difference lies in that in $\mbs H(\bm\theta)$ one needs to
  additionally compute
  the $F'(w)\bm\Phi(\bm\theta,\mbs x_p)$ ($\mbs x_p\in\mathbb{X}$)
  and $G'(w)\bm\Phi(\bm\theta,\mbs x_q)$ ($\mbs x_q\in\mathbb{X}_b$)
  based on the data for $\bm\Phi(\bm\theta,\mbs x)$ and its derivatives
  on the collocation points.
  In $\mbs J_0(\bm\theta)$ one needs to additionally compute
  the $\left.F'(w)\frac{\partial u}{\partial\bm\theta}\right|_{(\bm\theta,\mbs x_p)}$
  ($\mbs x_p\in\mathbb{X}$)
  and $\left.G'(w)\frac{\partial u}{\partial\bm\theta}\right|_{(\bm\theta,\mbs x_q)}$
  ($\mbs x_q\in\mathbb{X}_b$)
  based on the data for $u(\mbs x)$ and its derivatives with respect to $\mbs x$
  and $\bm\theta$ on the collocation points.
  
\end{remark}


To solve the nonlinear boundary value problem~\eqref{eq_15t},
we employ an overall Newton iteration.
Within each iteration we invoke the variable projection
method as given by Algorithm~\ref{alg_3} to solve
the system~\eqref{eq_17} for $u^{k+1}$, and the computed $u^{k+1}$ is
represented by the weight/bias coefficients of the artificial neural network.
Upon convergence of the Newton iteration,
the solution to the original nonlinear system~\eqref{eq_15t}
is given by the neural network, represented by
the neural-network coefficients.
In our implementation, we have considered two stopping criteria for
the Newton iterations, based on the Euclidean norms of
the residual vector $\mbs R$ and the increment vector $\Delta\mbs U$ defined by
\begin{equation}\label{eq_23}
  \mbs R = \begin{bmatrix}
    \vdots \\
    f(\mbs x_p) - Lu^k(\mbs x_p)-F(u^k(\mbs x_p)) \\
    \vdots \\ \hdashline[1pt/1pt] \vdots \\
    g(\mbs x_q) - Bu^k(\mbs x_q) - G(u^k(\mbs x_q)) \\
    \vdots
  \end{bmatrix}_{(N+N_b)\times 1},
  \qquad
  \Delta\mbs U=\begin{bmatrix}
  \vdots \\
  u^{k+1}(\mbs x_p) - u^k(\mbs x_p) \\
  \vdots
  \end{bmatrix}_{N\times 1}.
\end{equation}
The overall Newton-VarPro method with ANNs for solving
the nonlinear system~\eqref{eq_15t} is summarized
in the Algorithm~\ref{alg_4}.

\begin{algorithm}[tb]\label{alg_4}
  \DontPrintSemicolon
  \SetKwInOut{Input}{input}\SetKwInOut{Output}{output}
  \Input{input data $\mbs X$ to neural network;
    source data $f(\mbs x_p)$ ($\mbs x_p\in\mathbb{X}$)
    and $g(\mbs x_q)$ ($\mbs x_q\in\mathbb{X}_b$);
    initial guess $u^0(\mbs x)$.
  }
  \Output{solution $u(\mbs x)$, represented by the coefficients
    of the neural network.}
  \BlankLine\BlankLine
  \For{$k\leftarrow 0$ \KwTo maximum-number-of-newton-iterations}{
    compute the vector $\mbs R$ by~\eqref{eq_23}\;
    \If{$\|\mbs R \|$ is below a tolerance}{ break\; }
    \BlankLine
    compute $\mbs S$ by the equations~\eqref{eq_22} and~\eqref{eq_20}\;
    call Algorithm~\ref{alg_3}, with ``equation~\eqref{eq_4}'' on line 24
    therein replaced by ``equation~\eqref{eq_22}'', to obtain
    $(\bm\theta,\bm\beta)$, which are the neural-network
    representation of $u^{k+1}$ in the system~\eqref{eq_17}\;
    \BlankLine
    compute the vector $\Delta\mbs U$ by~\eqref{eq_23}\;
    \If{$\|\Delta \mbs U \|$ is below a tolerance}{ break\; }
  }
  
  \caption{Newton-variable projection algorithm for the
  nonlinear problem~\eqref{eq_15t}.}
\end{algorithm}


\begin{remark}\label{rem_7}
  When solving the nonlinear problem~\eqref{eq_15t} using the Newton-VarPro
  method (Algorithm~\ref{alg_4}), one can
  often turn off the initial-guess perturbations and
  sub-iterations in Algorithm~\ref{alg_3}
  when invoking this algorithm to solve the system~\eqref{eq_17}.
  This can be achieved by simply setting the ``maximum-number-of-sub-iterations''
  to zero on line $5$ of the Algorithm~\ref{alg_3}.
  In this case, if the converged cost from Algorithm~\ref{alg_3} is
  not very small (above some threshold), this means that
  the returned $u^{k+1}$ solution from that Newton step
  may not be that accurate. This inaccuracy, however, can be offset
  by the subsequent Newton iterations.

\end{remark}

\begin{remark}\label{rem_8}
  When solving nonlinear PDEs using the Newton-VarPro method,
  if the resolution is low (e.g.~using a small number of training collocation points),
  we observe from numerical experiments that the Newton iteration
  may have difficulty reaching convergence within a specified maximum number
  of iterations. In such a case, increasing the resolution (e.g.~increasing the number of
  collocation points) can typically improve the convergence of
  the Newton iteration.
  In the numerical experiments of Section \ref{sec:tests}
  we typically employ a relative tolerance $1E-8$ for the Newton
  iterations (see the lines $3$ and $9$ of Algorithm~\ref{alg_4}).

\end{remark}

\section{Numerical Examples}
\label{sec:tests}

We use several numerical examples involving linear
and nonlinear PDEs to illustrate the performance characteristics
of the VarPro method.
These problems are in two spatial dimensions or in one spatial
dimension plus time.
We also compare the simulation results of the current VarPro method
and the ELM method from~\cite{DongY2021,DongL2021} to demonstrate
the superior accuracy of the current method.

As stated previously, the current VarPro method is implemented in Python,
using the Tensorflow (www.tensorflow.org) and Keras (keras.io)
libraries. For the linear least squares method we employ
the scipy routine ``scipy.linalg.lstsq'' in our implementation,
which invokes the corresponding routine from the LAPACK library.
For the nonlinear least squares method, we employ
the scipy routine ``scipy.optimize.least\_squares'' in our application code,
which implements the Gauss-Newton method together with
a trust region algorithm~\cite{BranchCL1999}.
The differential operators acting on the output fields of
the last hidden layer are computed by a forward-mode auto-differentiation
in our implementation, as discussed in Remark~\ref{rem_1}.
The data for the Jacobian matrix are
computed by a reverse-mode auto-differentiation using
a vectorized map together with the GradientTape in Tensorflow,
as discussed in Remark~\ref{rem_2}.
We would like to mention that in our implementation of the neural network,
between the input layer and the first hidden layer,
we have incorporated a lambda layer from Keras to normalize
the input data $\mbs X$ from the rectangular
domain $\Omega=[a_1,b_1]\times\dots\times[a_d,b_d]$
to the standard domain $[-1,1]^d$.

As in our previous works~\cite{DongL2021,DongL2021bip,DongY2021},
we employ a fixed seed value for the random number generators
in the numerical experiments in each subsection, so that
the reported results here can be exactly reproducible.
We use the same seed for the random number generators from
the Tensorflow library and from the numpy package.
These seed values are $1$ in Sections~\ref{sec:poisson} and~\ref{sec:nonl_helm},
$10$ in Sections~\ref{sec:advec} and~\ref{sec:burgers},
and $22$ in Section~\ref{sec:sg}.

\subsection{Linear Examples}

\subsubsection{Poisson Equation}
\label{sec:poisson}

\begin{figure}
  \centerline{
    \includegraphics[width=2.5in]{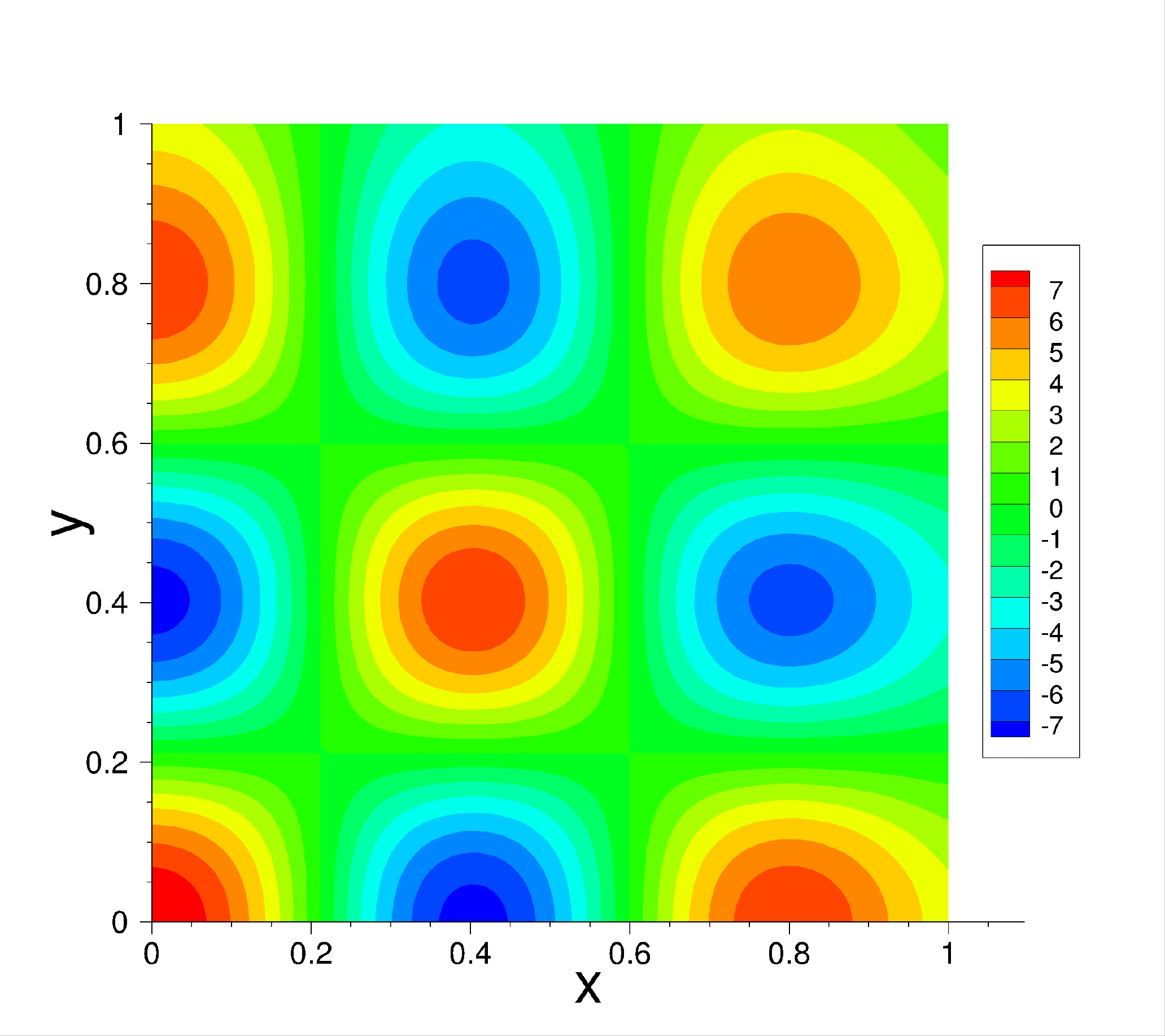}(a)
    \includegraphics[width=2.5in]{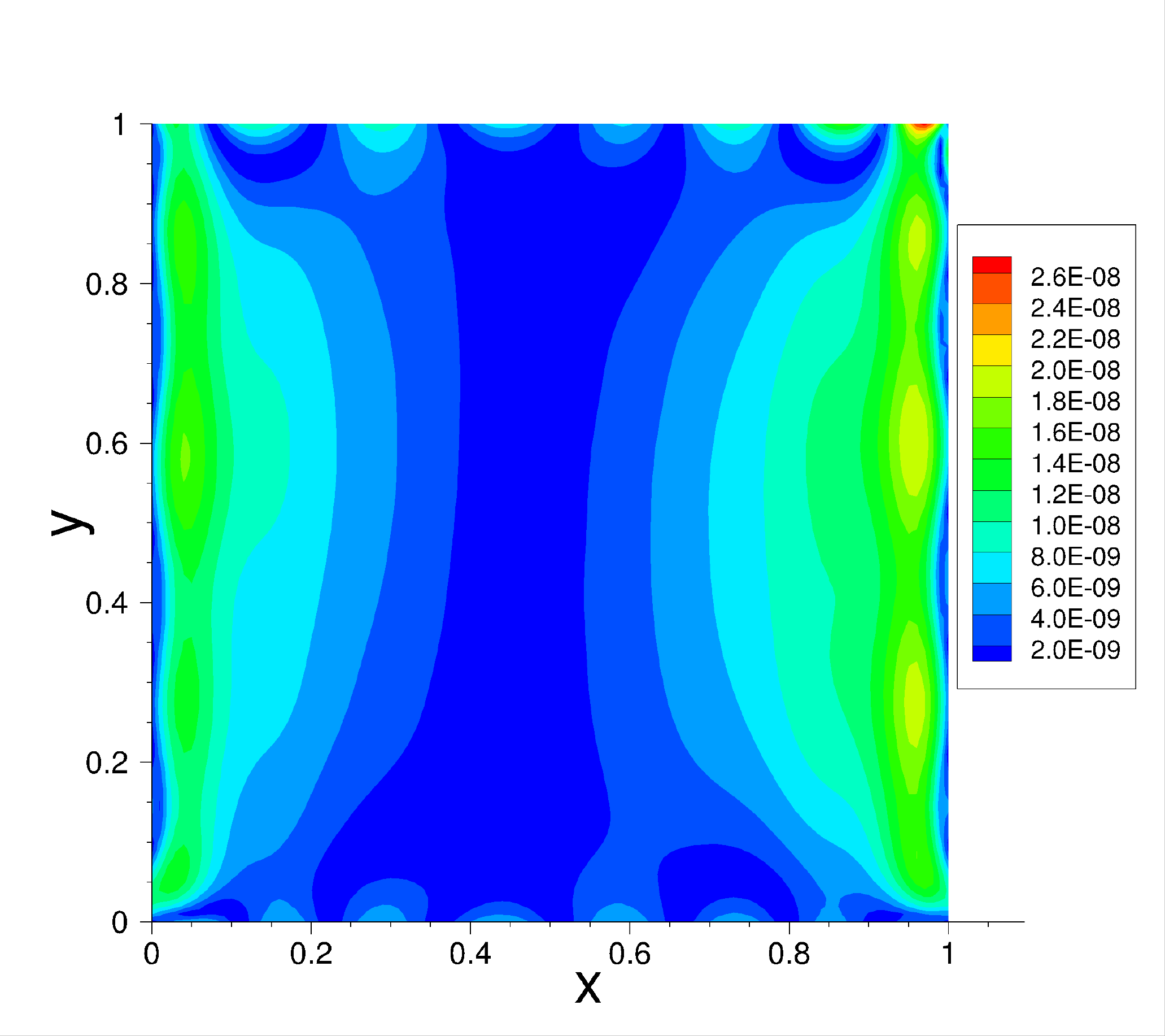}(b)
  }
  \caption{Poisson equation:
    (a) Distribution of the exact solution.
    (b) Distribution of the absolute error of the VarPro solution.
    In (b), neural network [2, 125, 1], ``$\cos$'' activation function,
    $Q=15\times 15$ uniform training collocation points.
  }
  \label{fg_1}
\end{figure}


We first consider the canonical two-dimensional (2D) Poisson equation on
a unit square domain, $(x,y)\in[0,1]\times[0,1]$,
\begin{subequations}\label{eq_24}
  \begin{align}
    & \frac{\partial^2 u}{\partial x^2} + \frac{\partial^2u}{\partial y^2} = f(x,y), \\
    & u(x,0) = g_1(x), \quad u(x,1) = g_2(x), \quad
    u(0,y) = g_3(y), \quad u(1,y) = g_4(y),
  \end{align}
\end{subequations}
where $u(x,y)$ is the field function to be solved for, $f(x,y)$ is
a prescribed source term, and $g_i$ ($1\leqslant i\leqslant 4$) denote
the boundary data.
We employ the following analytic solution to this problem in the tests,
\begin{multline}\label{eq_25}
  u = \left[2\cos\left(\frac32\pi x+\frac25\pi \right)
    +\frac32\cos\left(3\pi x-\frac{\pi}{5} \right)
    + \frac{1}{1+x^2} \right]
  \left[2\cos\left(\frac32\pi y+\frac25\pi \right) \right. \\
    \left.
    +\frac32\cos\left(3\pi y-\frac{\pi}{5} \right)
    + \frac{1}{1+y^2} \right],
\end{multline}
by choosing the source term $f$ and the boundary data $g_i$ ($1\leqslant i\leqslant 4$)
appropriately. Figure~\ref{fg_1}(a) shows the
distribution of this analytic solution in the $xy$ plane.


We employ feed-forward neural networks with one or two hidden layers, with
the architecture given by $[2, M, 1]$ or $[2, 20, M, 1]$,
where $M$ is varied systematically or fixed at $M=100$, $125$ or $200$.
The two input nodes represent the coordinates $(x,y)$, and the single output node
represents the solution field $u(x,y)$.
The activation function for the hidden nodes is either
the cosine function, $\sigma(x)=\cos(x)$, or the Gaussian function,
$\sigma(x) = e^{-x^2}$.
The output layer is required to be linear (no activation function)
and contain no bias.

We employ a uniform set of $Q=Q_1\times Q_1$ grid points on the domain as
the training collocation points, where $Q_1$ denotes the number of uniform
grid points in each direction (including the two end points)
and is varied systematically
between around $5$ and $30$ in the tests.
After the neural network is trained by the VarPro method on
the $Q_1\times Q_1$ collocation points,
the neural network is evaluated on a much larger uniform set of
$Q_2\times Q_2$ grid points, where $Q_2=101$ for this problem,
to obtain the solution $u(x,y)$. This solution is compared with
the analytic solution~\eqref{eq_25} to compute the maximum ($L^{\infty}$)
and the root-mean-squares (rms, or $L^2$) errors.
These maximum/rms errors are then recorded and referred
to as the errors associated with the given neural network architecture and
the training collocation points $Q=Q_1\times Q_1$ for the VarPro method.


\begin{table}[tb]
  \centering
  \begin{tabular}{ll | ll}
    \hline
    parameter & value & parameter & value \\ \hline
    neural network & $[2, M, 1]$ or $[2, 20, M, 1]$ & training points $Q$ & $Q_1\times Q_1$ \\
    $M$ & varied & $Q_1$ & varied \\
    activation function & $\cos$, Gaussian & testing points & $Q_2\times Q_2$ \\
    random seed & $1$ & $Q_2$ & 101 \\
    initial guess $\bm\theta_0$ & random values on $[-R_m,R_m]$ & $R_m$ & $1.0$ \\
    $\delta$ (Algorithm~\ref{alg_3}) & $0.5$, $1.0$, $2.0$, $5.0$, or $7.0$ 
    & $p$ (Algorithm~\ref{alg_3}) & $0.5$ \\
    max-subiterations & $5$ & threshold (Algorithm~\ref{alg_3}) & $1E-12$ \\
    \hline
  \end{tabular}
  \caption{Poisson equation: Main simulation parameters of the VarPro method.}
  \label{tab_1}
\end{table}

The main simulation parameters of the VarPro method are summarized in
Table~\ref{tab_1}. The last three rows of this table pertain to
the parameters in Algorithm~\ref{alg_3}.
$\bm\theta_0$ denotes the initial guess to
the hidden-layer coefficients in Algorithm~\ref{alg_3},
which are set to uniform random values
generated on $[-R_m,R_m]$ with $R_m=1.0$. When comparing the VarPro method
and the ELM method, we also employ a value $R_m=R_{m0}$ with VarPro, where $R_{m0}$
is the optimal $R_m$ value corresponding to ELM computed using
the method from~\cite{DongY2021}.
$\delta$ and $p$ in this table are the maximum perturbation magnitude and
the preference probability in Algorithm~\ref{alg_3}, respectively.
The ``max-subiterations'' here refers to the maximum-number-of-sub-iterations
in Algorithm~\ref{alg_3}. The ``threshold'' here refers to the threshold on
the lines $3$ and $19$ in Algorithm~\ref{alg_3}.

\begin{figure}
  \centerline{
    \includegraphics[width=2.5in]{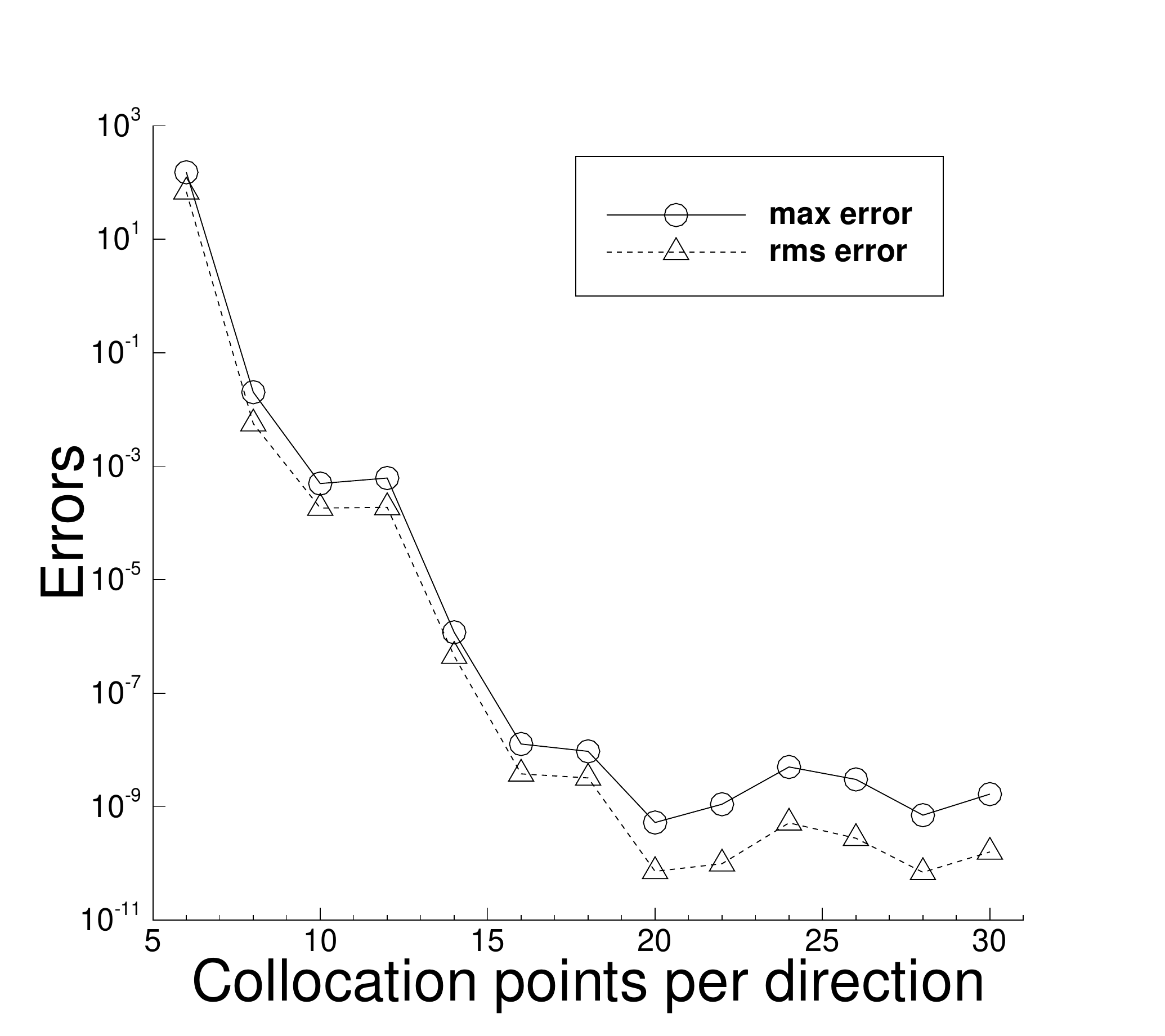}(a)
    \includegraphics[width=2.5in]{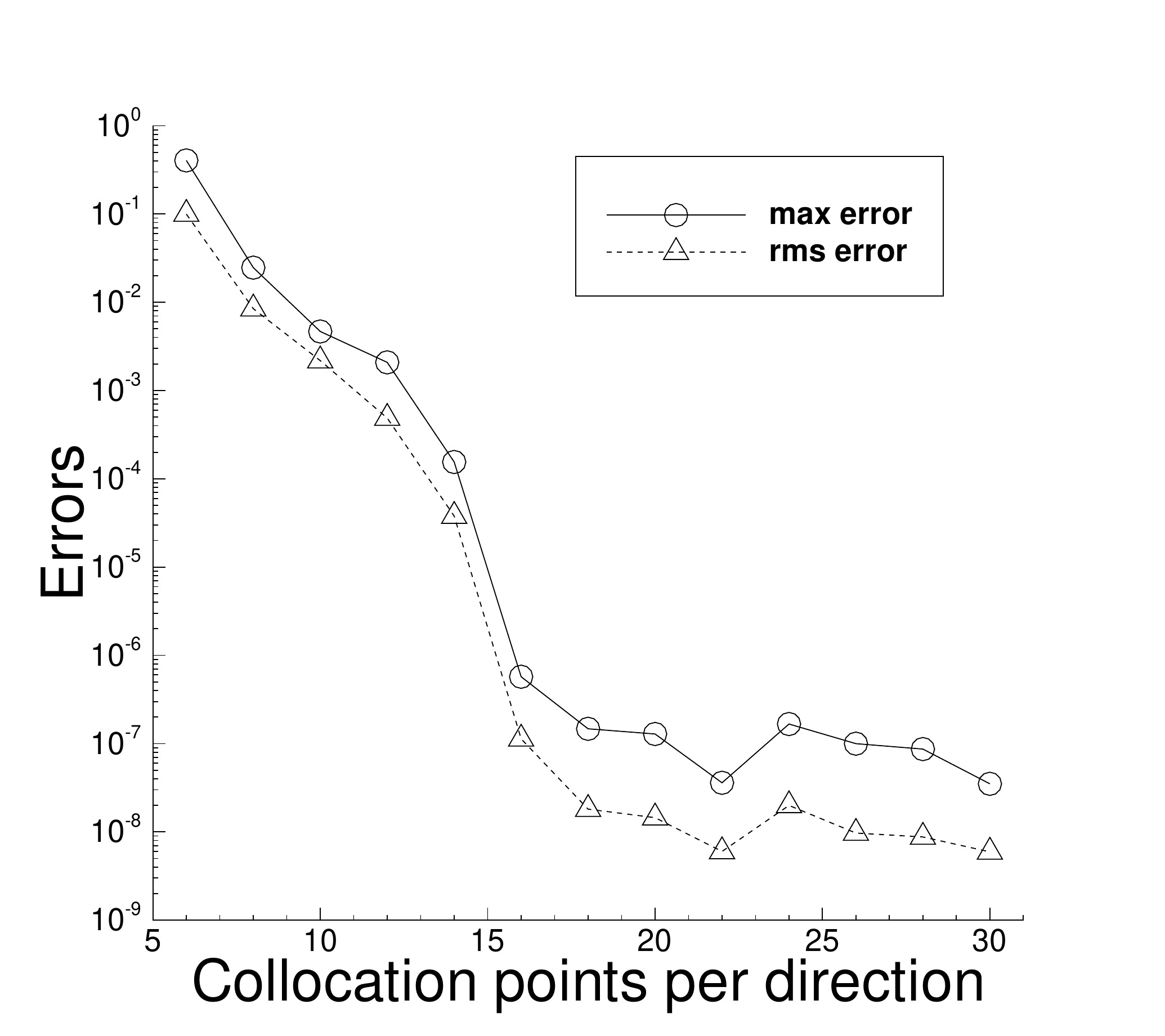}(b)
  }
  \caption{Poisson equation:
    Maximum/rms errors of the VarPro solution
    versus the number of collocation
    points per direction ($Q_1$) obtained with (a) the $\cos$ activation
    function, and (b) the Gaussian activation function.
    Neural network [2, 200, 1] in (a,b); $Q_1$ is varied in (a,b);
    $\delta = 7.0$ in (a), and $\delta = 1.0$ in (b).
  }
  \label{fg_2}
\end{figure}

Let us first consider the VarPro results obtained with neural networks containing
a single hidden layer.
Figure~\ref{fg_1}(b) shows the distribution of the absolute error of
the VarPro solution in the $xy$ plane.
This result corresponds to the neural network architecture
$[2, 125, 1]$, with the $\cos$ activation
function, a uniform set of $Q=15\times 15$ training collocation points,
and $\delta = 7.0$ in Algorithm~\ref{alg_3}.
The VarPro solution is highly accurate, with a maximum error around $10^{-8}$
in the domain.

Figure~\ref{fg_2} illustrates the convergence behavior of the VarPro solution
as a function of the number of training collocation points in the domain.
Here we employ a neural network $[2, 200, 1]$, with the $\cos$ and
the Gaussian activation functions. The number of collocation points
in each direction ($Q_1$) is varied systematically.
Figure \ref{fg_2} shows the maximum and the rms errors of the VarPro solution
in the domain as a function of $Q_1$, obtained using the
$\cos$ activation function (plot (a)) and the Gaussian activation function (plot (b)).
The VarPro errors decrease approximately exponentially when $Q_1$
is below around $20$, and then appear to stagnate as $Q_1$ further increases.
The VarPro errors reach a level around $10^{-11}\sim 10^{-9}$ with
the $\cos$ activation function and a level around $10^{-9}\sim 10^{-7}$
with the Gaussian activation function.

\begin{figure}
  \centerline{
    \includegraphics[width=2.5in]{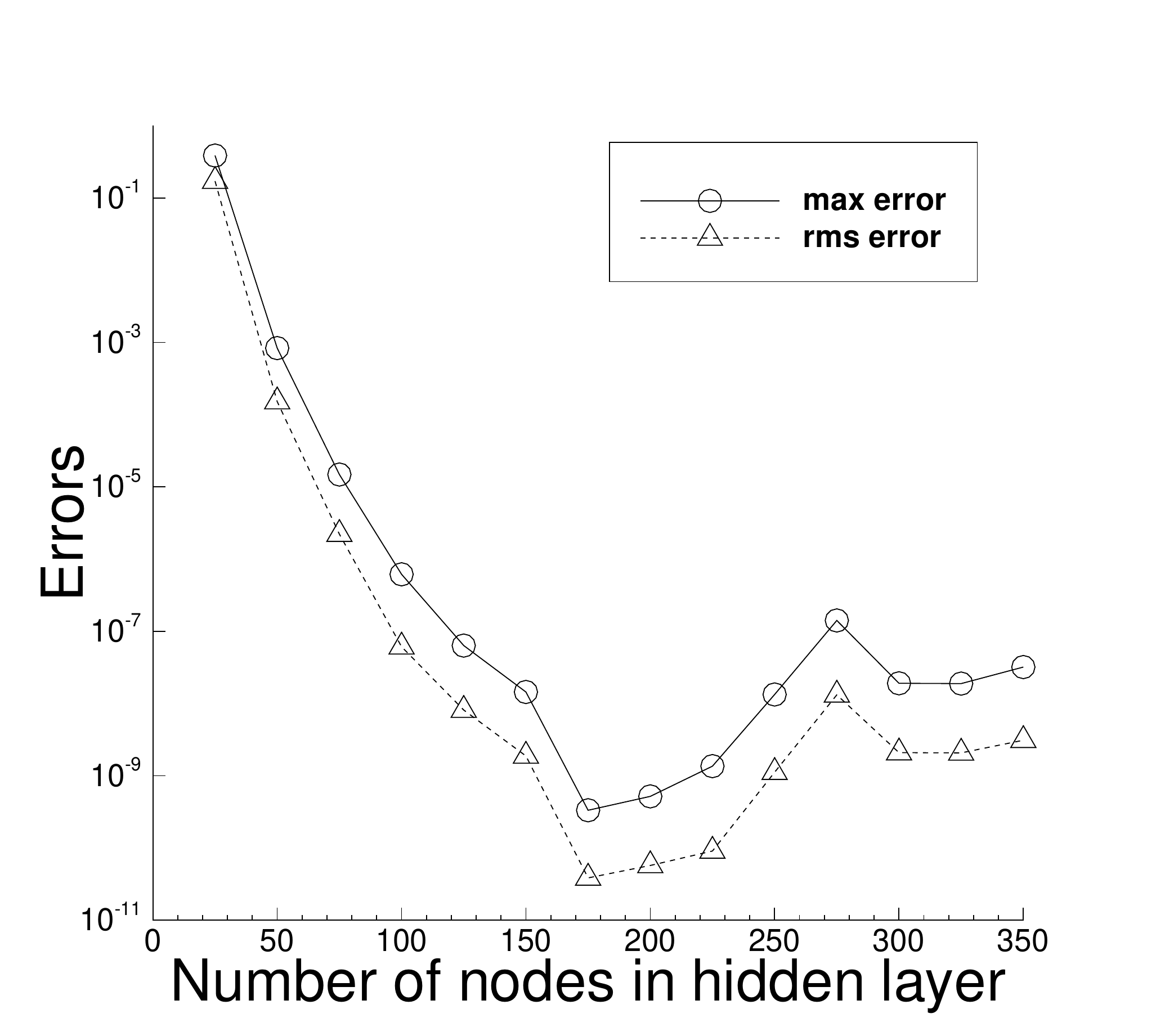}(a)
    \includegraphics[width=2.5in]{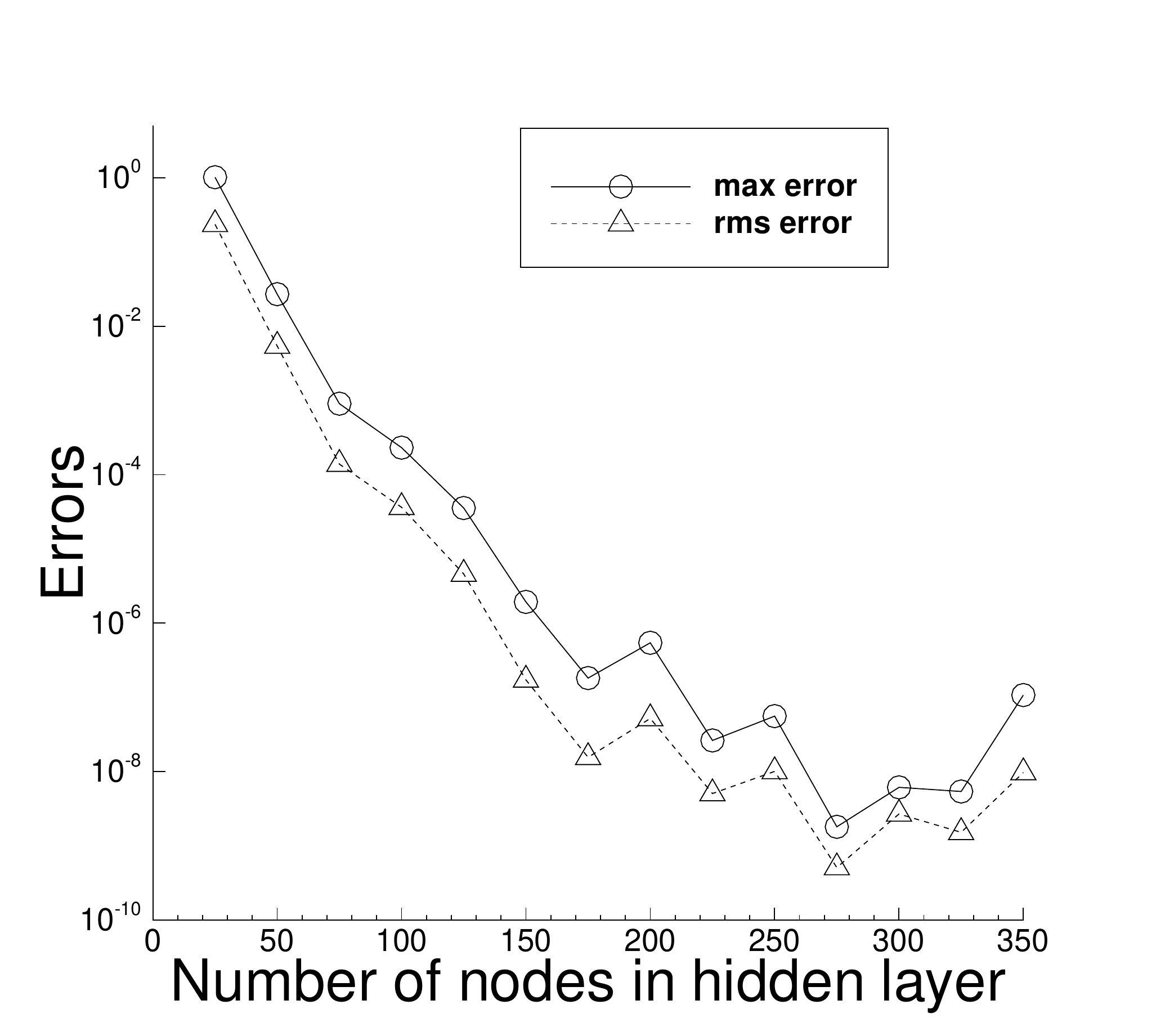}(b)
  }
  \caption{Poisson equation:
    The maximum/rms errors of the VarPro solution
    versus the number of nodes in the hidden layer ($M$),
    computed using (a) the $\cos$  and (b) the Gaussian
    activation functions.
    Neural network [2, $M$, 1], where $M$ is varied in (a,b);
    $Q=21\times 21$ in (a,b);
    $\delta = 7.0$ in (a) and $\delta = 2.0$ in (b).
  }
  \label{fg_3}
\end{figure}

Figure \ref{fg_3} illustrates the convergence behavior of the VarPro accuracy
with respect to the number of nodes in the hidden layer ($M$) of the network.
Here we consider neural networks with the architecture $[2,M,1]$,
where $M$ is varied systematically, with the $\cos$ and Gaussian activation
functions. A fixed uniform set of $Q=21\times 21$ training collocation points is used.
Figure \ref{fg_3} shows the maximum/rms errors of the VarPro solution in the domain
as a function of $M$, obtained with the $\cos$ (plot(a)) and the Gaussian (plot (b))
activation functions.
One can observe an approximately exponential decrease in the VarPro errors
with increasing $M$ (when $M$ is below a certain value),
and then the errors appear to stagnate (or increase slightly) as
$M$ further increases.

\begin{figure}
  \centerline{
    \includegraphics[width=2.5in]{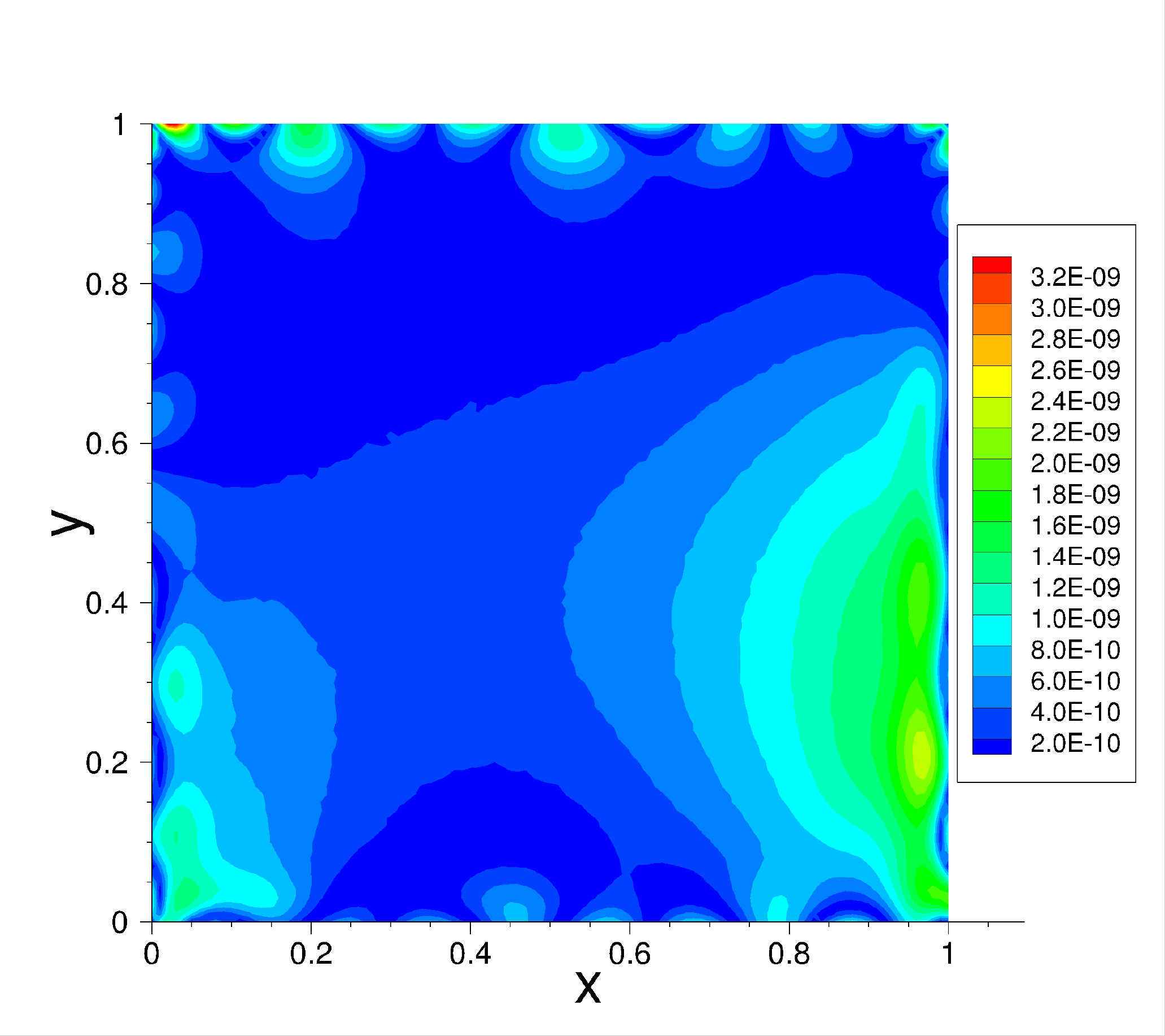}(a)
    \includegraphics[width=2.5in]{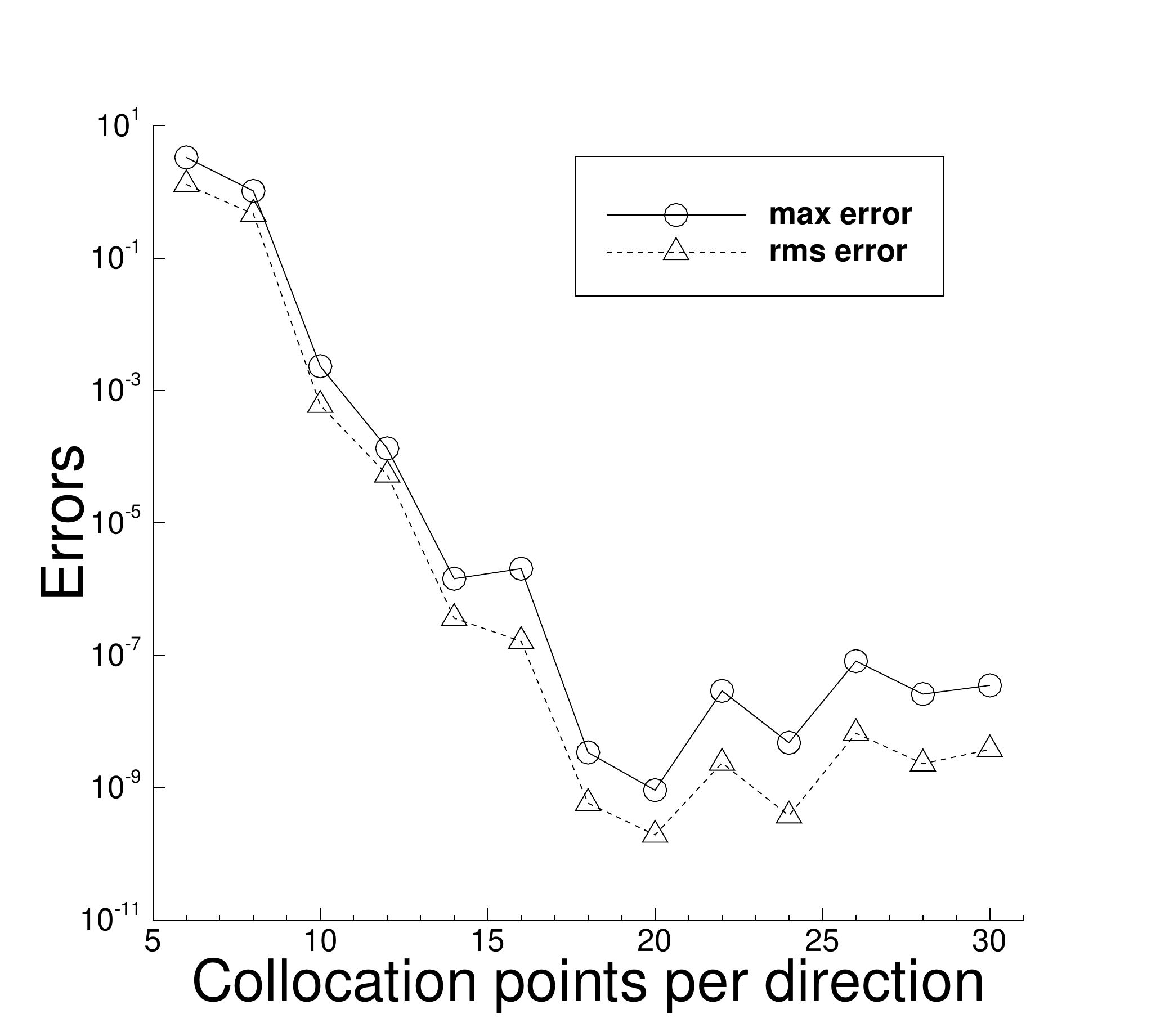}(b)
  }
  \caption{Poisson equation ($2$ hidden layers in neural network):
    (a) Error distribution of the VarPro solution.
    (b) The maximum/rms errors of the VarPro solution
    versus the number of collocation points per direction ($Q_1$).
    Neural network [2, 20, 100, 1], with the $\cos$ activation function.
    $Q=18\times 18$ in (a), and is varied in (b).
    $\delta = 0.5$ in (a,b).
  }
  \label{fg_4}
\end{figure}

Figure \ref{fg_4} illustrates the VarPro solution using a neural network
containing two hidden layers.
Here we have employed a neural network with the architecture
$[2, 20, 100, 1]$ and the $\cos$ activation function for all the
hidden nodes. Figure~\ref{fg_4}(a) shows the distribution of
the absolute error of the VarPro solution obtained with
a set of $Q=18\times 18$ uniform collocation points in
the domain. The maximum error is on the level $10^{-9}$,
indicating a high accuracy.
Figure \ref{fg_4}(b) depicts the maximum/rms errors of the VarPro
solution as a function of the number of collocation points
in each direction ($Q_1$). One can again
observe an exponential decrease in the errors (before saturation)
with increasing number of collocation points.
All these results suggest that the VarPro method produces
highly accurate results for solving the Poisson equation.


\begin{table}[tb]
  \centering
  \begin{tabular}{l | l | l| ll| ll}
    \hline
    Neural & $[-R_m,R_m]$ & collocation & VarPro & & ELM & \\ \cline{4-7}
    network &  & points & max-error & rms-error & max-error & rms-error \\ \hline
    [2, 100, 1] & $R_m=1$ & $5\times 5$ & $6.470E+1$ & $2.422E+1$ & $1.247E+2$ & $2.745E+1$  \\
    & & $10\times 10$ & $8.388E-3$ & $3.941E-3$ & $1.402E+1$ & $2.575E+0$  \\
    & & $15\times 15$ & $6.018E-7$ & $8.241E-8$ & $1.475E+1$ & $1.938E+0$  \\
    & & $20\times 20$ & $3.693E-7$ & $4.216E-8$ & $1.690E+1$ & $2.527E+0$ \\
    & & $25\times 25$ & $5.845E-7$ & $8.054E-8$ & $1.777E+1$ & $2.752E+0$ \\
    & & $30 \times 30$ & $2.688E-7$ & $2.867E-8$ & $1.864E+1$ & $2.916E+0$ \\ \cline{2-7}
    & $R_m=R_{m0}$ & $5\times 5$ & $2.156E+0$ & $6.114E-1$ & $2.156E+0$ & $6.114E-1$ \\
    & $\qquad=6$ & $10\times 10$ & $7.735E-4$ & $1.497E-4$ & $1.353E-1$ & $2.497E-2$ \\
    & & $15\times 15$ & $4.175E-7$ & $4.614E-8$ & $3.019E-1$ & $6.004E-2$ \\
    & & $20\times 20$ & $1.753E-7$ & $1.974E-8$ & $3.859E-1$ & $7.575E-2$ \\
    & & $25\times 25$ & $8.443E-7$ & $1.227E-7$ & $4.336E-1$ & $8.489E-2$ \\
    & & $30\times 30$ & $8.709E-8$ & $1.088E-8$ & $4.673E-1$ & $9.157E-2$ \\
    \hline
    [2, 200, 1] & $R_m=1$ & $5\times 5$ & $5.948E+0$ & $2.102E+0$ & $9.325E+1$ & $2.024E+1$  \\
    & & $10\times 10$ & $2.127E-2$ & $4.398E-3$ & $4.417E+0$ & $7.083E-1$ \\
    & & $15\times 15$ & $6.082E-8$ & $1.983E-8$ & $5.615E+0$ & $8.019E-1$ \\
    & & $20\times 20$ & $1.459E-9$ & $1.203E-10$ & $4.979E+0$ & $7.610E-1$ \\
    & & $25\times 25$ & $1.782E-7$ & $7.978E-8$ & $5.077E+0$ & $8.565E-1$ \\
    & & $30 \times 30$ & $3.000E-9$ & $3.420E-10$ & $5.633E+0$ & $8.804E-1$  \\ \cline{2-7}
    & $R_m=R_{m0}$ & $5\times 5$ & $7.292E-1$ & $2.733E-1$ & $7.292E-1$ & $2.732E-1$ \\
    & $\qquad=6$ & $10\times 10$ & $1.283E-4$ & $3.705E-5$ & $1.283E-4$ & $3.705E-5$  \\
    & & $15\times 15$ & $2.315E-9$ & $4.746E-10$ & $9.822E-6$ & $7.481E-7$  \\
    & & $20\times 20$ & $3.449E-10$ & $3.722E-11$ & $1.245E-5$ & $1.518E-6$  \\
    & & $25\times 25$ & $6.379E-9$ & $4.980E-10$ & $1.174E-5$ & $1.677E-6$  \\
    & & $30\times 30$ & $4.221E-10$ & $4.086E-11$ & $1.242E-5$ & $1.800E-6$  \\
    \hline
  \end{tabular}
  \caption{Poisson equation: comparison of the maximum/rms errors 
    obtained using the VarPro and ELM methods.
    $\cos$ activation function.
    In both VarPro and ELM, the hidden-layer coefficients are
    initialized/set to uniform random values
    generated on $[-R_m,R_m]$, with $R_m=1.0$ or with $R_m=R_{m0}$.
    $R_{m0}$ is the optimal $R_m$ for ELM computed using the method
    from~\cite{DongY2021}, and in this case $R_{m0}=6.0$.
    $\delta = 5.0$ in VarPro.
  }
  \label{tab_2}
\end{table}

Table~\ref{tab_2} compares the errors of 
the current VarPro method and the ELM method~\cite{DongY2021,DongL2021}
for solving the Poisson equation.
We have considered two neural networks having the
architecture $[2, M, 1]$, with $M=100$ and $M=200$.
A uniform set of $Q=Q_1\times Q_1$ training collocation points are
employed on the domain, where $Q_1$ is varied between $Q_1=5$ and $Q_1=30$.
In ELM the hidden-layer coefficients are set (and fixed) to uniform random values
generated on $[-R_m,R_m]$, and in VarPro the hidden-layer coefficients
are initialized (i.e.~the initial guess $\bm\theta_0$ in Algorithm~\ref{alg_3})
to the same random values from $[-R_m,R_m]$.
So the random hidden-layer coefficients in ELM and the
initial hidden-layer coefficients in VarPro are identical.
We have considered two $R_m$ values, $R_m=1.0$ and $R_m=R_{m0}$,
where $R_{m0}$ is the optimal $R_m$ for ELM computed using the method
from~\cite{DongY2021} and in this case $R_{m0}= 6.0$.
We can make the following observations:
\begin{itemize}
\item
  The VarPro method in general produces considerably more accurate
  results than ELM, under the same settings and conditions,
  especially when the size of the neural network is still not quite large.

\item
  The ELM accuracy has a fairly strong dependence on the $R_m$ value.
  On the other hand, the VarPro accuracy is less sensitive or insensitive to
  the $R_m$ value.

\end{itemize}
In these tests the VarPro method has produced errors on
the order $10^{-8}\sim 10^{-10}$ with the given neural networks.
We should point out that the ELM method can also achieve numerical
errors on such levels, but it requires neural networks with
a larger number of nodes in the hidden layer.

\subsubsection{Advection Equation}
\label{sec:advec}

\begin{figure}
  \centering
    \includegraphics[width=3.5in]{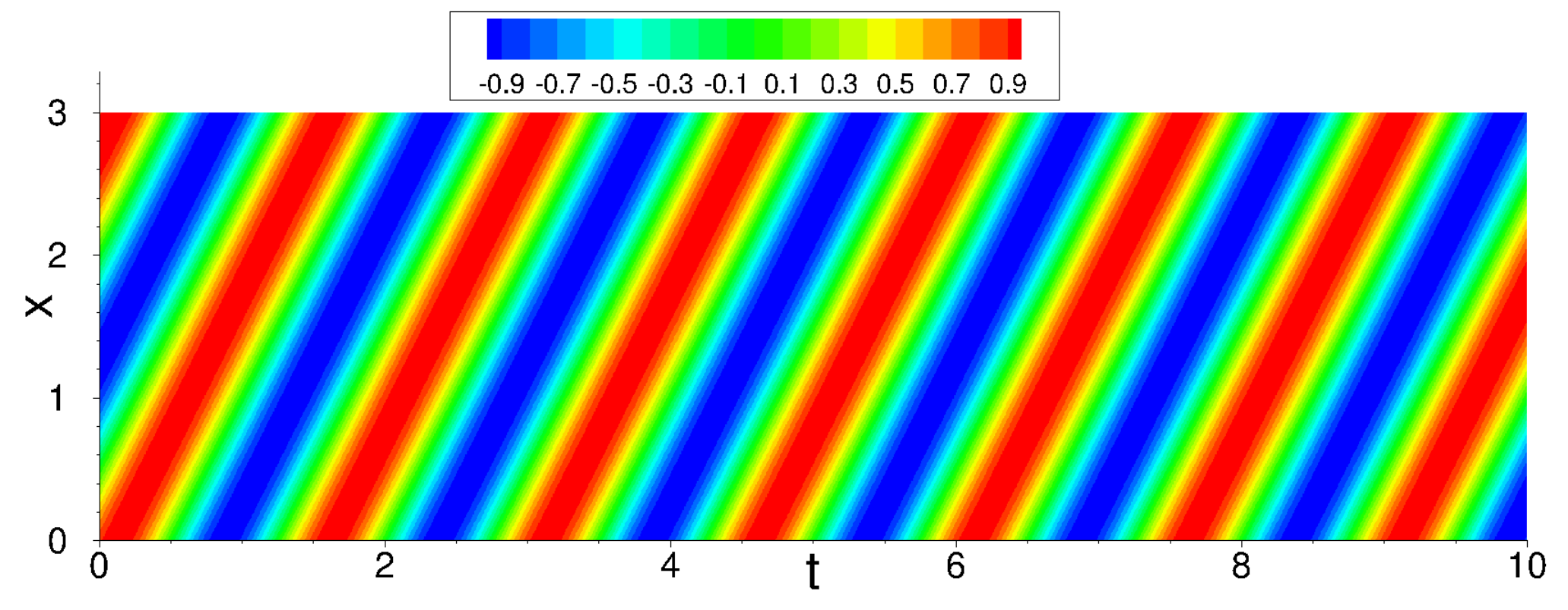}(a)
    \includegraphics[width=3.5in]{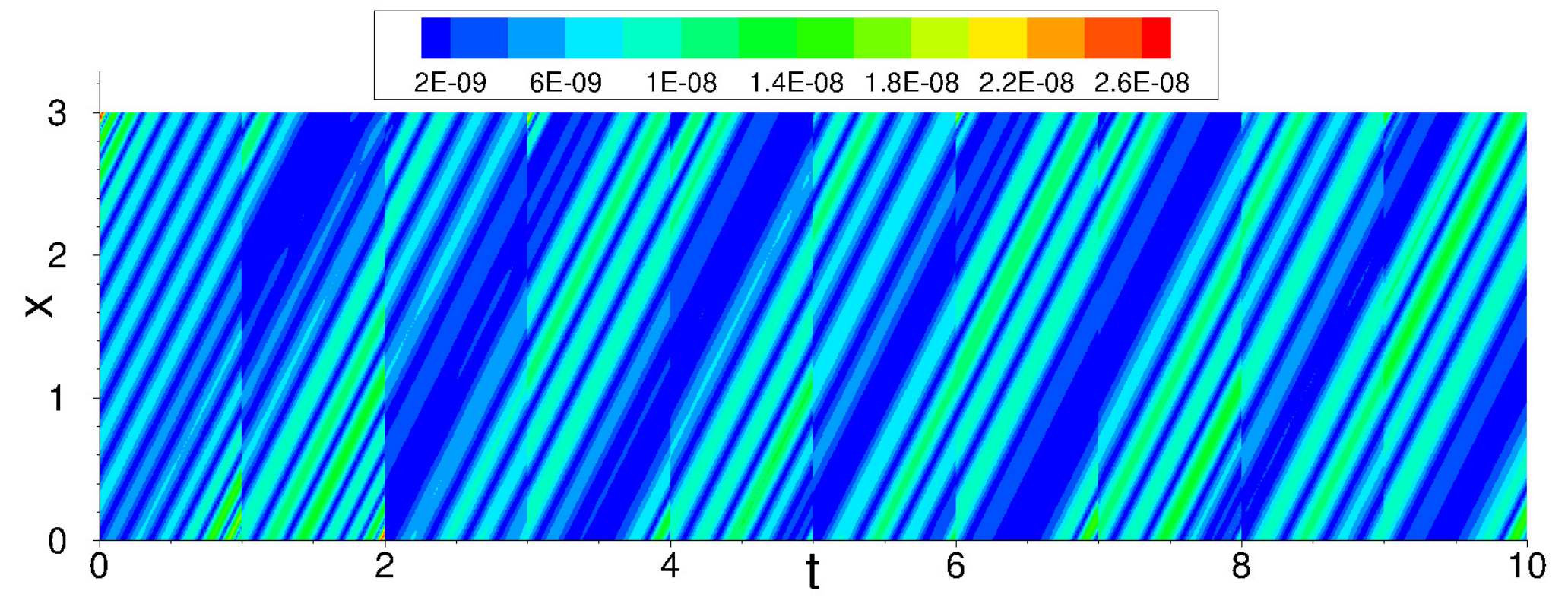}(b)
    \caption{Advection equation:
      Distributions of (a) the exact solution, and (b) the absolute error of
      the VarPro solution, in the spatial-temporal plane.
      $t_f=10$. In (b), $10$ time blocks, $Q=21\times 21$ uniform collocation points per time block,
      neural network [2, 100, 1], Gaussian activation function,
      the max sub-iterations is 2 and $\delta = 3.0$ in VarPro.
    }
    \label{fg_5}
\end{figure}

As another linear example we consider the spatial-temporal domain
$\Omega=\{(x,t)\ |\ x\in[0,3],\ t\in[0,t_f]\}$ in this test, where
the temporal dimension $t_f$ is to be specified below.
We consider the initial/boundary value problem with
the advection equation on $\Omega$,
\begin{subequations}\label{eq_26}
  \begin{align}
    & \frac{\partial u}{\partial t} - c\frac{\partial u}{\partial x} = 0, \\
    & u(0,t) = u(3,t), \\
    & u(x,0) = \sin\frac{2\pi}{3}(x-2),
  \end{align}
\end{subequations}
where $u(x,t)$ is the field function to be solved for,
and $c=-2.0$ is the wave speed.
This problem has the following exact solution,
\begin{equation}\label{eq_28}
  u(x,t) = \sin\frac{2\pi}{3}(x-2t-2).
\end{equation}
Figure \ref{fg_5}(a) shows the distribution of this exact
solution in the spatial-temporal domain with $t_f=10$.


To solve this problem with the VarPro method, we
employ a feed-forward neural network with one or two hidden layers,
with the architecture given by $[2, M, 1]$ or $[2, 10, M, 1]$,
where $M$ is varied systematically in the tests.
The two input nodes represent the spatial/temporal
coordinates $(x,t)$, and the output node represents
the solution field $u(x,t)$.
We employ the Gaussian function, $\sigma(x)=e^{-x^2}$,
or the Gaussian error linear unit (GELU)~\cite{HendrycksG2020},
$\sigma(x)=\frac12 x\left[1+\text{erf}\left(\frac{x}{\sqrt{2}} \right) \right]$,
as the activation function for all the hidden nodes.
The output layer is linear and with zero bias.

\begin{table}[tb]
  \centering
  \begin{tabular}{ll | ll}
    \hline
    parameter & value & parameter & value \\ \hline
    $t_f$ & $10$, or $100$ & number of time blocks & $10$, or $100$ \\
    neural network & $[2, M, 1]$, or $[2, 10, M, 1]$ & training points $Q$ & $Q_1\times Q_1$ \\
    $M$ &  varied & $Q_1$ & varied \\
    activation function & Gaussian, GELU & testing points & $Q_2\times Q_2$ \\
    random seed & $10$ & $Q_2$ & 101 \\
    initial guess $\bm\theta_0$ & random values on $[-R_m,R_m]$ & $R_m$ & $1.0$ \\
    $\delta$ (Algorithm~\ref{alg_3}) & $0.0$, $0.05$, $1.0$, or $3.0$ 
    & $p$ (Algorithm~\ref{alg_3}) & $0.5$ \\
    max-subiterations & $0$, or $2$ & threshold (Algorithm~\ref{alg_3}) & $1E-12$ \\
    \hline
  \end{tabular}
  \caption{Advection equation: main simulation parameters of the VarPro method.}
  \label{tab_3}
\end{table}

We primarily consider a temporal dimension $t_f=10$ for the domain $\Omega$.
We employ the block time marching (BTM) scheme from~\cite{DongL2021}
together with the VarPro method for this problem; see
Remark~\ref{rem_4}.
We employ $10$ uniform time blocks in time.
and in each time block employ a uniform set of $Q=Q_1\times Q_1$
training collocation points with the VarPro method,
where $Q_1$ is varied systematically.
Following Section~\ref{sec:poisson},
we employ a much larger uniform set of $Q_2\times Q_2$
grid points within each time block
to evaluate the trained neural network
for the solution field and compute its errors by comparing
with the exact solution~\eqref{eq_28}.
We have also considered another spatial-temporal domain
with a much larger temporal dimension $t_f=100$. Correspondingly,
$100$ uniform time blocks are employed in simulations of this case.
The main simulation parameters for this problem are
summarized in Table~\ref{tab_3}.

\begin{figure}
  \centerline{
    \includegraphics[width=2.5in]{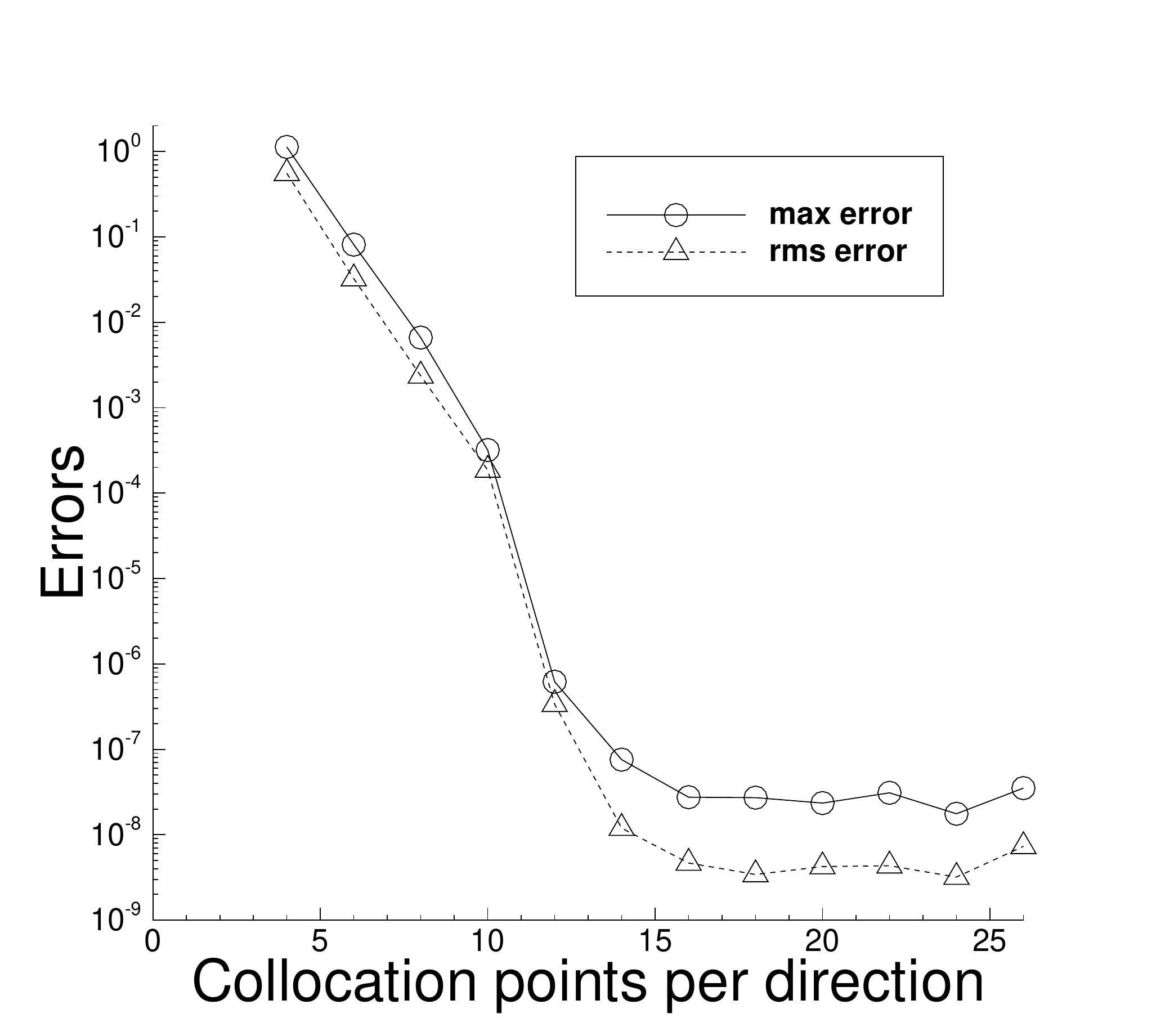}(a)
    \includegraphics[width=2.5in]{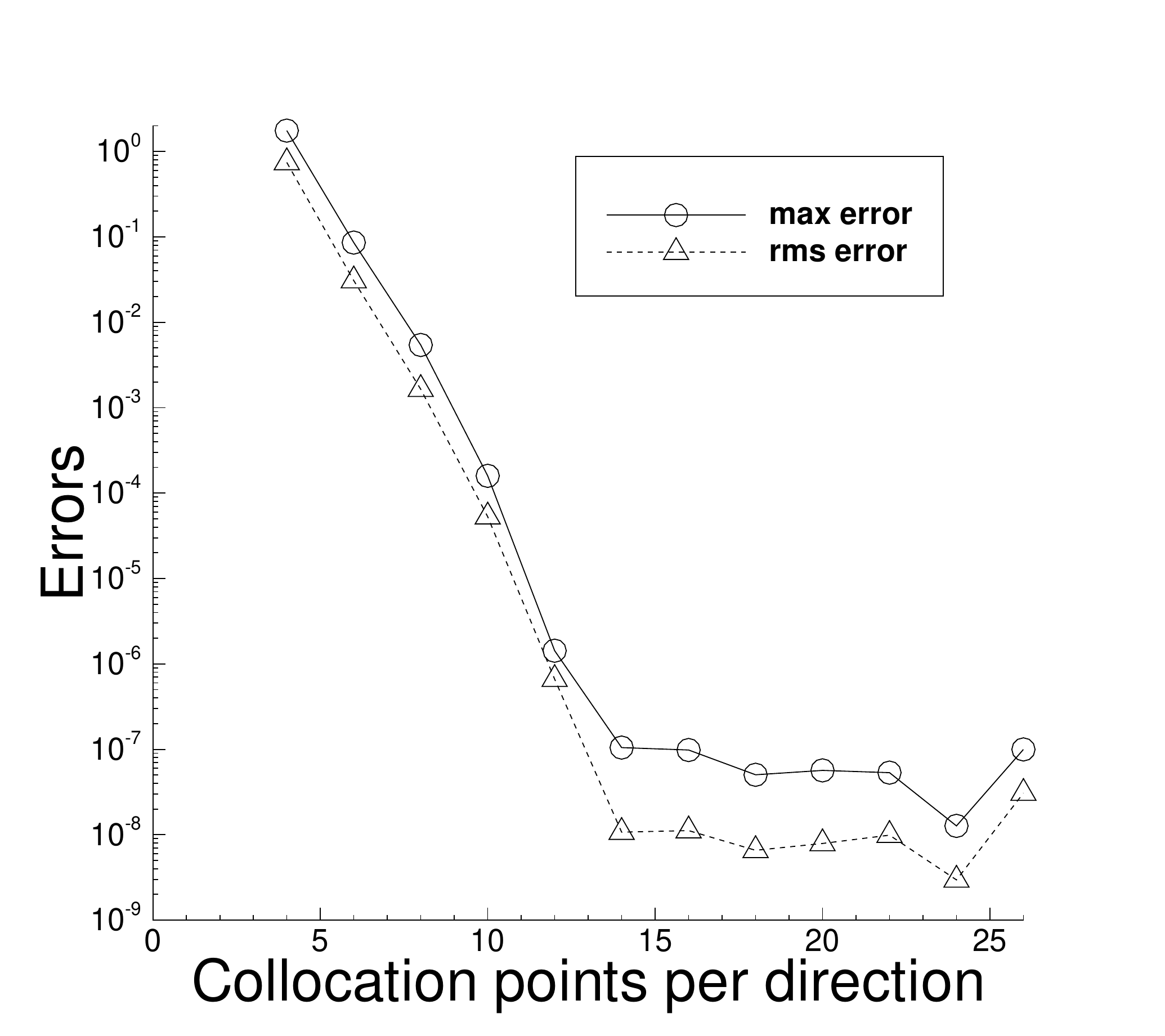}(b)
  }
  \caption{Advection equation: the maximum/rms errors of the
    VarPro solution versus the number of collocation points per direction in each
    time block,
    obtained with (a) the Gaussian and (b) the GELU activation functions.
    In (a,b), $t_f=10$, $10$ time blocks, neural network $[2, 100, 1]$,
    max-subiterations = 2 and $\delta=1.0$ in VarPro.
  }
  \label{fg_6}
\end{figure}

Let us first look into the VarPro errors obtained using neural networks
with one hidden layer. Figure \ref{fg_5}(b) shows the distribution
of the absolute error of the VarPro result in the spatial-temporal
domain. This result is for the temporal dimension $t_f=10$,
and is obtained using a neural network $[2, 100, 1]$ with
the Gaussian activation function and a uniform set of $Q=21\times 21$
training collocation points in the domain.
The VarPro result is highly accurate, with a maximum error on
the order $10^{-8}$ in the overall domain.

Figure \ref{fg_6} illustrates the convergence behavior of the VarPro
solution with respect to the number of collocation points
per direction ($Q_1$) in each time block.
This is for the temporal dimension $t_f=10$ computed with a neural
network $[2, 100, 1]$ and
$10$ time blocks . The plot (a) shows the maximum and rms errors
in the overall domain of the VarPro solution as a function
of $Q_1$ obtained with the Gaussian activation function.
The plot (b) shows the corresponding result obtained with
the GELU activation function.
The exponential decrease in the errors with increasing
number of collocation points (before saturation) is unmistakable.

\begin{figure}
  \centerline{
    \includegraphics[width=2.5in]{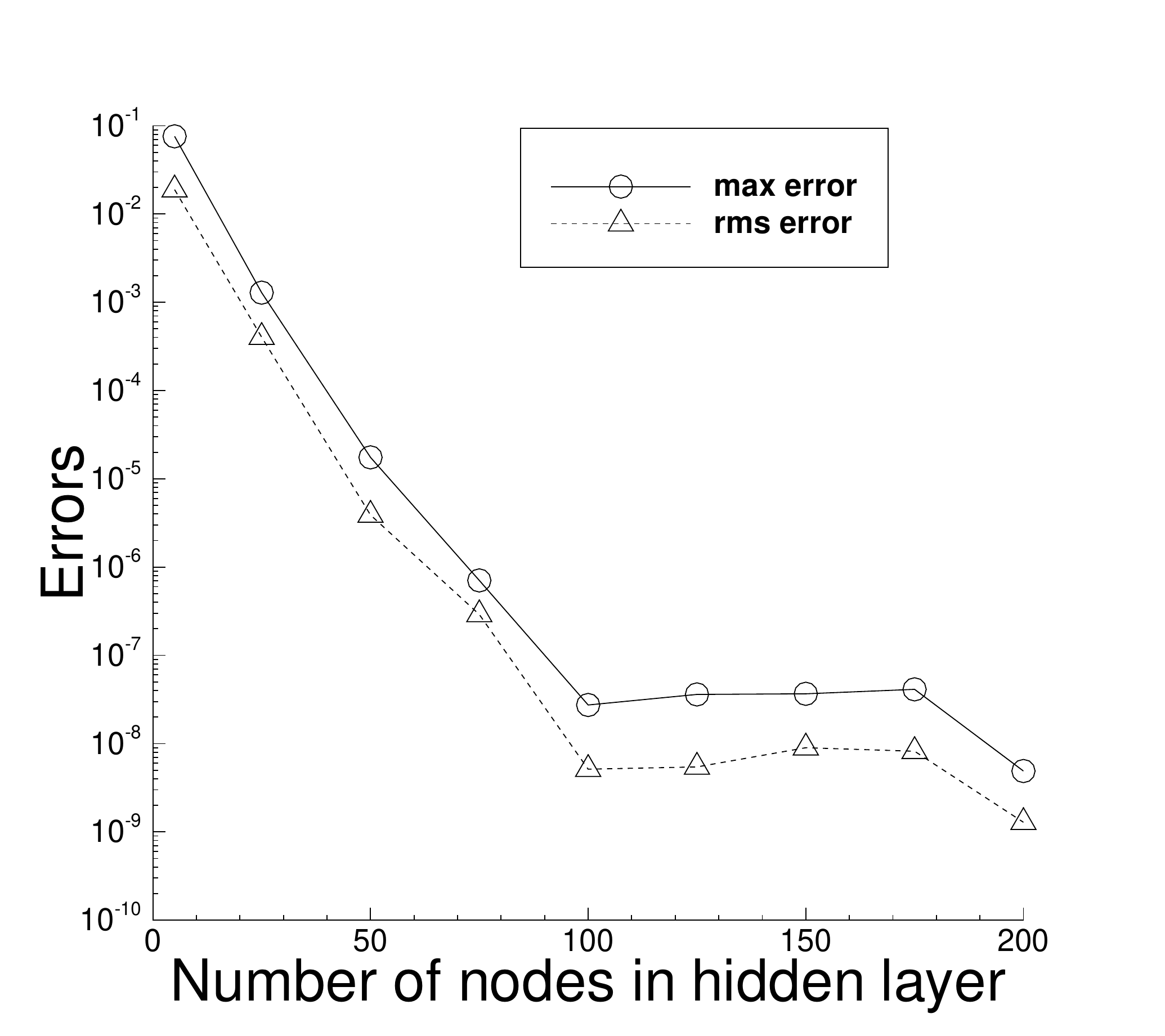}(a)
    \includegraphics[width=2.5in]{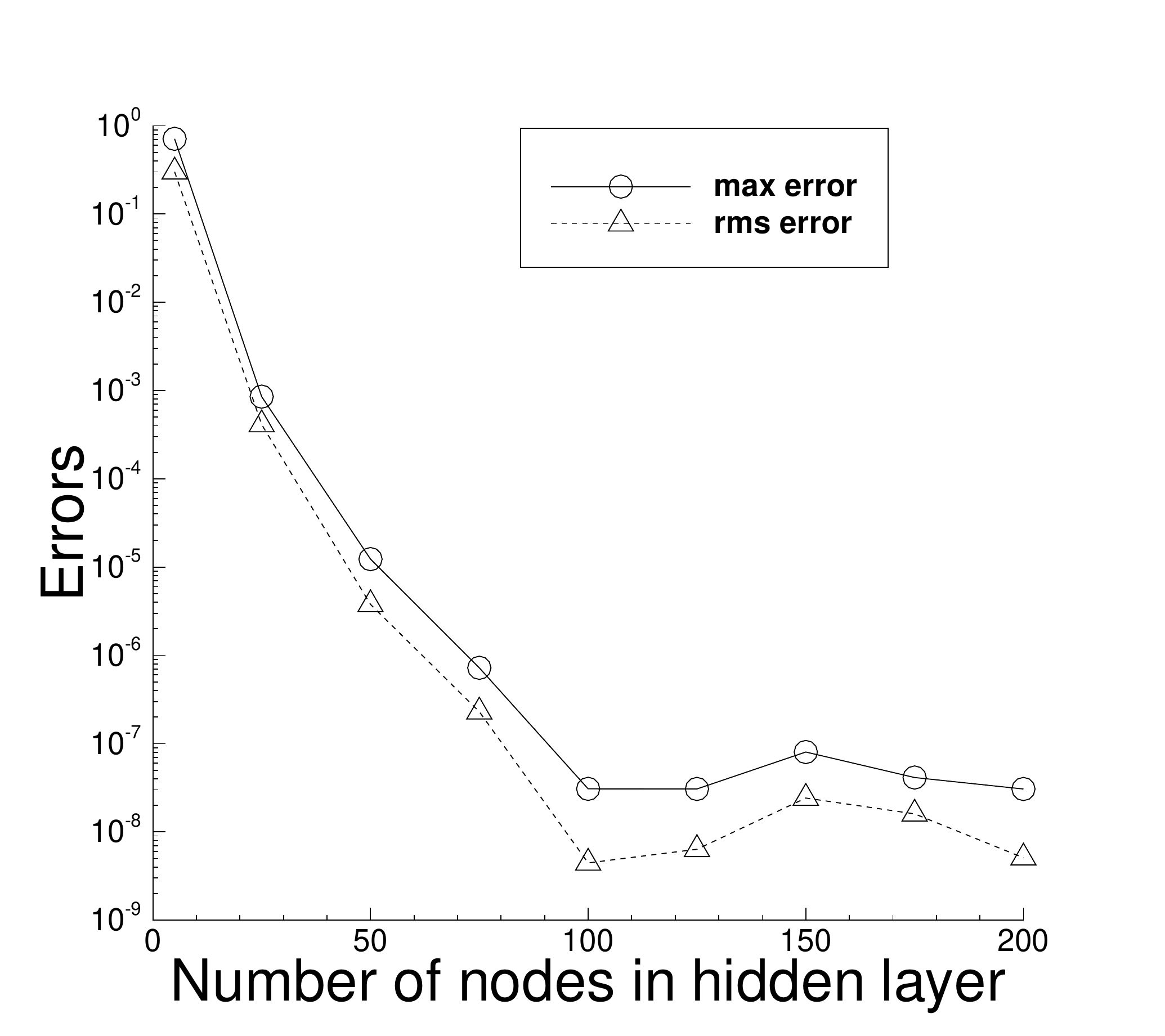}(b)
  }
  \caption{Advection equation: the maximum/rms errors of the VarPro solution
    versus the number of nodes in the hidden layer ($M$) obtained with
    (a) the Gaussian and (b) the GELU activation functions.
    In (a,b), $t_f=10$, $10$ time blocks, $Q=21\times 21$ uniform collocation points,
    neural network $[2, M, 1]$ with $M$ varied,
    max-subiterations$=2$ and $\delta=1.0$ in VarPro.
  }
  \label{fg_7}
\end{figure}

Figure~\ref{fg_7} illustrates the convergence behavior of the VarPro solution
with respect to the number of nodes in the hidden layer ($M$).
Here the domain corresponds to $t_f=10$, with $10$ time blocks and
a uniform set of $Q=21\times 21$ training collocation points per time block
in the VarPro simulation.
The neural network is given by $[2, M, 1]$, where $M$ is varied systematically.
Figures \ref{fg_7}(a) and (b) shows the maximum/rms errors in the overall
domain as a function of $M$ obtained using the Gaussian and the GELU activation
functions, respectively.
The exponential decrease in the errors with increasing $M$ (before saturation)
is evident.

\begin{figure}
  \centerline{
    \includegraphics[height=1.5in]{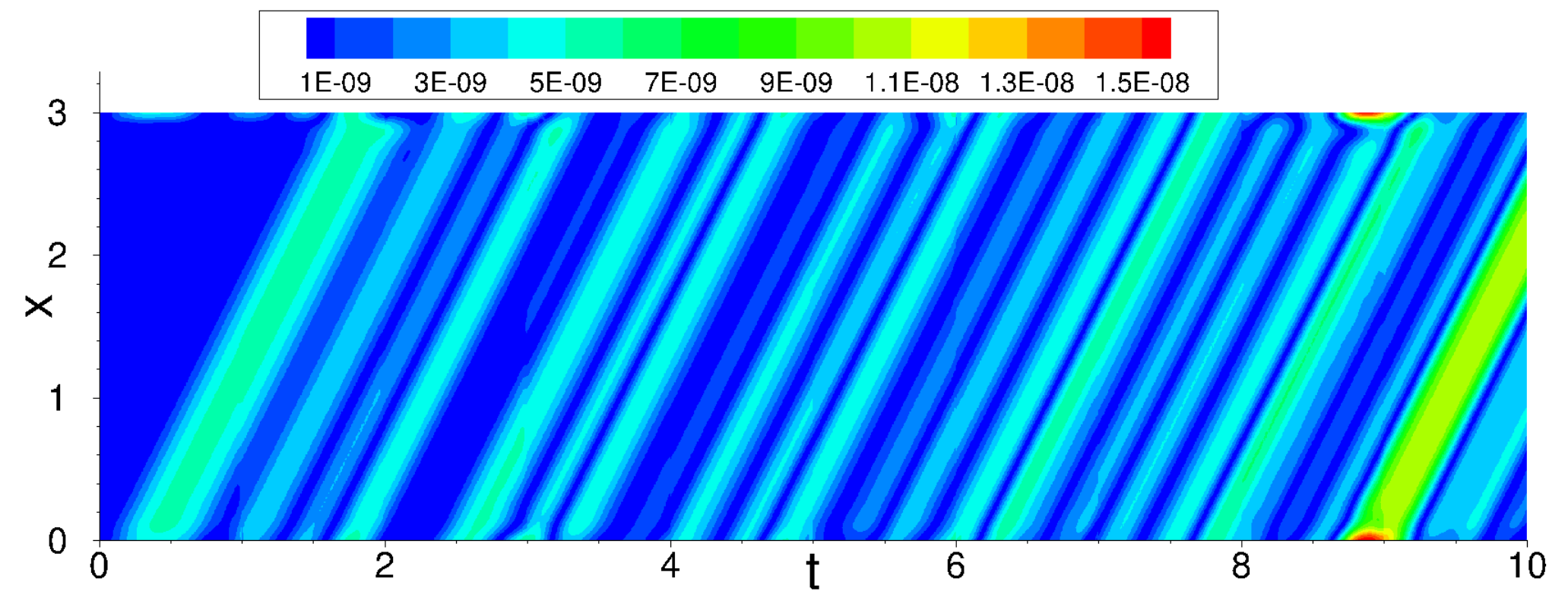}(a)
    \includegraphics[height=1.8in]{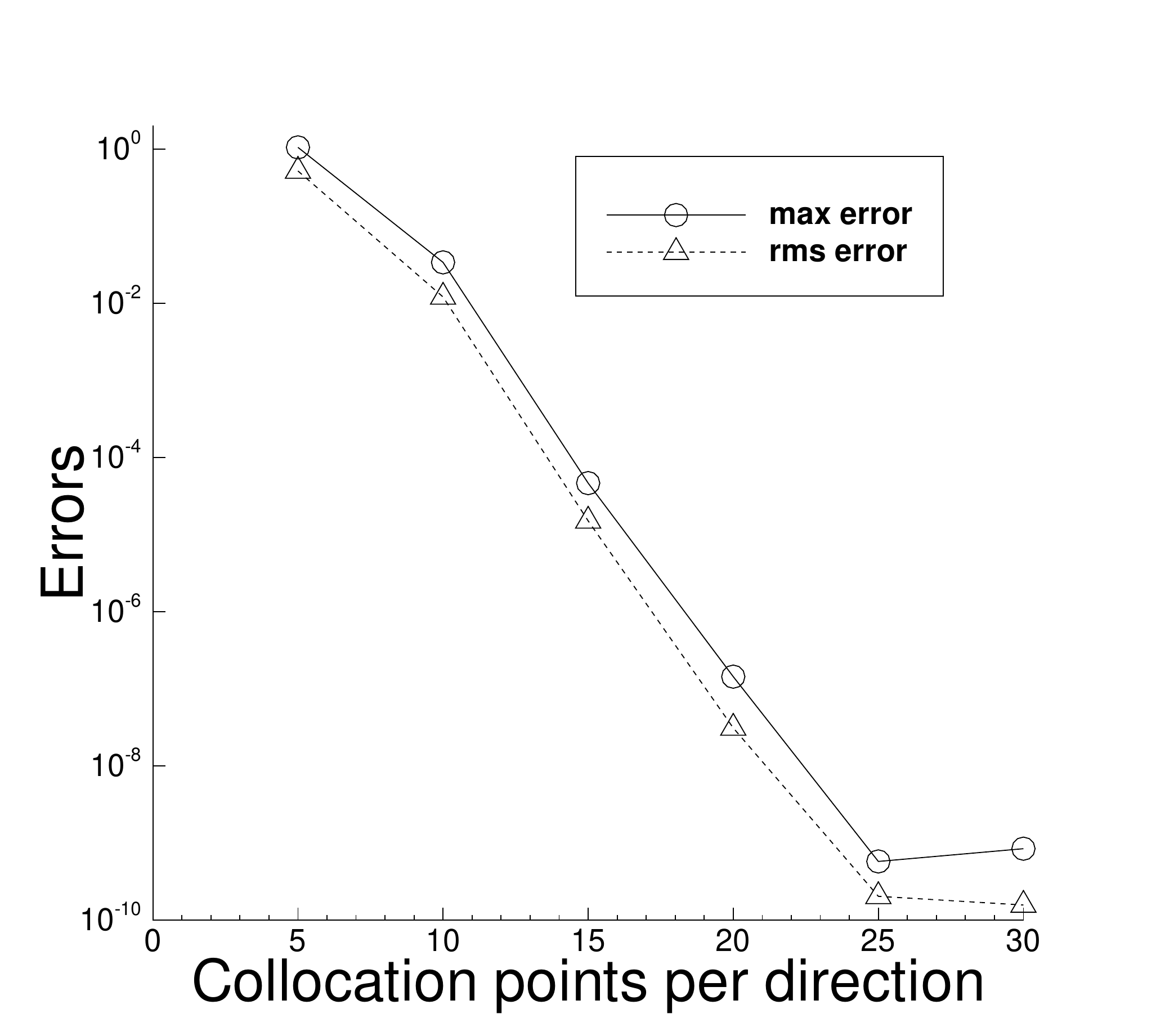}(b)
  }
  \caption{Advection equation (two hidden layers in NN):
    (a) Error distribution of the VarPro solution
    in the spatial-temporal plane, and
    (b) the maximum/rms errors of the VarPro solution versus
    the number of collocation points per direction, obtained with
    $2$ hidden layers in the neural network.
    In (a,b), $t_f=10$, $10$ uniform time blocks, neural network
    $[2, 10, 100, 1]$, Gaussian activation function,
    max-subiterations = 2 and $\delta=0.05$ in VarPro.
    $Q=21\times 21$ in (a) and is varied in (b).
  }
  \label{fg_8}
\end{figure}

Figure \ref{fg_8} illustrates the VarPro results computed using
a neural network containing two hidden layers.
The domain corresponds to $t_f=10$, and the neural network
has the architecture $[2, 10, 100, 1]$ with the Gaussian activation
function. Figure~\ref{fg_8}(a) shows the VarPro error distribution
in the overall spatial-temporal plane, obtained with a uniform
$Q=21\times 21$ training collocation points per time block.
Figure~\ref{fg_8}(b) depicts the maximum/rms VarPro errors
in the overall domain as a function of the collocation points
per direction in each time block, demonstrating the exponential
convergence behavior.


\begin{table}[tb]
  \centering
  \begin{tabular}{l | l| ll| ll}
    \hline
    $[-R_m,R_m]$ & collocation & VarPro & & ELM & \\ \cline{3-6}
     & points & max-error & rms-error & max-error & rms-error \\ \hline
    $R_m=1$ & $5\times 5$ & $2.505E-1$ & $1.005E-1$ & $2.505E-1$ & $1.005E-1$  \\
    & $10\times 10$ & $3.194E-4$ & $1.864E-4$ & $3.758E-4$ & $1.245E-4$  \\
    & $15\times 15$ & $4.254E-8$ & $6.954E-9$ & $4.232E-5$ & $1.269E-5$  \\
    & $20\times 20$ & $2.348E-8$ & $4.240E-9$ & $4.618E-5$ & $1.402E-5$ \\
    & $25\times 25$ & $1.493E-8$ & $2.889E-9$ & $5.471E-5$ & $1.598E-5$  \\
    & $30 \times 30$ & $1.354E-8$ & $2.528E-9$ & $6.419E-5$ & $1.773E-5$  \\ \hline
    $R_m=R_{m0}$ & $5\times 5$ & $1.038E-1$ & $3.252E-2$ & $1.038E-1$ & $3.252E-2$ \\
    $\qquad=0.7$ & $10\times 10$ & $1.996E-4$ & $9.437E-5$ & $2.240E-4$ & $6.519E-5$  \\
    & $15\times 15$ & $7.934E-8$ & $9.495E-9$ & $3.881E-5$ & $9.642E-6$ \\
    & $20\times 20$ & $9.261E-8$ & $3.253E-8$ & $3.471E-5$ & $1.060E-5$ \\
    & $25\times 25$ & $3.444E-8$ & $4.931E-9$ & $3.314E-5$ & $1.117E-5$ \\
    & $30\times 30$ & $4.670E-8$ & $7.428E-9$ & $3.283E-5$ & $1.164E-5$ \\
    \hline
  \end{tabular}
  \caption{Advection equation: comparison of the maximum/rms errors of the solutions
    obtained using the VarPro and ELM methods.
    $t_f=10$, $10$ time blocks,
    Neural network [2, 100, 1], Gaussian activation function.
    In VarPro, the max-subiterations is 2 and $\delta = 1.0$.
    In ELM/VarPro, the hidden-layer coefficients are
    set/initialized to uniform random values
    from $[-R_m,R_m]$, with $R_m=1.0$ or with $R_m=R_{m0}=0.7$.
  }
  \label{tab_4}
\end{table}

A comparison between the current VarPro method and
the ELM method~\cite{DongY2021,DongL2021} for solving the advection equation
is provided in Table~\ref{tab_4}.
The maximum and rms errors of the VarPro and the ELM methods obtained on
a series of training collocation points are listed.
These results are for the domain $t_f=10$, with $10$ time blocks
in block time marching. We have employed
a neural network $[2, 100, 1]$ with the Gaussian activation function.
The random hidden-layer coefficients in ELM and
the initial hidden-layer coefficients in VarPro are
generated by $R_m=1$ and $R_m=R_{m0}=0.7$.
It is evident that VarPro generally leads to significantly more accurate
results than ELM.


\begin{figure}
  \centering
  \includegraphics[width=6in]{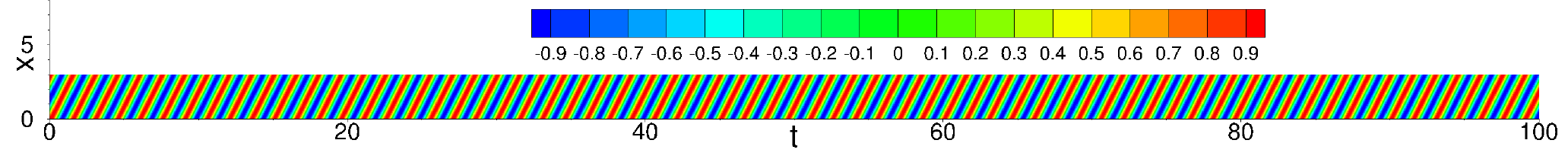}(a)
  \includegraphics[width=6in]{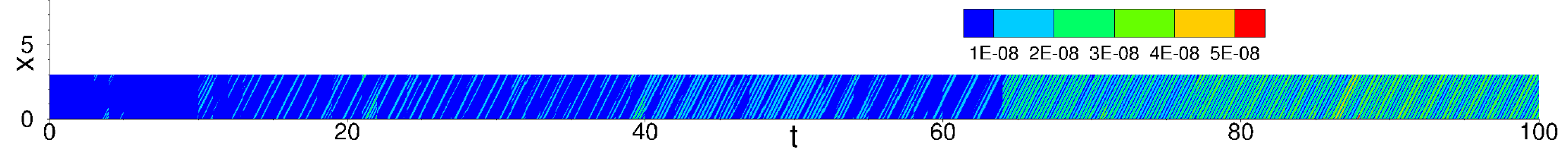}(b)
  \caption{Advection equation (long-time simulation): Distributions of (a) the VarPro
    solution and (b) its absolute error in the spatial-temporal domain.
    Domain: $(x,t)\in[0,3]\times[0,100]$, 100 uniform time blocks,
    $Q=25\times 25$ uniform collocation points per time block,
    neural network $[2,150,1]$, Gaussian activation function,
    no subiteration (max-subiterations=0) in VarPro.
  }
  \label{fg_9}
\end{figure}

\begin{figure}
  \centering
  \includegraphics[width=5in]{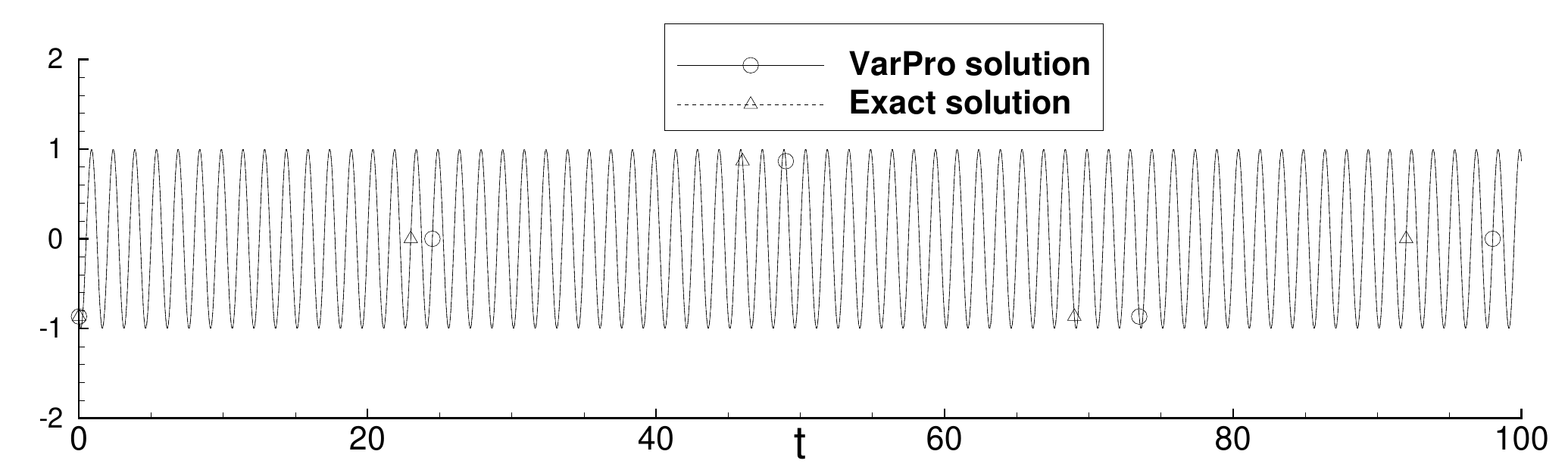}(a)
  \includegraphics[width=5in]{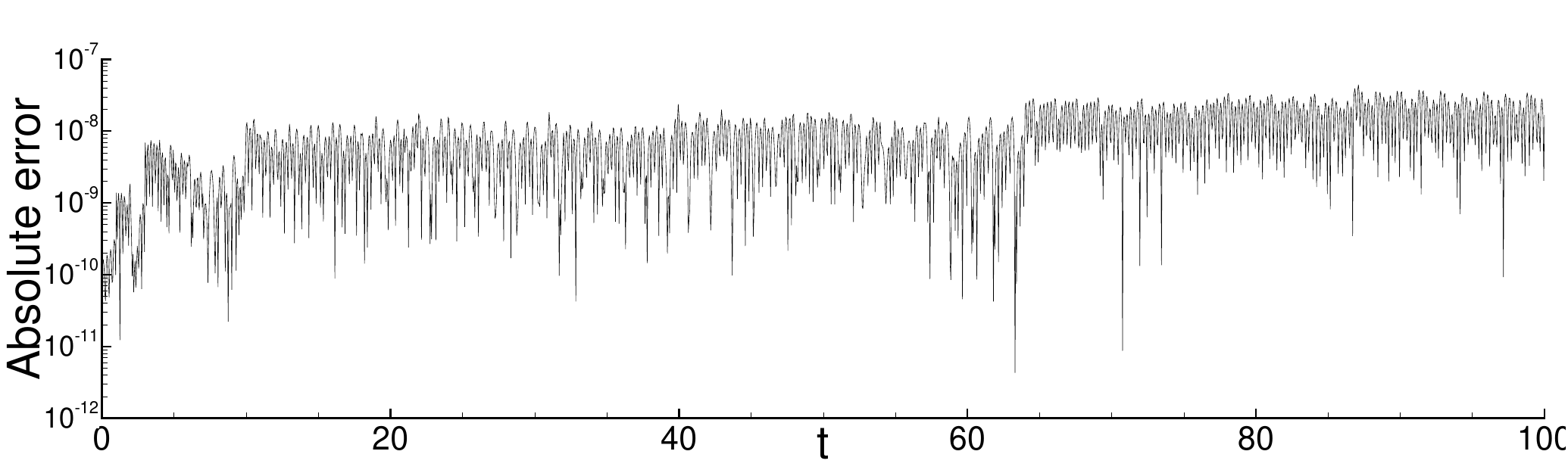}(b)
  \caption{Advection equation (long-time simulation): (a) Comparison of the time histories
    between the VarPro solution and the exact solution at the mid-point of the
    domain ($x=1.5$).
    (b) Time history of the absolute error of the VarPro solution at the mid-point
    of the domain ($x=1.5$).
    Simulation parameters and configurations here follow those of
    Figure~\ref{fg_9}.
  }
  \label{fg_10}
\end{figure}

Figures~\ref{fg_9} and~\ref{fg_10} illustrate a longer-time simulation of
the advection equation using the VarPro method and the block time
marching scheme.
Here the domain corresponds to $t_f=100$. We have employed
$100$ uniform time blocks, a set of $Q=25\times 25$ uniform collocation
points per time time block, and a neural network $[2, 150, 1]$
with the Gaussian activation function.
Figures~\ref{fg_9}(a) and (b) show the distributions of the VarPro solution
and its absolute errors in the overall spatial-temporal plane.
It can be observed that the VarPro method has produced
highly accurate results, with the maximum error on the order $10^{-8}$
in this long-time simulation.
Figure~\ref{fg_10}(a) compares the time histories of the VarPro solution
and the exact solution~\eqref{eq_28} at the mid-point of
the domain ($x=1.5$),
and Figure~\ref{fg_10}(b) shows the corresponding VarPro error
history  at this point.
These results indicate that the VarPro method together with
the block time marching scheme can produce highly accurate results
in long-time simulations.

\subsection{Nonlinear Examples}

\subsubsection{Nonlinear Helmholtz Equation}
\label{sec:nonl_helm}

\begin{figure}
  \centerline{
    \includegraphics[width=2.5in]{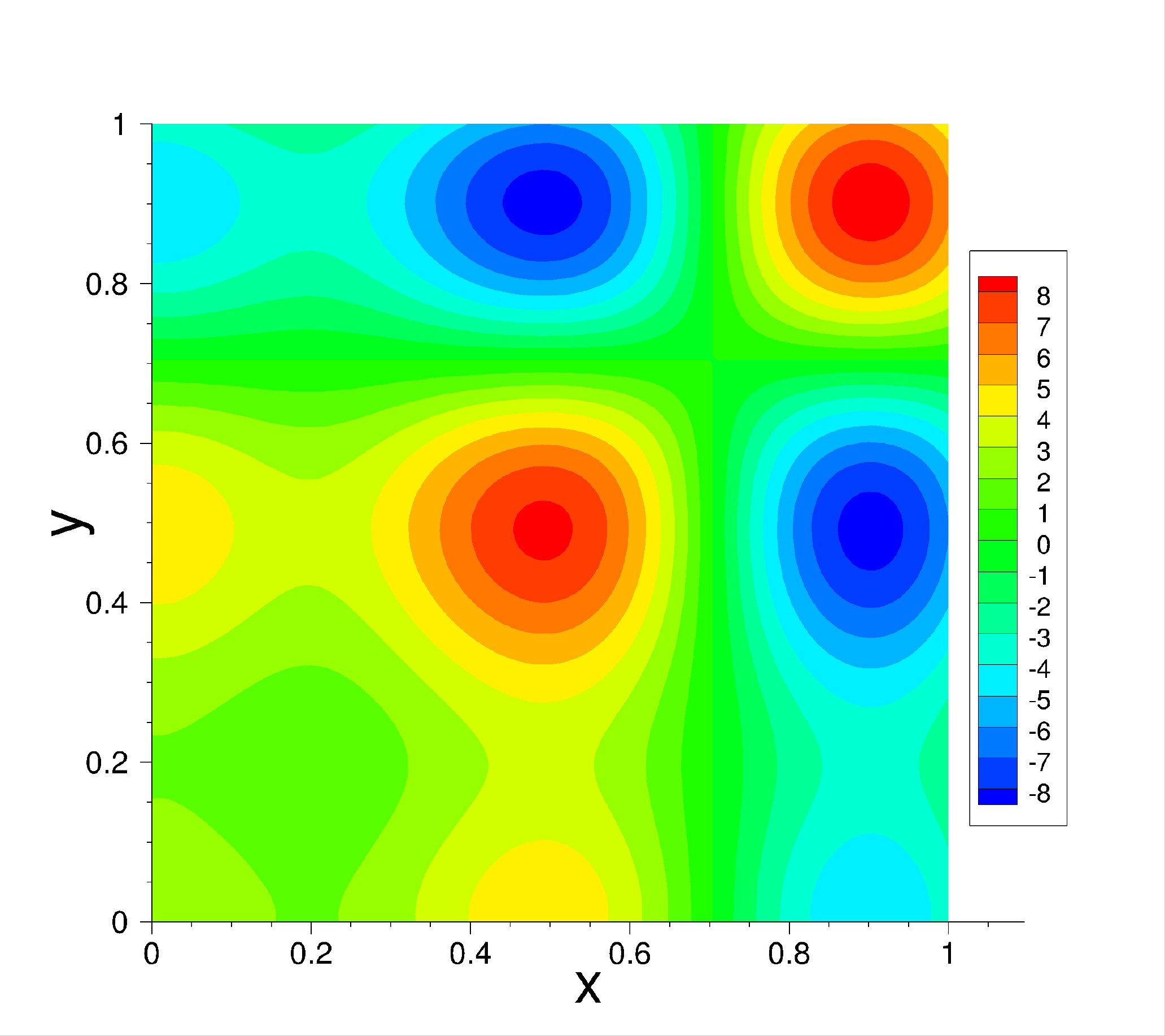}(a)
    \includegraphics[width=2.5in]{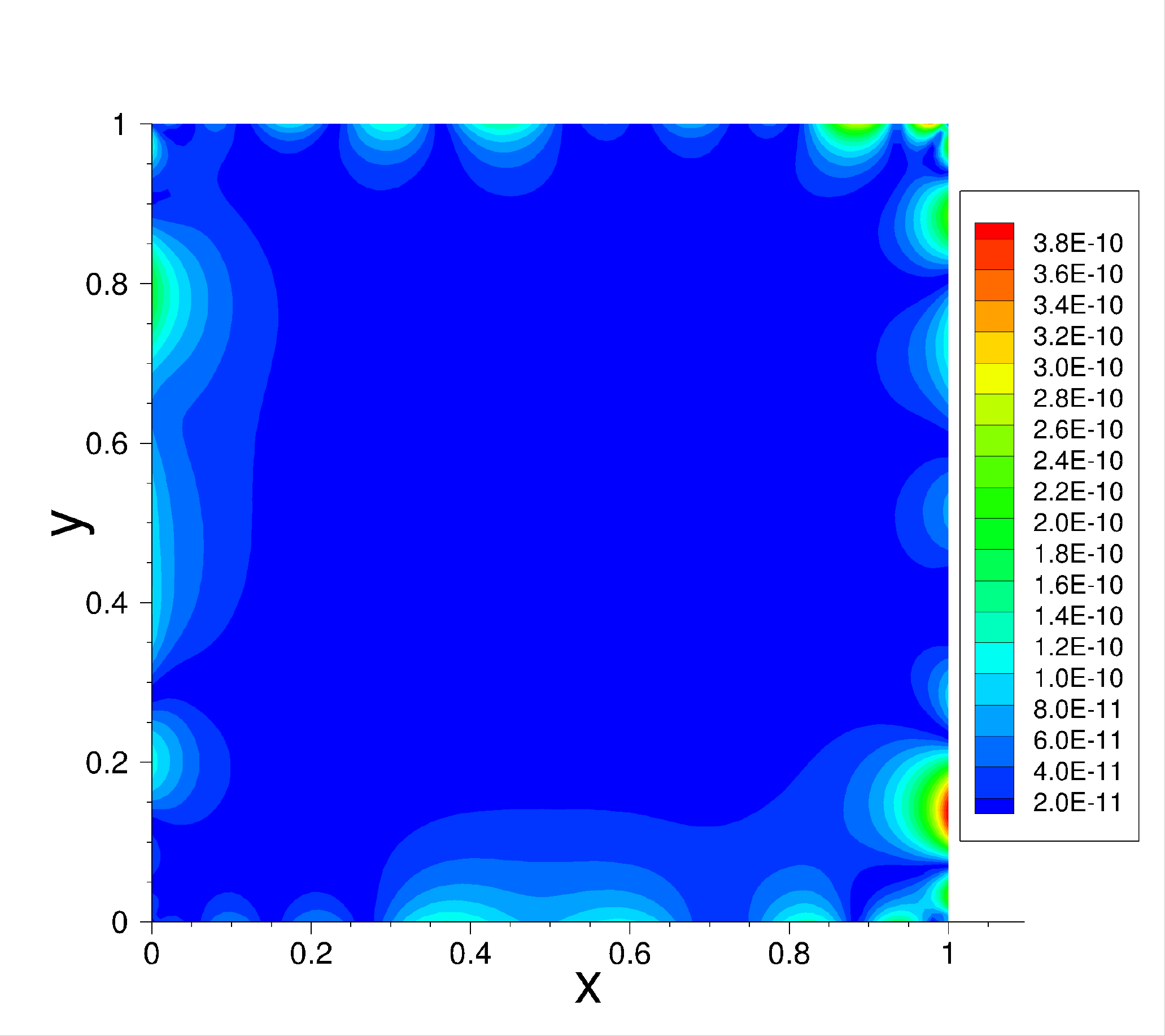}(b)
  }
  \caption{Nonlinear Helmholtz equation: Distributions of
    (a) the exact solution and (b) the absolute error of the VarPro
    solution.
    In (b), neural network [2, 200, 1], ``$\sin$'' activation function,
    $Q=21\times 21$ uniform collocation points, $\delta = 0.05$ in VarPro.
  }
  \label{fg_11}
\end{figure}

In the first nonlinear example, we consider
the boundary value problem
with a nonlinear Helmholtz equation on the unit square domain $[0,1]\times[0,1]$,
\begin{subequations}
  \begin{align}
    &
    \frac{\partial^2u}{\partial x^2} + \frac{\partial^2u}{\partial y^2} - 100u
    + 5\cos(2u) = f(x,y), \\
    &
    u(x,0) = g_1(x), \quad
    u(x,1) = g_2(x), \quad
    u(0,y) = g_3(y), \quad
    u(1,y) = g_4(y),
  \end{align}
\end{subequations}
where $u(x,y)$ is the field function to be solved for,
$f(x,y)$ is a prescribed source term, and $g_i$ ($1\leqslant i\leqslant 4$)
denote the boundary data.
With $f$ and $g_i$ ($1\leqslant i\leqslant 4$) chosen appropriately,
this problem admits the following analytic solution,
\begin{multline}
  u = \left[\frac52\cos\left(\frac32\pi x - \frac25\pi \right)
  + \frac32\cos\left(3\pi x + \frac{3\pi}{10} \right)
  +\frac12(e^{x} - e^{-x})\right]
  \left[\frac52\cos\left(\frac32\pi y - \frac25\pi \right) \right. \\
    \left.
  + \frac32\cos\left(3\pi y + \frac{3\pi}{10} \right)
  +\frac12(e^{y} - e^{-y})\right].
\end{multline}
We employ this analytic solution in the following tests.
Figure~\ref{fg_11}(a) shows the distribution of this analytic
solution in the $xy$ plane.

\begin{table}[tb]
  \centering
  \begin{tabular}{ll | ll}
    \hline
    parameter & value & parameter & value \\ \hline
    neural network & $[2, M, 1]$, or $[2, 5, M, 1]$ & training points $Q$ & $Q_1\times Q_1$ \\
    $M$ &  varied & $Q_1$ & varied \\
    activation function & $\sin$, Gaussian & testing points & $Q_2\times Q_2$ \\
    random seed & $1$ & $Q_2$ & 101 \\
    initial guess $\bm\theta_0$ & random values on $[-R_m,R_m]$ & $R_m$ & $1.0$ \\
    $\delta$ (Algorithm~\ref{alg_3}) &  $0.02$, $0.05$, $0.1$, or $0.2$
    & $p$ (Algorithm~\ref{alg_3}) & $0.5$ \\
    max-subiterations & $2$ & threshold (Algorithm~\ref{alg_3}) & $1E-12$ \\
    max-iterations-newton & $20$ & tolerance-newton & $1E-8$ \\
    \hline
  \end{tabular}
  \caption{Nonlinear Helmholtz equation:
    main simulation parameters of the VarPro method.}
  \label{tab_5}
\end{table}

We employ neural networks with the architectures
$[2, M, 1]$ and $[2, 5, M, 1]$ in the VarPro simulations,
where $M$ is varied in the tests.
The sine function, $\sigma(x)=\sin(x)$, or
the Gaussian function, $\sigma(x)=e^{-x^2}$,
is employed as the activation functions for the hidden nodes.
A uniform set of $Q=Q_1\times Q_1$
training collocation points, where $Q_1$ is varied,
is used to train the neural network.
The VarPro solution is computed on a larger set of $Q_2\times Q_2$
(with $Q_2=101$)
uniform grid points by evaluating the trained neural network,
and compared with the analytic solution to compute its errors.
Table~\ref{tab_5} provides the main simulation parameters
for this problem and the VarPro method.
In this table ``max-iterations-newton'' denotes the maximum number of
Newton iterations, and ``tolerance-newton'' denotes
the relative tolerance for the Newton iteration (see lines
$3$ and $9$ of Algorithm~\ref{alg_4}).

Figure \ref{fg_11}(b) illustrates the error distribution of
a VarPro solution in the $xy$ plane,
computed using a neural network $[2, 200, 1]$ with the $\sin$
activation function and a uniform set of $Q=21\times 21$
collocation points.
The result is observed to be highly accurate, with
a maximum error on the order $10^{-10}$ in the domain.

\begin{figure}
  \centerline{
    \includegraphics[width=2.5in]{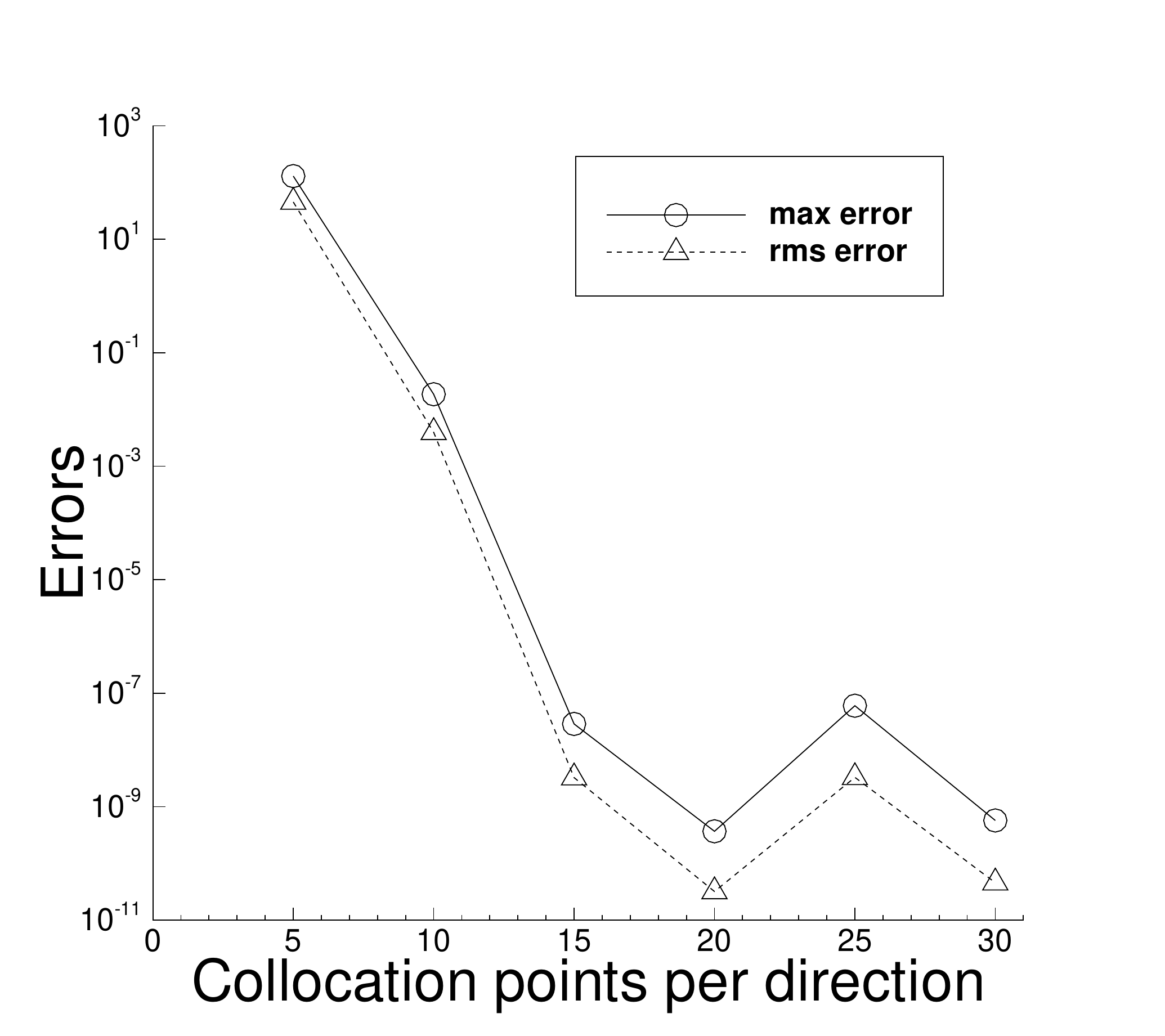}(a)
    \includegraphics[width=2.5in]{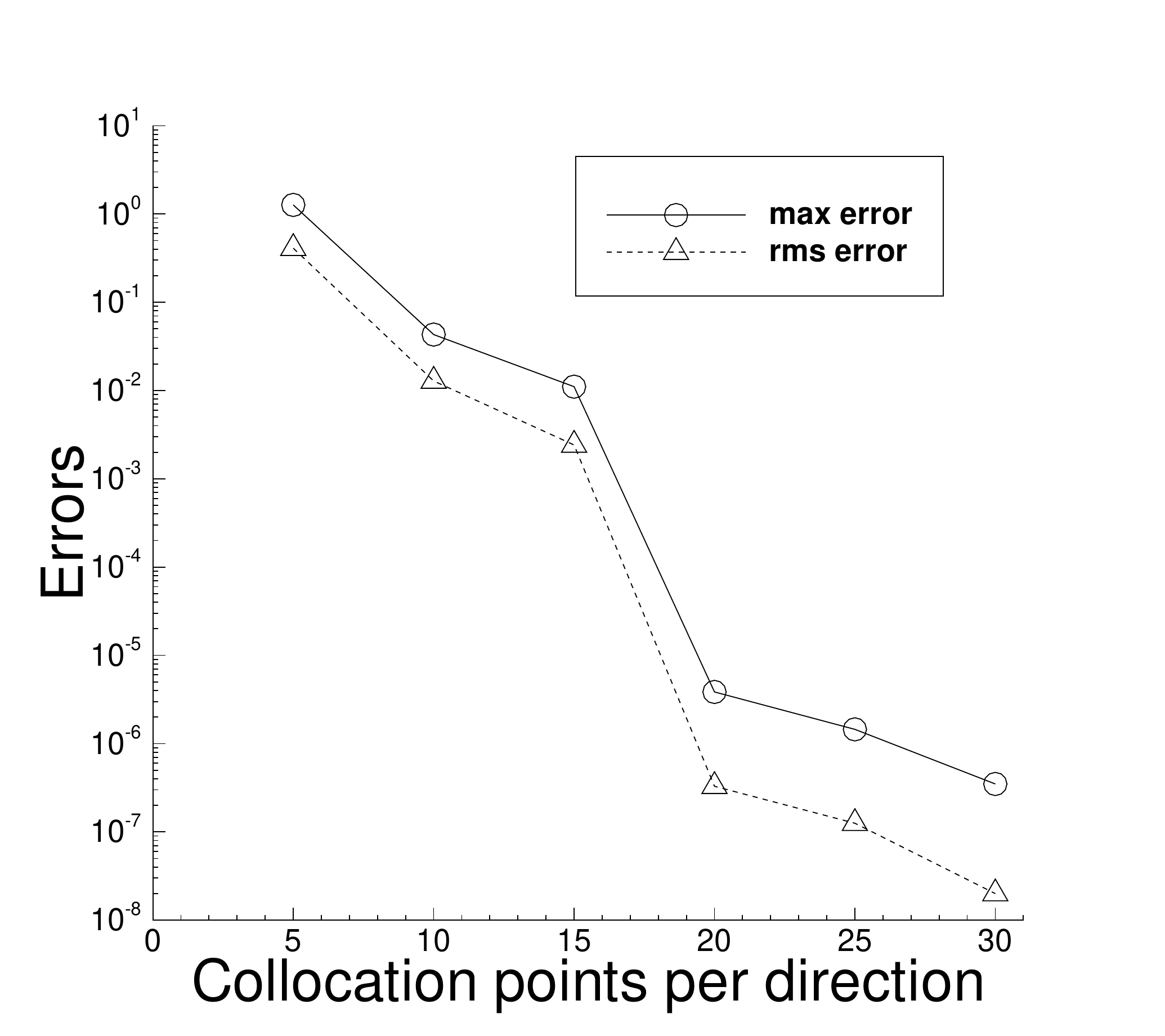}(b)
  }
  \caption{Nonlinear Helmholtz equation: the maximum/rms errors
    of the VarPro solution versus the number of collocation points
    per direction, obtained using (a) the sine and (b) the Gaussian
    activation functions.
    Neural network [2, 200, 1],
    $\delta=0.1$ in (a) and $\delta=0.2$ in (b) with VarPro.
  }
  \label{fg_12}
\end{figure}

Figure~\ref{fg_12} illustrates the convergence behavior of the VarPro method
with respect to the number of training collocation points in the domain.
In these tests the number of collocation points per direction ($Q_1$)
is varied systematically.
The two plots show the maximum/rms errors in the domain of the VarPro
solution as a function of $Q_1$, obtained using
the $\sin$ (plot (a)) and the Gaussian (plot (b)) activation functions.
These VarPro results are attained using a neural network $[2, 200, 1]$.
We observe an exponential decrease in the VarPro errors (before saturation)
with increasing number of collocation points.

\begin{figure}
  \centerline{
    \includegraphics[width=2.5in]{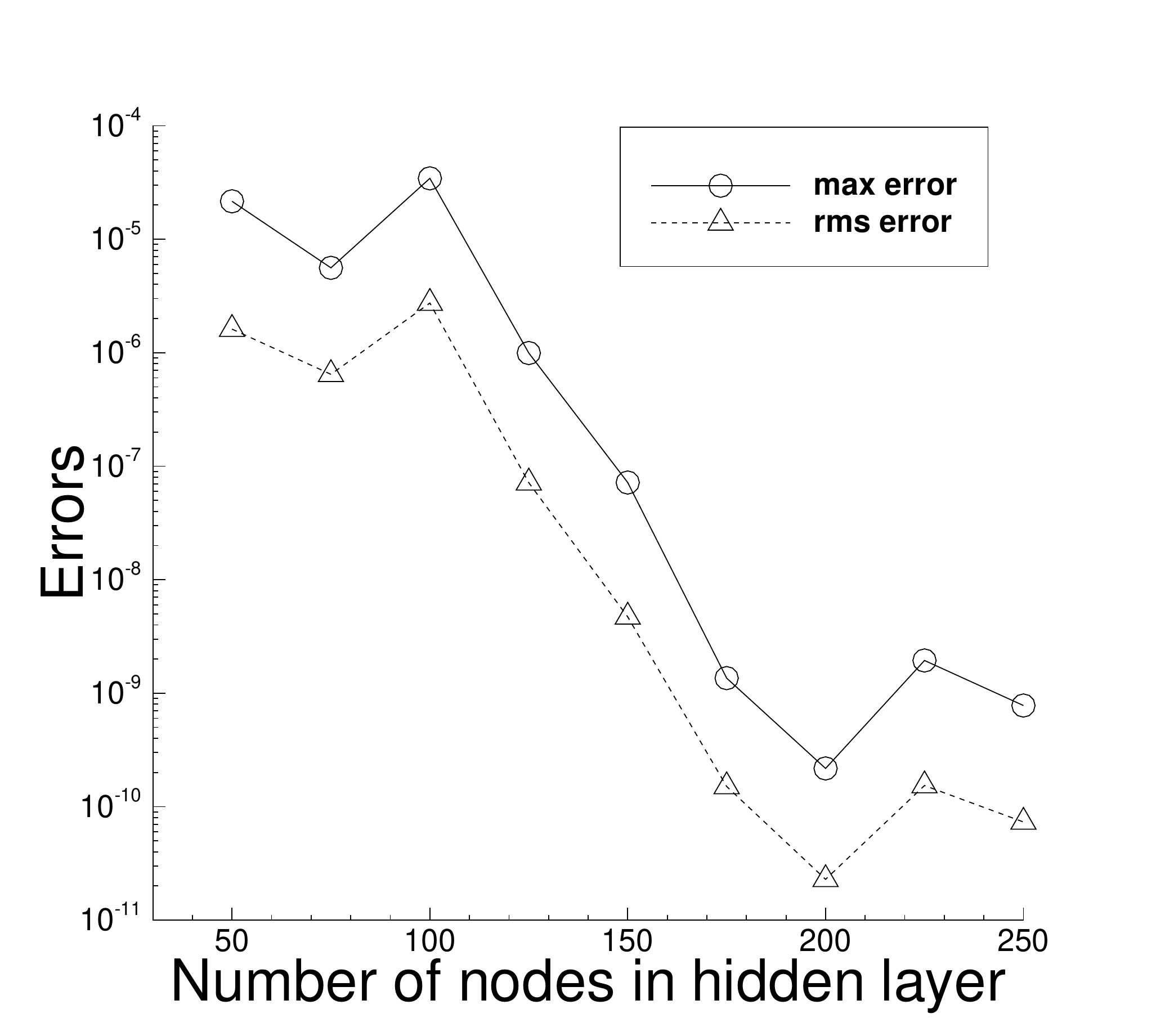}(a)
    \includegraphics[width=2.5in]{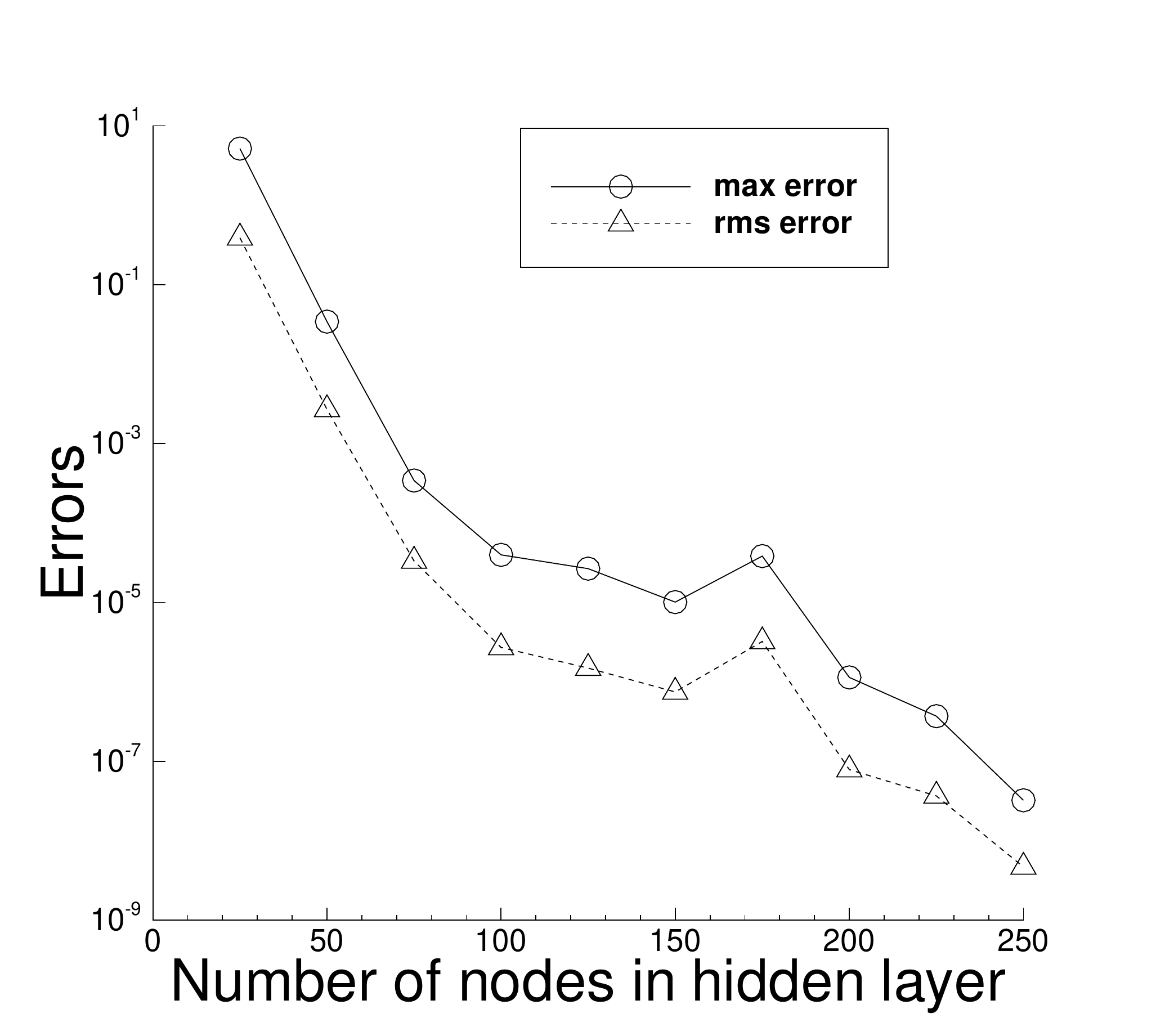}(b)
  }
  \caption{Nonlinear Helmholtz equation: the maximum/rms errors of
    the VarPro solution versus the number of nodes in the hidden layer,
    obtained with (a) the sine and (b) the Gaussian activation functions.
    Neural network $[2, M, 1]$, where $M$ is varied,
    $Q=21\times 21$, $\delta=0.1$ in (a) and $\delta=0.2$ in (b)
    with VarPro.
  }
  \label{fg_13}
\end{figure}

Figure \ref{fg_13} illustrates the convergence behavior of the VarPro solution
with respect to the number of nodes in the hidden layer ($M$) for the nonlinear
Helmholtz equation.
Here the neural network has an architecture $[2, M, 1]$, where $M$ is
varied systematically, and a uniform set of $Q=21\times 21$ training collocation
points is employed in the simulation.
This figure shows the maximum/rms errors in the domain as a function of $M$,
obtained using the sin (plot (a)) and the Gaussian (plot (b)) activation functions.
The errors computed with the $\sin$ activation function appear not quite regular
as $M$ increases.
But overall all these errors appear to decrease approximately
exponentially with increasing $M$. 

\begin{figure}
  \centerline{
    \includegraphics[width=2.5in]{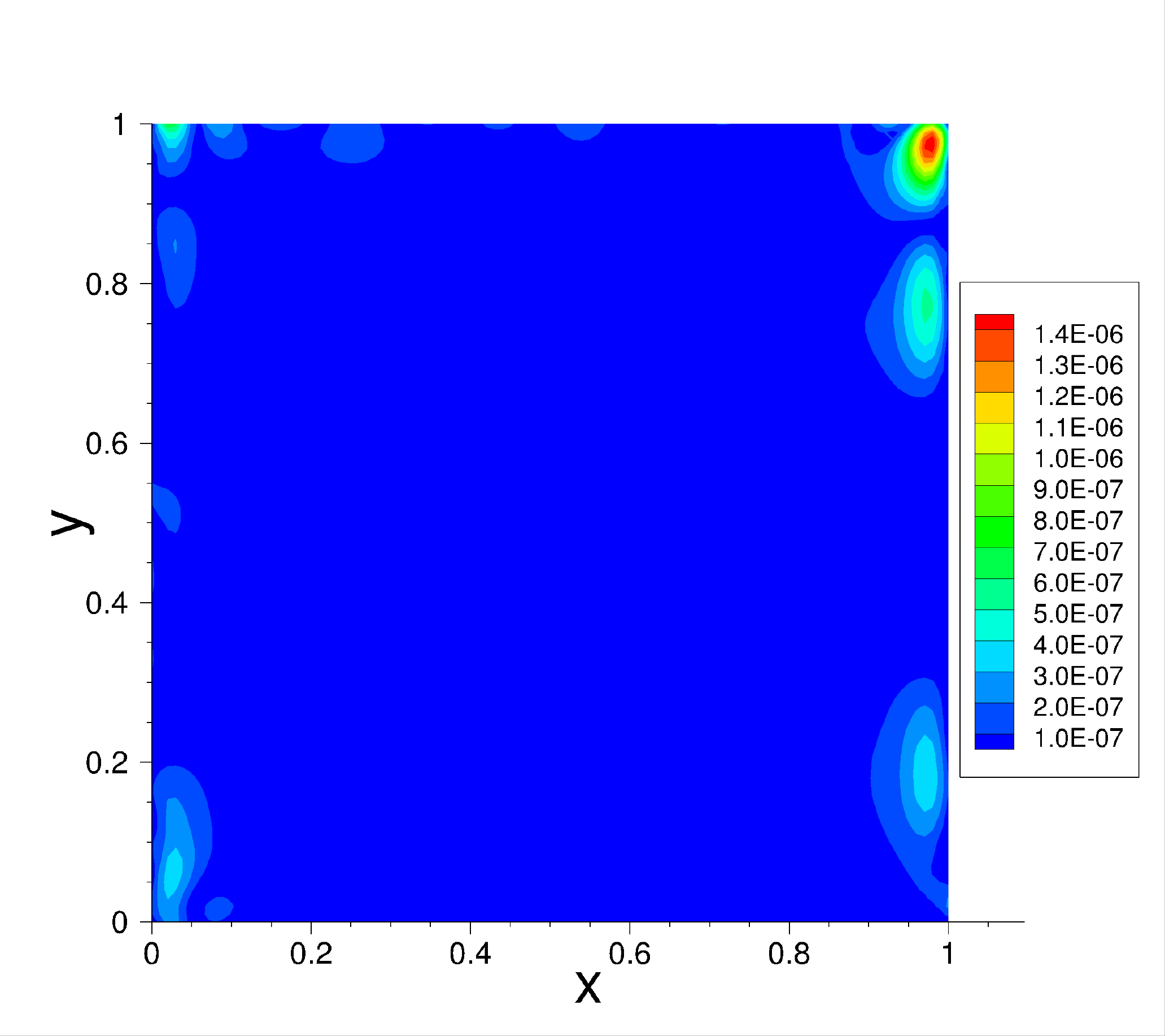}(a)
    \includegraphics[width=2.5in]{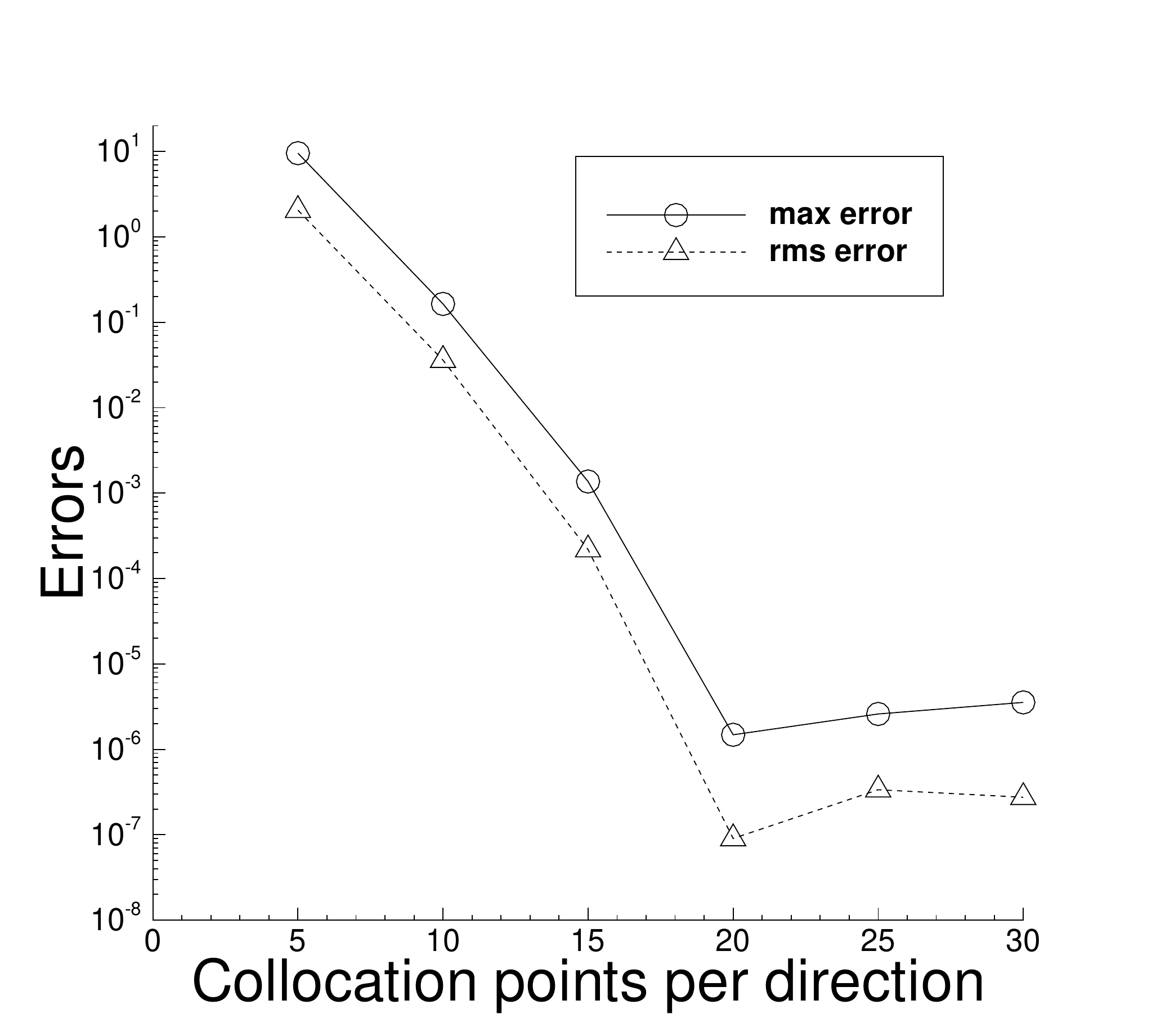}(b)
  }
  \caption{Nonlinear Helmholtz equation (two hidden layers in NN):
    (a) Error distribution of the VarPro
    solution. (b) The VarPro maximum/rms errors versus
    the number of collocation points per direction ($Q_1$).
    Neural network [2, 5, 200, 1], $\sin$ activation function,
    $Q=20\times 20$ in (a) and is varied in (b), $\delta=0.02$ in (a,b) with VarPro.
  }
  \label{fg_14}
\end{figure}

Figure~\ref{fg_14} illustrates the VarPro results obtained using two
hidden layers in the neural network.
Here we consider a neural network with the architecture $[2, 5, 200, 1]$,
with the $\sin$ activation function.
Figure \ref{fg_14}(a) shows the error distribution of the VarPro solution
obtained using $Q=20\times 20$ training collocation points.
In Figure \ref{fg_14}(b) the number of collocation points per direction ($Q_1$)
is varied systematically, and the maximum/rms errors are plotted as
a function of $Q_1$. An exponential decrease in the errors (before saturation)
can be observed. 


\begin{table}[tb]
  \centering
  \begin{tabular}{l | l| ll| ll}
    \hline
    $[-R_m,R_m]$ & collocation & VarPro & & ELM & \\ \cline{3-6}
     & points & max-error & rms-error & max-error & rms-error \\ \hline
    $R_m=1$ & $5\times 5$ & $1.297E+2$ & $4.536E+1$ & $4.388E+0$ & $8.195E-1$  \\
    & $10\times 10$ & $1.855E-2$ & $3.955E-3$ & $7.701E+0$ & $1.210E+0$  \\
    & $15\times 15$ & $2.868E-8$ & $3.252E-9$ & $3.743E-1$ & $4.767E-2$ \\
    & $20\times 20$ & $3.679E-10$ & $3.203E-11$ & $1.280E+0$ & $1.864E-1$ \\
    & $25\times 25$ & $6.014E-8$ & $3.281E-9$ & $1.434E+0$ & $2.198E-1$  \\
    & $30 \times 30$ & $5.709E-10$ & $4.576E-11$ & $8.365E-1$ & $1.123E-1$  \\ \hline
    $R_m=R_{m0}=4.4$ & $5\times 5$ & $6.742E-1$ & $2.408E-1$ & $1.182E-1$ & $3.589E-2$  \\
    & $10\times 10$ & $5.869E-3$ & $1.169E-3$ & $8.987E-7$ & $1.861E-7$  \\
    & $15\times 15$ & $7.130E-9$ & $1.024E-9$ & $6.690E-9$ & $8.655E-10$ \\
    & $20\times 20$ & $2.851E-9$ & $2.364E-10$ & $2.384E-8$ & $2.835E-9$ \\
    & $25\times 25$ & $4.193E-10$ & $5.173E-11$ & $3.133E-8$ & $3.611E-9$  \\
    & $30\times 30$ & $1.128E-9$ & $1.029E-10$ & $3.813E-8$ & $4.194E-9$  \\
    \hline
  \end{tabular}
  \caption{Nonlinear Helmholtz equation:
    comparison of the maximum/rms errors of the VarPro and ELM solutions.
    Neural network [2, 200, 1], $\sin$ activation function.
    In VarPro, $\delta=0.1$, and the tolerance-newton is set to $1E-8$ with $R_m=1$
    and to $1E-14$ with $R_m=4.4$.
  }
  \label{tab_6}
\end{table}

Table~\ref{tab_6} is a comparison of the errors between the
VarPro method and the ELM method for solving the nonlinear
Helmholtz equation.
These results are for a neural network $[2, 200, 1]$
with the $\sin$ activation function.
In ELM the random hidden-layer coefficients are set to,
and in VarPro the hidden-layer coefficients are initialized
to, uniform random values from $[-R_m,R_m]$, where
$R_m=1.0$ or $R_m=R_{m0}=4.4$.
Several sets of uniform training collocation points
 are tested, ranging from $Q=5\times 5$ to $Q=30\times 30$.
The VarPro results are in general markedly more accurate than those of
the ELM results. This is especially pronounced
for those cases corresponding to $R_m=1.0$.

\subsubsection{Viscous Burgers' Equation}
\label{sec:burgers}

\begin{figure}
  \centerline{
    \includegraphics[width=2.5in]{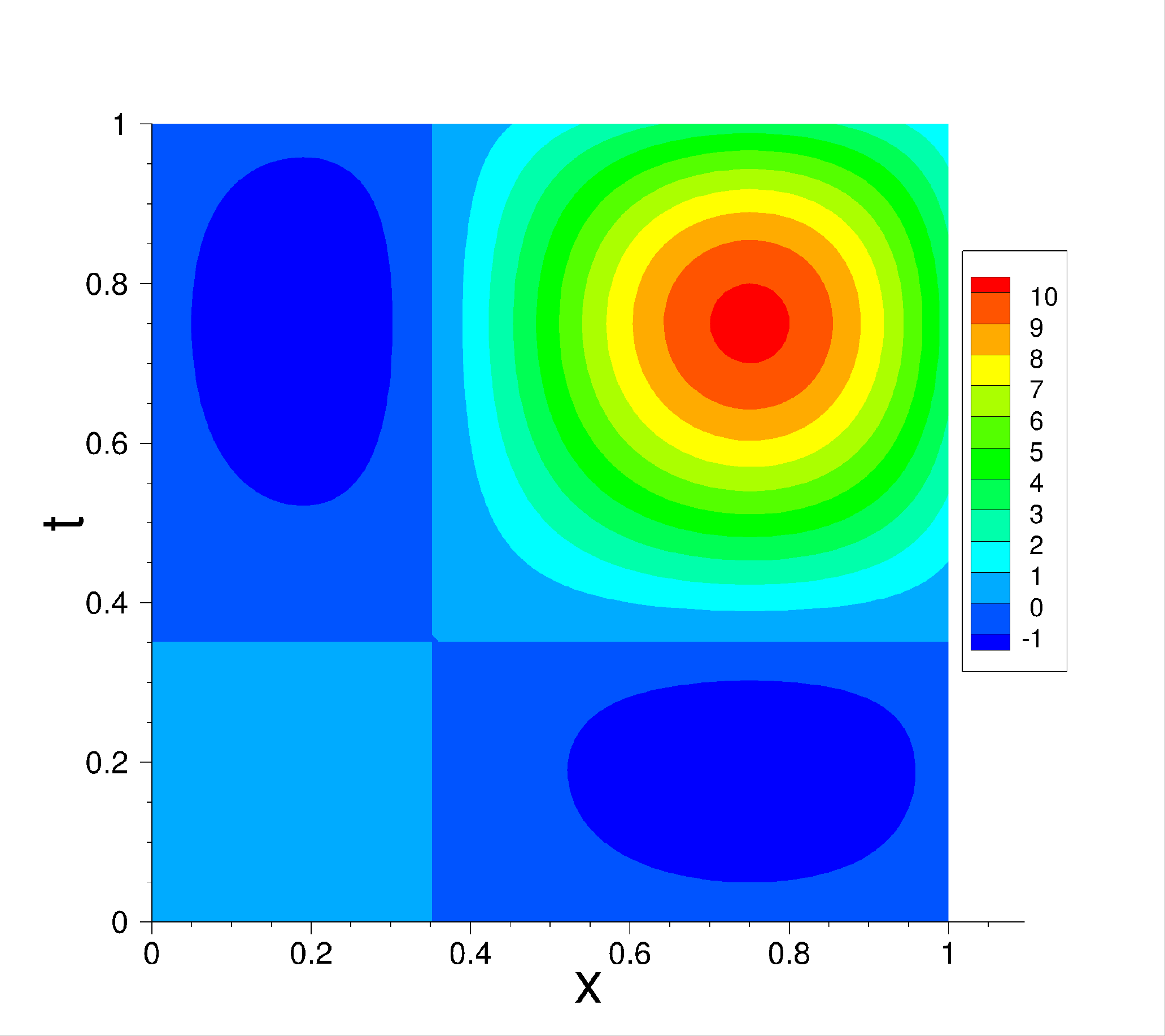}(a)
    \includegraphics[width=2.5in]{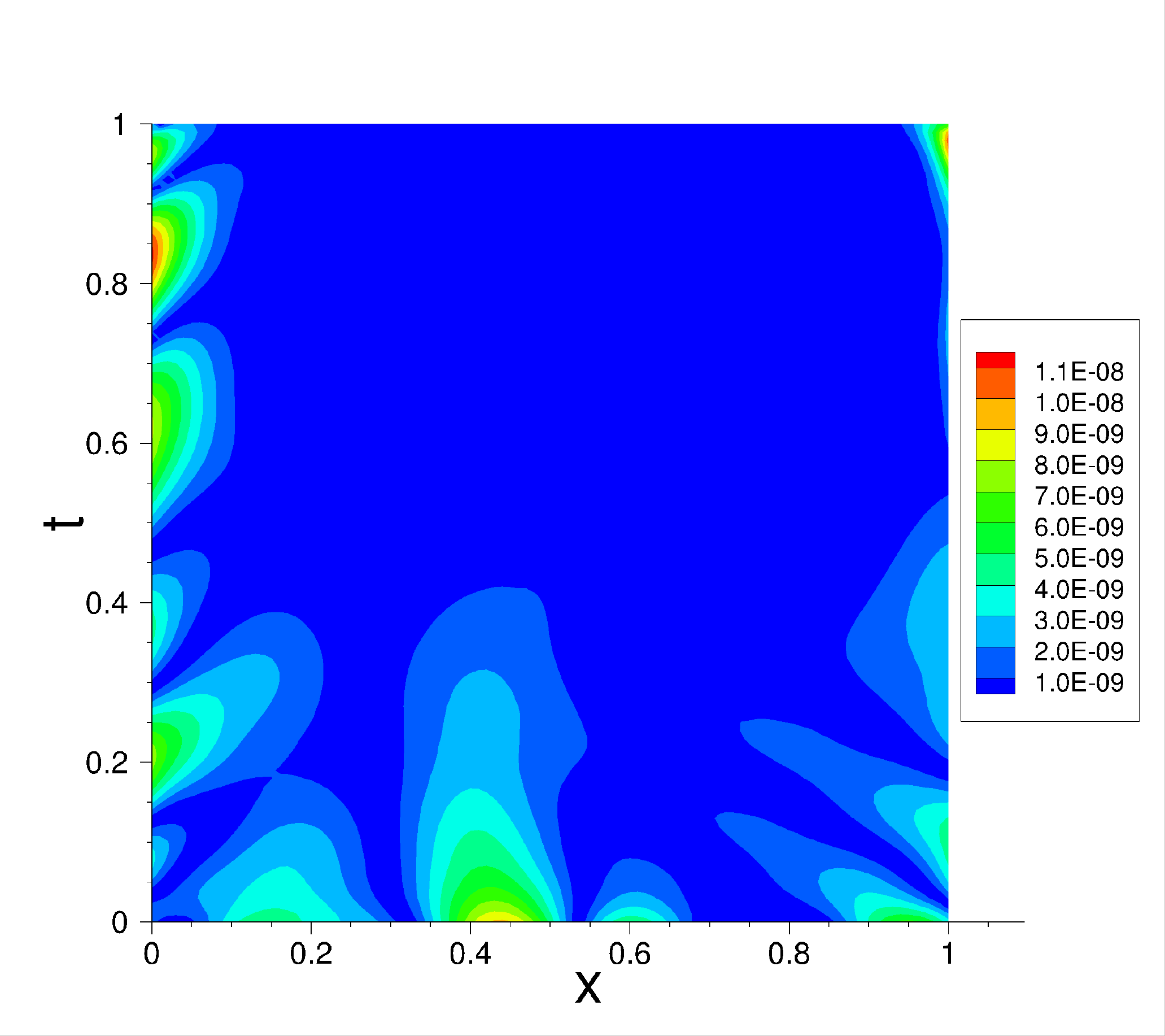}(b)
  }
  \caption{Burgers' equation: Distributions of (a) the exact solution
    and (b) the absolute error of the VarPro solution in the spatial-temporal plane.
    In (b), neural network [2, 150, 1], and $Q=31\times 31$ training collocation points.
  }
  \label{fg_15}
\end{figure}

In the second nonlinear example we use the viscous Burgers' equation to
test the VarPro method. Consider the spatial-temporal domain,
$Q=\{(x,t)\ |\ x\in[0,1],\ t\in[0,1] \}$, and
the following initial/boundary value problem on $\Omega$,
\begin{subequations}\label{eq_31}
  \begin{align}
    & \frac{\partial u}{\partial t} + u\frac{\partial u}{\partial x}
    - \nu\frac{\partial^2u}{\partial x^2} = f(x,t), \\
    &
    u(0,t) = g_1(t), \quad
    u(1,t) = g_2(t), \\
    &
    u(x,0) = h(x).
  \end{align}
\end{subequations}
In the above equations, $u(x,t)$ is the field function to be solved for,
$\nu=0.05$, $f(x,t)$ is a prescribed source term, $g_1$ and $g_2$
are the boundary conditions, and $h$ is the initial condition.
We choose $f$, $g_1$, $g_2$ and $h$ such that this problem has the the following
analytic solution,
\begin{equation}\label{eq_32}
  u(x,t) = \left[2\cos\left(\pi x+\frac{2\pi}{5} \right)
    + \frac32\cos\left(2\pi x-\frac{3\pi}{5} \right)
    \right]
  \left[2\cos\left(\pi t+\frac{2\pi}{5} \right)
    + \frac32\cos\left(2\pi t-\frac{3\pi}{5} \right)
    \right].
\end{equation}
Figure~\ref{fg_15}(a) shows the distribution of this analytic solution
in the spatial-temporal plane.

\begin{table}[tb]
  \centering
  \begin{tabular}{ll | ll}
    \hline
    parameter & value & parameter & value \\ \hline
    domain & $(x,t)\in[0,1]\times[0,1]$ & block time marching & none \\
    neural network & $[2, M, 1]$ & training points $Q$ & $Q_1\times Q_1$ \\
    $M$ &  varied & $Q_1$ & varied \\
    activation function & Gaussian & testing points & $Q_2\times Q_2$ \\
    random seed & $10$ & $Q_2$ & 101 \\
    initial guess $\bm\theta_0$ & random values on $[-R_m,R_m]$ & $R_m$ & $1.0$ \\
    $\delta$ (Algorithm~\ref{alg_3}) &  un-used
    & $p$ (Algorithm~\ref{alg_3}) & un-used \\
    max-subiterations & $0$ (no subiteration) & threshold (Algorithm~\ref{alg_3}) & $1E-12$ \\
    max-iterations-newton & $50$ & tolerance-newton & $1E-8$ \\
    \hline
  \end{tabular}
  \caption{Burgers' equation:
    main simulation parameters of the VarPro method.}
  \label{tab_7}
\end{table}

We employ neural networks with an architecture $[2, M, 1]$ in
the VarPro simulations, where $M$ is varied systematically in the tests.
The two input nodes represent $(x,t)$ and the linear output node
represents the solution field $u(x,t)$.
The Gaussian activation function, $\sigma(x)=e^{-x^2}$, is employed
in all the hidden nodes.
A uniform set of $Q=Q_1\times Q_1$ collocation points in
the spatial-temporal domain is used to train the neural network
with the VarPro method, and $Q_1$ is varied systematically in the tests.
The trained neural network is evaluated on a larger set of $Q_2\times Q_2$
uniform grid points to attain the field solution, which
is then compared with the analytic solution~\eqref{eq_32} to compute
the errors. The main simulation parameters for this problem
are summarized in Table~\ref{tab_7}.

Figure~\ref{fg_15}(b) illustrates the distribution of the absolute error
of a VarPro solution in the spatial-temporal plane.
This is computed using a neural network $[2, 150, 1]$ with
a uniform set of $Q=31\times 31$ training collocation points.
The VarPro solution can be observed to be quite accurate,
with a maximum error on the order $10^{-8}$ in the domain.

\begin{figure}
  \centerline{
    \includegraphics[width=2.5in]{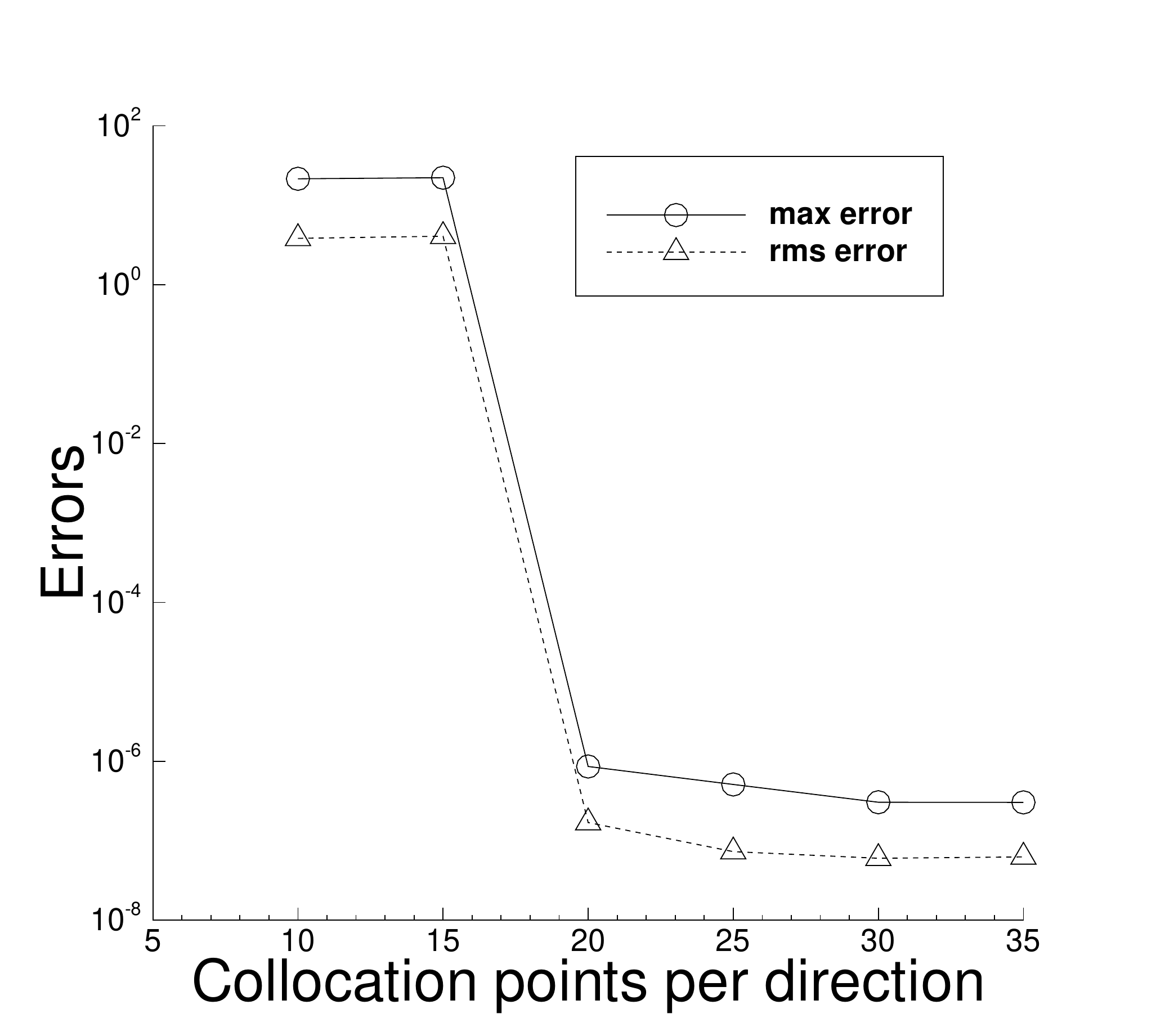}(a)
    \includegraphics[width=2.5in]{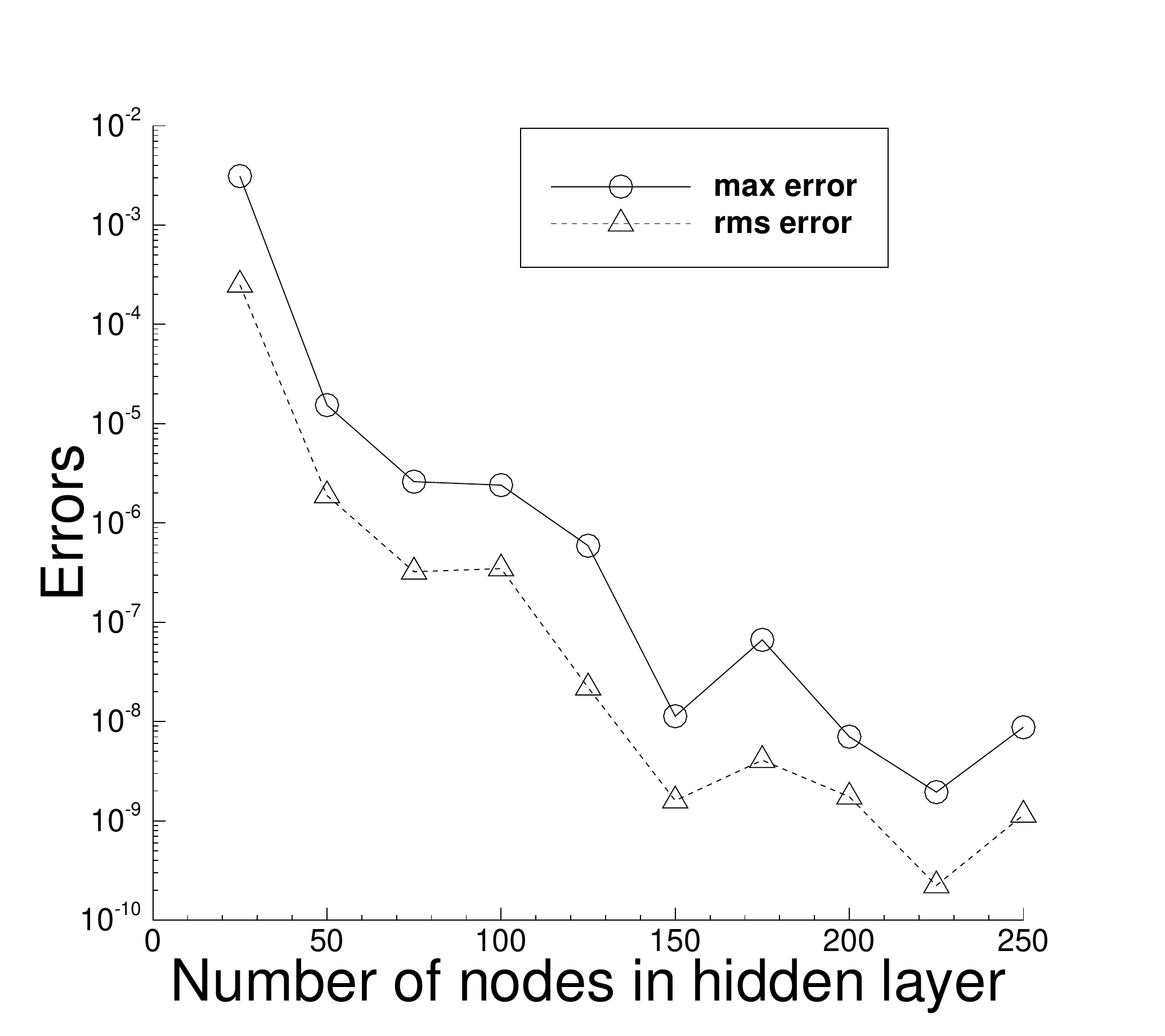}(b)
  }
  \caption{Burgers' equation: The maximum/rms errors of
    the VarPro solution versus (a) the number of collocation points
    in each direction ($Q_1$), and (b) the number of nodes in the hidden layer ($M$)
    of the neural network.
    Neural network $[2, M, 1]$, $Q=Q_1\times Q_1$ uniform collocation points.
    $M=100$ in (a) and is varied in (b). $Q_1=31$ in (b) and is varied in (a).
  }
  \label{fg_16}
\end{figure}

Figure~\ref{fg_16} illustrates the convergence behavior of the VarPro method for
solving the Burgers' equation.
In these tests the neural network is given by $[2, M, 1]$, where $M$ is either fixed at $M=100$
or varied between $M=25$ and $M=250$.
A set of $Q=Q_1\times Q_1$ uniform training collocation points
is used, where $Q_1$ is either fixed at $Q_1=31$ or varied 
between $Q_1=10$ and $Q_1=35$.
Figure~\ref{fg_16}(a) shows the maximum/rms errors of
the VarPro solution as a function of $Q_1$, corresponding to a fixed $M=100$ for
the neural network. The error behavior is not quite regular.
With a smaller $Q_1$ (e.g.~$10$ or $15$) the errors are at a level $1\sim 10$,
while with a larger $Q_1$ ($20$ and beyond) the errors abruptly drop to
a level around $10^{-8}\sim 10^{-6}$. We observe that with the smaller $Q_1=10$ and $15$
the Newton iteration fails to converge to the prescribed tolerance within
the prescribed maximum number of iterations.
Figure~\ref{fg_16}(b) shows the maximum/rms errors as a function of $M$,
corresponding to a fixed $Q=31\times 31$ for the collocation points.
The errors can be observed to decrease approximately exponentially
with increasing $M$.


\begin{table}[tb]
  \centering
  \begin{tabular}{l | l| ll| ll}
    \hline
    $[-R_m,R_m]$ & $M$ & VarPro & & ELM & \\ \cline{3-6}
     &  & max-error & rms-error & max-error & rms-error \\ \hline
    $R_m=1$ & $25$ & $3.111E-3$ & $2.499E-4$ & $6.382E+0$ & $8.382E-1$ \\
    & $50$ & $1.538E-5$ & $1.890E-6$ & $6.669E-2$ & $1.016E-2$ \\
    & $75$ & $2.603E-6$ & $3.222E-7$ & $1.216E-2$ & $1.406E-3$ \\
    & $100$ & $2.406E-6$ & $3.471E-7$ & $4.189E-4$ & $6.540E-5$  \\
    & $125$ & $5.894E-7$ & $2.193E-8$ & $1.088E-4$ & $1.099E-5$ \\
    & $150$ & $1.131E-8$ & $1.599E-9$ & $4.387E-6$ & $5.623E-7$  \\ \hline
    $R_m=R_{m0}=0.9$ & $25$ & $1.592E-3$ & $2.150E-4$ & $3.501E+0$ & $6.245E-1$  \\
    & $50$ & $1.069E-5$ & $1.794E-6$ & $1.390E-1$ & $7.607E-3$  \\
    & $75$ & $3.049E-6$ & $3.576E-7$ & $1.948E-2$ & $1.524E-3$ \\
    & $100$ & $4.831E-6$ & $3.311E-7$ & $7.639E-4$ & $4.780E-5$  \\
    & $125$ & $4.163E-7$ & $7.029E-8$ & $5.864E-5$ & $6.823E-6$  \\
    & $150$ & $6.063E-8$ & $7.242E-9$ & $2.000E-6$ & $2.877E-7$  \\
    \hline
  \end{tabular}
  \caption{Burgers' equation:
    comparison of the maximum/rms errors of the VarPro and ELM solutions.
    Neural network $[2, M, 1]$ with $M$ varied;
    fixed $Q=31\times 31$ training collocation points.
  }
  \label{tab_8}
\end{table}

Table~\ref{tab_8} provides a comparison of the solution errors
obtained using the VarPro method and the ELM method~\cite{DongY2021,DongL2021}
for the Burgers' equation.
These are computed using a fixed uniform set of
$Q=31\times 31$ collocation points and
a series of neural networks with the architecture $[2, M, 1]$, where
$M$ is varied between $M=25$ and $M=150$.
The random hidden-layer coefficients in ELM are set,
and the hidden-layer coefficients in VarPro are initialized,
by using $R_m=1$ and $R_m=R_{m0}=0.9$ in the tests.
It is evident that the VarPro method produces
significantly more accurate results than the ELM method.

\subsubsection{Nonlinear Klein-Gordon Equation}
\label{sec:sg}

\begin{figure}
  \centerline{
    \includegraphics[width=2.0in]{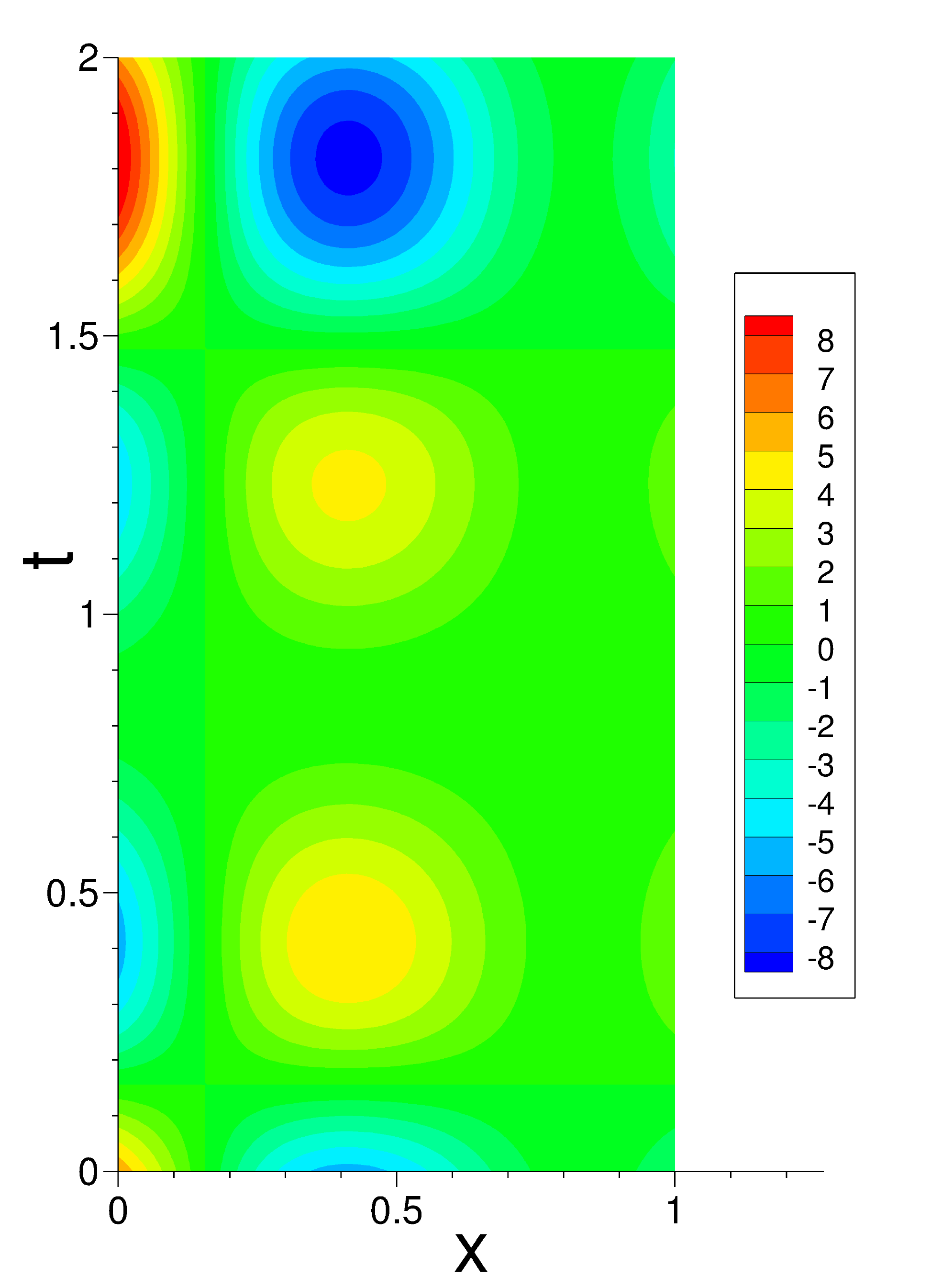}(a)
    \includegraphics[width=2.0in]{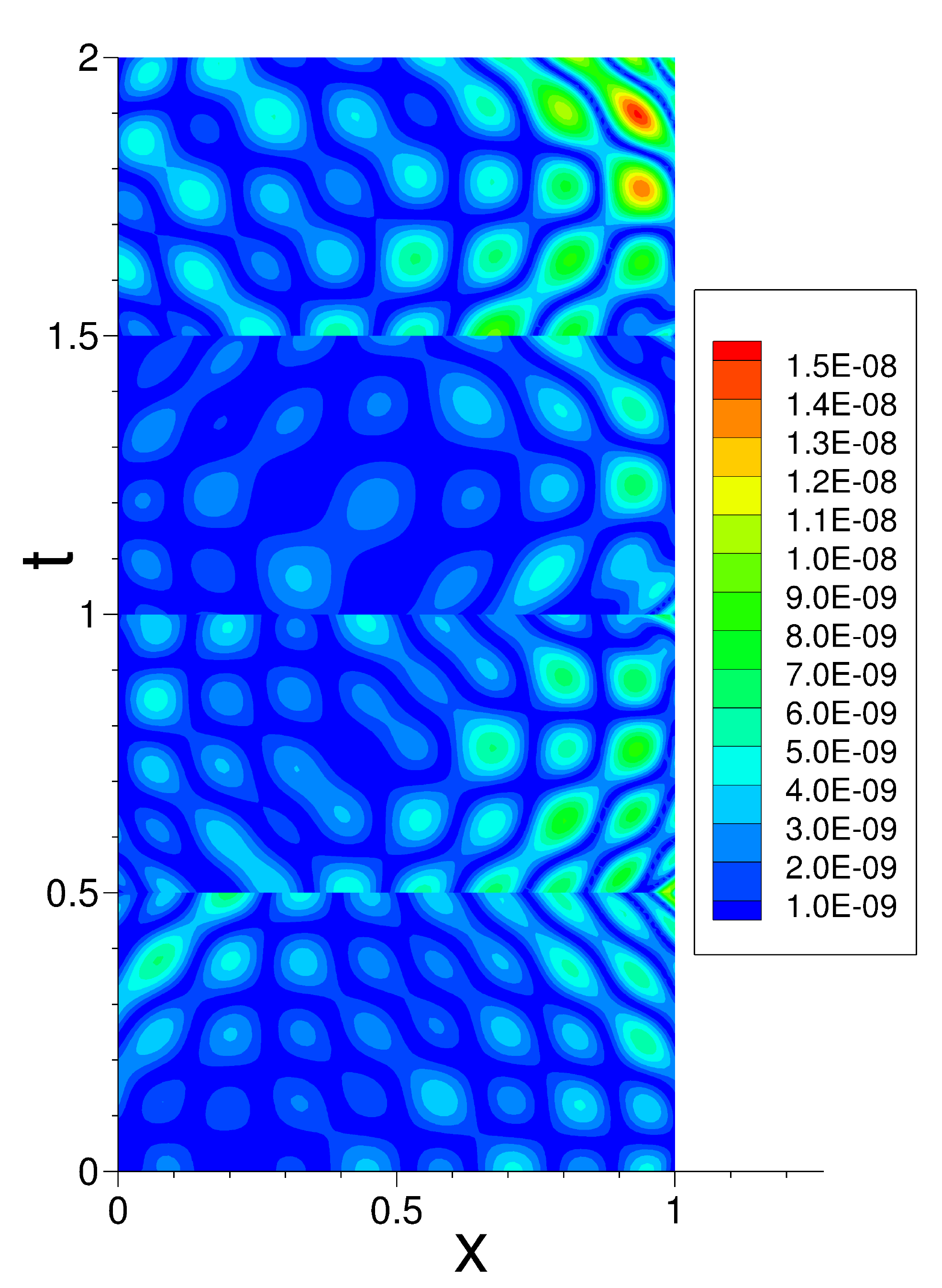}(b)
  }
  \caption{Nonlinear Klein-Gordon equation: Distributions of (a) the exact solution
    and (b) the absolute error of the VarPro solution in the spatial-temporal plane.
    In (b), $4$ uniform time blocks in domain, neural network $[2, 200, 1]$,
    $Q=21\times 21$ uniform collocation points per time block.
  }
  \label{fg_17}
\end{figure}

In the last example we use the nonlinear Klein-Gordon equation~\cite{Strauss1978}
to test the VarPro method.
Consider the spatial-temporal domain,
$\Omega=\{(x,t)\ |\ x\in[0,1],\ t\in[0,2] \}$, and the initial/boundary value
problem with the nonlinear Klein-Gordon equation on $\Omega$,
\begin{subequations}
  \begin{align}
    &
    \frac{\partial^2u}{\partial t^2} - \frac{\partial^2u}{\partial x^2} + u
    + \sin(u) = f(x,t), \\
    & u(0,t) = g_1(t), \quad
    u(1,t) = g_2(t), \\
    & u(x,0) = h_1(x), \quad
    \left.\frac{\partial u}{\partial t}\right|_{(x,0)} = h_2(x),
  \end{align}
\end{subequations}
where $u(x,t)$ is the field function to be solved for, $f(x,t)$ is a
prescribed source term, $g_1$ and $g_2$ are the boundary conditions,
and $h_1$ and $h_2$ are the initial conditions.
We choose the source term, the boundary and initial conditions appropriately such
that this problem has the following analytic solution,
\begin{equation}\label{eq_34}
  u(x,t) = \left[2\cos\left(\pi x + \frac{\pi}{5} \right)
    + \frac95\cos\left(2\pi x + \frac{7\pi}{20}  \right)  \right]
  \left[2\cos\left(\pi t + \frac{\pi}{5} \right)
    + \frac95\cos\left(2\pi t + \frac{7\pi}{20}  \right)  \right].
\end{equation}
We employ this analytic solution to test the accuracy of the VarPro method.
Figure~\ref{fg_17}(a) shows the distribution of this analytic
solution in the spatial-temporal plane.

\begin{table}[tb]
  \centering
  \begin{tabular}{ll | ll}
    \hline
    parameter & value & parameter & value \\ \hline
    domain & $(x,t)\in[0,1]\times[0,2]$ & time blocks & $4$ \\
    neural network & $[2, M, 1]$ & training points $Q$ & $Q_1\times Q_1$ \\
    $M$ &  varied & $Q_1$ & varied \\
    activation function & Gaussian & testing points & $Q_2\times Q_2$ \\
    random seed & $22$ & $Q_2$ & 101 \\
    initial guess $\bm\theta_0$ & random values on $[-R_m,R_m]$ & $R_m$ & $1.0$ \\
    $\delta$ (Algorithm~\ref{alg_3}) &  un-used
    & $p$ (Algorithm~\ref{alg_3}) & un-used \\
    max-subiterations & $0$ (no subiteration) & threshold (Algorithm~\ref{alg_3}) & $1E-12$ \\
    max-iterations-newton & $20$ & tolerance-newton & $1E-8$ \\
    \hline
  \end{tabular}
  \caption{Nonlinear Klein-Gordon equation:
    main simulation parameters of the VarPro method.}
  \label{tab_9}
\end{table}

We employ the block time marching scheme~\cite{DongL2021} together with
the VarPro method to solve this problem.
We use $4$ uniform time blocks in the domain, and on each time block
employ a neural network with the architecture $[2,M,1]$,
where $M$ is varied in the tests.
The two input nodes represent $(x,t)$,
and the linear output node represents the field solution $u(x,t)$.
The Gaussian activation function $\sigma(x)=e^{-x^2}$ is employed
for all the hidden nodes. 
On each time block a uniform set of $Q=Q_1\times Q_1$ collocation points,
where $Q_1$ is varied,
is used to train the neural network.
The trained neural network is evaluated on a larger set of $Q_2\times Q_2$
uniform grid points to obtain the field solution $u(x,t)$,
which is compared with the exact solution~\eqref{eq_34}
to compute the errors of the VarPro simulation.
Table~\ref{tab_9} summarizes the main simulation parameters for
this problem.

Figure~\ref{fg_17}(b) shows the distribution of the absolute error
of a VarPro simulation in the spatial-temporal plane.
In this simulation the neural network architecture
is given by $[2, 200, 1]$, and a uniform set of $Q=21\times 21$ training
collocation points is used on each time block.
The maximum error level is around $10^{-8}$ on the overall domain, indicating that
the VarPro result is quite accurate.

\begin{figure}
  \centerline{
    \includegraphics[width=2.5in]{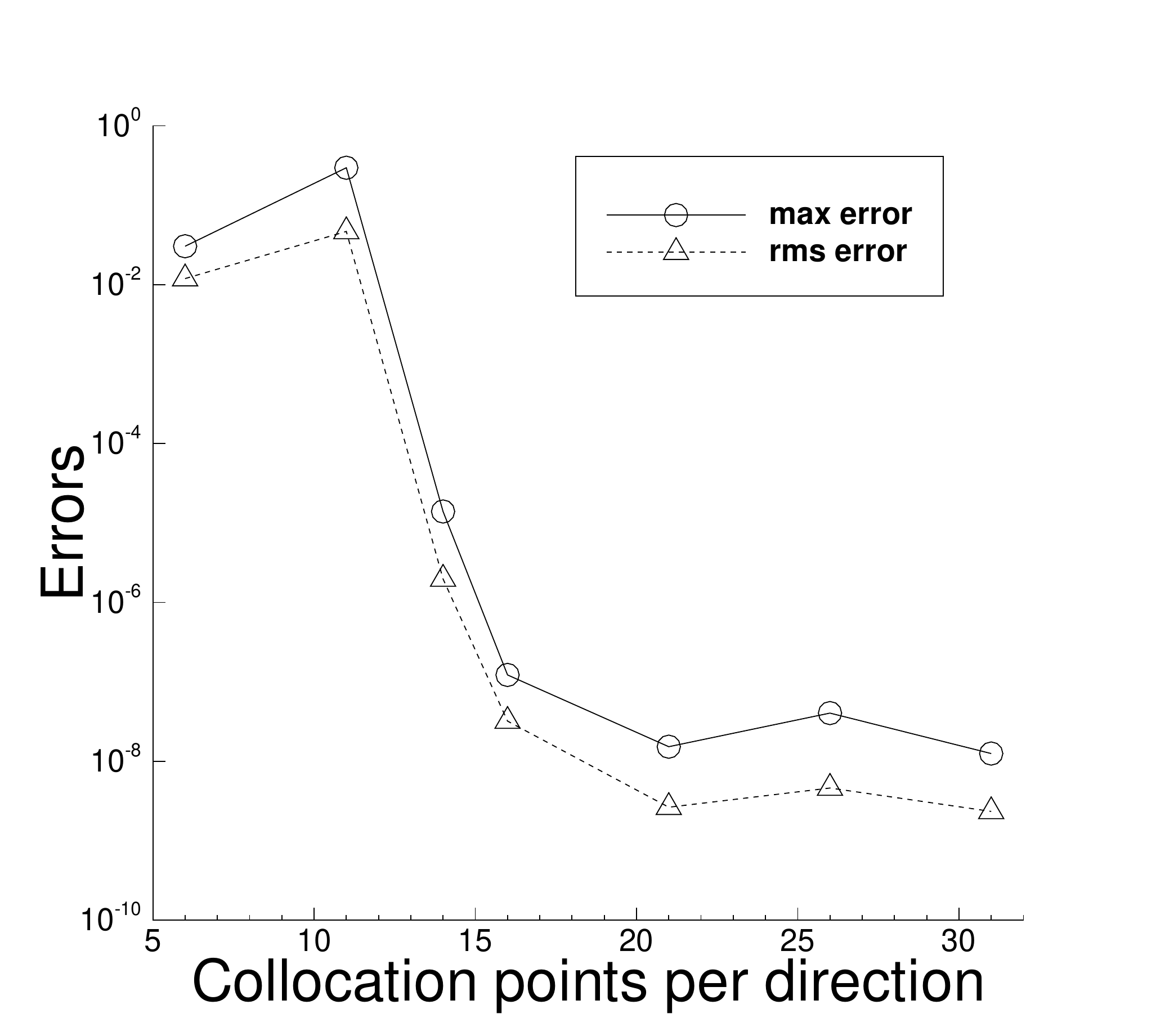}(a)
    \includegraphics[width=2.5in]{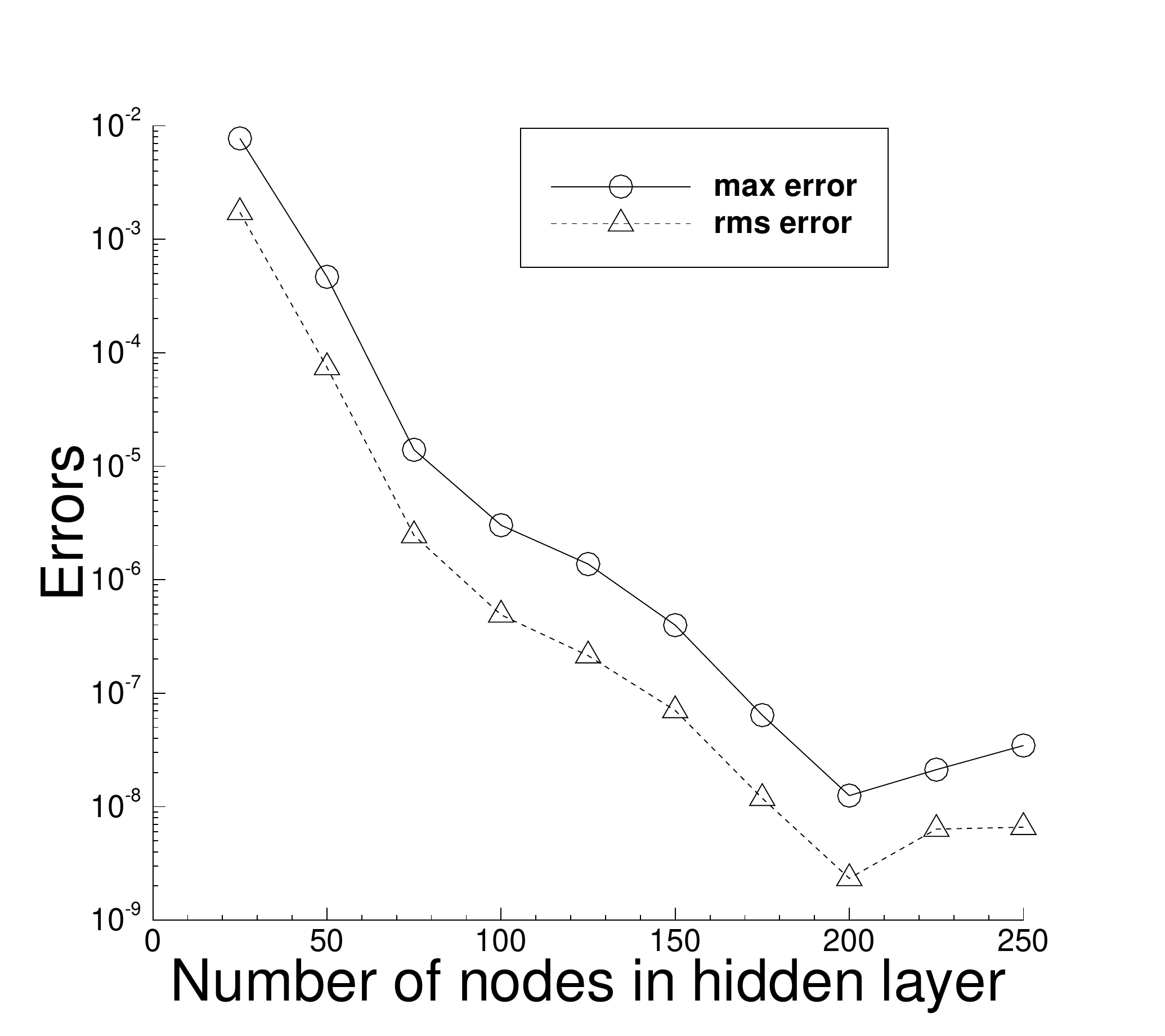}(b)
  }
  \caption{Nonlinear Klein-Gordon equation: the maximum/rms errors of the VarPro solution
    in the overall domain versus (a) the number of collocation points
    per direction ($Q_1$) in each time block, and (b) the number of nodes in the hidden layer ($M$)
    of the neural network.
    In (a,b), neural network $[2, M, 1]$, with $Q=Q_1\times Q_1$ training collocation points.
    $M=200$ in (a) and is varied in (b).
    $Q_1=31$ in (b) and is varied in (a).
  }
  \label{fg_18}
\end{figure}

Figure~\ref{fg_18} illustrates the convergence behavior of the VarPro method
for solving the nonlinear Klein-Gordon equation.
Here we employ the neural network $[2, M, 1]$, and $Q=Q_1\times Q_1$ uniform
collocation points in each time block.
In the first group of tests we fix $M=200$ and vary $Q_1$ systematically.
In the second group of tests we fix $Q_1=31$ and vary $M$ systematically.
The maximum and rms errors of the VarPro solution in the overall domain are computed
for each case.
Figure \ref{fg_18}(a) shows these errors as a function of $Q_1$ for the first group of tests,
and Figure \ref{fg_18}(b) shows these errors as a function of $M$ for
the second group of tests.
These results indicate that the VarPro errors decrease approximately exponentially
with increasing number of collocation points or with increasing number of nodes in
the hidden layer.
We also notice some irregularity in the errors of Figure~\ref{fg_18}(a)
when the number of collocation points is small.

\section{Concluding Remarks}
\label{sec:summary}



In this paper we have presented a variable projection-based method
together with artificial neural networks for numerically approximating
linear and nonlinear partial differential equations.
The basic idea of variable projection (VarPro) for solving separable nonlinear
least squares problems is to distinguish
the linear parameters from the nonlinear parameters,
and then eliminate the linear parameters 
to attain a reduced formulation of the problem. One can then
solve the reduced problem for the nonlinear parameters first, and then
compute the linear parameters by using the linear least squares method afterwards.

Approximating linear PDEs (with linear boundary/initial conditions)
by variable projection and artificial neural networks
is conceptually straightforward. In this case,
in the resultant nonlinear least squares problem
the coefficients in the linear output layer
are the linear parameters, and those in the hidden layers are the nonlinear
parameters. The output-layer coefficients
can be expressed in terms of the hidden-layer coefficients by solving
a linear least squares problem, and
they can be eliminated from the problem.
The reduced  problem
involves only the hidden-layer coefficients, and
it can be solved by the nonlinear least squares method.
The main issues with the VarPro implementation lie in
the computations of the residual function and the Jacobian matrix of the reduced
problem. We have discussed in some detail how to implement these computations
with neural networks
in Algorithms~\ref{alg_1} and~\ref{alg_2} and in the Remarks~\ref{rem_1}
and~\ref{rem_2}.

For approximating nonlinear PDEs,
or linear PDEs with nonlinear boundary/initial conditions,
the variable projection approach with the artificial
neural networks  cannot be directly used.
This is because the resultant nonlinear least squares problem is not separable.
In this case, all the weight/bias coefficients in the neural network become nonlinear
parameters, even if the
output layer contains no activation function.

To overcome this issue, we have presented a Newton/variable projection (Newton-VarPro)
method for approximating
nonlinear PDEs, or linear PDEs with nonlinear boundary/initial conditions.
We first linearize the problem, with a particular linearized form
for the Newton iteration. The linearization is formulated  in terms of
the updated approximation field, not the increment field.
This linearization is critical to the accuracy of the current Newton-VarPro method.
The linearized system
can be solved by the variable projection approach together with
artificial neural networks.
Therefore for solving nonlinear PDEs, our method involves an overall
Newton iteration, and within each iteration the variable projection method
is used to solve the linearized problem to attain the updated approximation field.
Upon convergence of the Newton iteration, the solution to
the nonlinear problem is represented by the weight/bias coefficients
of the neural network.


We have presented ample numerical examples with linear
and nonlinear PDEs to test the accuracy of the variable projection method.
It is observed that, for smooth field solutions,
the errors of the VarPro method decrease exponentially
or nearly exponentially with increasing number of collocation points
or with increasing number of output-layer coefficients.
The test results unequivocally show that the VarPro method is highly accurate.
Even with a fairly small number of nodes in the neural network,
or with a fairly small set of collocation points,
the VarPro method can produce very accurate simulation results.

In particular, we have compared extensively
the current VarPro method with the extreme learning
machine (ELM) method~\cite{DongY2021},
which is arguably the most accurate neural network-based PDE solver so
far~\cite{DongL2021,DongY2021}.
Under the same simulation conditions and settings,
the VarPro method generally leads to significantly more accurate results
than the ELM method, especially in cases with a fairly small or a moderate number
of nodes in the neural network.

\begin{table}[tb]
  \centering
  \begin{tabular}{l| l|  l | l | l}
    \hline
    Test problem & neural  & collocation
    & VarPro training & ELM training \\
     & network & points & time (seconds) & time (seconds)    \\ \hline
    Advection equation & $[2, 100, 1]$  & $10\times 10$ & $7.4$ & $0.29$ \\ \cline{3-5}
    &  &  $15\times 15$ & $15.4$ & $0.31$ \\ \cline{3-5}
    & & $20\times 20$ & $78.2$ & $0.32$ \\
    \hline
    Nonlinear Helmholtz & $[2, 200, 1]$  & $10\times 10$ & $1.7$ & $2.0$ \\ \cline{3-5}
    equation &  &  $15\times 15$ & $36.5$ & $2.8$ \\ \cline{3-5}
    & & $20\times 20$ & $78.1$ & $4.0$ \\
    \hline
  \end{tabular}
  \caption{Comparison of the computational cost (network training time)
    between VarPro and ELM for solving the advection equation (Section~\ref{sec:advec})
    and the nonlinear Helmholtz equation (Section~\ref{sec:nonl_helm}).
    The hidden-layer coefficients in ELM are set/fixed to, and the hidden-layer
    coefficients in VarPro are initialized to, uniform random values
    from $[-R_m,R_m]$ with $R_m=1.0$.
    The cases in this table are selected from and correspond to those cases
    in Tables~\ref{tab_4} and \ref{tab_6}.
    Please refer to the corresponding cases in Tables~\ref{tab_4} and~\ref{tab_6}
    for the VarPro/ELM errors.
  }
  \label{tab_10}
\end{table}


While the VarPro method is significantly superior to ELM in terms of
the accuracy, its computational cost (i.e.~training time of the neural network)
is generally much higher than that of the ELM method.
This is because in VarPro one needs to solve the reduced problem for
the hidden-layer coefficients
by a nonlinear least squares computation, apart from the computation for
the linear output-layer coefficients.
In contrast, in ELM only the linear output-layer coefficients are computed,
while the hidden-layer coefficients in the neural network are fixed to
the random values that are pre-set.
Table~\ref{tab_10} illustrates this points with a list of the network training time for the
VarPro method and the ELM method with selected cases in solving
the the advection equation (Section~\ref{sec:advec})
and the nonlinear Helmholtz equation (Section~\ref{sec:nonl_helm}).


The variable projection method is a powerful technique
for training artificial neural networks, providing a considerably
superior accuracy for scientific machine learning, as demonstrated
by ample numerical examples in the current paper.
The Newton-VarPro method developed herein
provides an effective tool and enables the use of the variable projection
strategy to tackle nonlinear problems in scientific machine learning.
The application potential of this technique is enormous.
This and related aspects, as well as further studies and improvements,
of this technique will be pursued in a future endeavor.

\section*{Acknowledgement}
This work was partially supported by
NSF (DMS-2012415).

\bibliographystyle{plain}
\bibliography{varpro,elm,elm1,mypub,dnn,dnn1,sem,obc}

\begin{thebibliography}{10}

\bibitem{AravkinL2012}
A.Y. Aravkin and T.~van Leeuwen.
\newblock Estimating nuisance parameters in inverse problems.
\newblock {\em Inverse Problems}, 28:115016, 2012.

\bibitem{AskhamK2018}
T.~Askham and J.N. Kutz.
\newblock Variable projection methods for an optimized dynamic mode
  decomposition.
\newblock {\em SIAM J. Applied Dynamical Systems}, 17:380--416, 2018.

\bibitem{Bjorck2015}
A.~Bjorck.
\newblock {\em Numerical Methods in Matrix Computations}.
\newblock Springer, 2015.

\bibitem{BranchCL1999}
M.A. Branch, T.F. Coleman, and Y.~Li.
\newblock A subspace, interior, and conjugate gradient method for large-scale
  bound-constrained minimization problems.
\newblock {\em SIAM Journal on Scientific Computing}, 21:1--23, 1999.

\bibitem{ByrdSS1988}
R.H. Byrd, R.B. Schnabel, and G.A. Shultz.
\newblock Approximate solution of the trust region problem by minimization over
  two-dimensional subspaces.
\newblock {\em Math. Programming}, 1988.

\bibitem{ChenGCL2019}
G.-Y. Chen, M.~Gan, C.L.P. Chen, and H.-X. Li.
\newblock A regularized variable projection algorithm for separable nonlinear
  least-squares problems.
\newblock {\em IEEE Transactions on Automatic Control}, 64:526--537, 2019.

\bibitem{ChungHN2006}
J.~Chung, E.~Haber, and J.~Nagy.
\newblock Numerical methods for coupled super-resolution.
\newblock {\em Inverse Problems}, 22:1261--1272, 2006.

\bibitem{CornelioPN2014}
A.~Cornelio, E.L. Piccolomini, and J.G. Nagy.
\newblock Constraned numerical optimization methods for blind deconvolution.
\newblock {\em Numer. Algor.}, 65:23--42, 2014.

\bibitem{CyrGPPT2020}
E.C. Cyr, M.A. Gulian, R.G. Patel, M.~Perego, and N.A. Trask.
\newblock Robust training and initialization of deep neural networks: an
  adaptive basis viewpoint.
\newblock {\em Proceedings of Machine Learning Research}, 107:1--26, 2020.

\bibitem{DennisS1996}
J.E. Dennis and R.B. Schnabel.
\newblock {\em Numerical Methods for Unconstrained Optimization and Nonlinear
  Equations}.
\newblock SIAM, 1996.

\bibitem{Dong2015clesobc}
S.~Dong.
\newblock A convective-like energy-stable open boundary condition for
  simulations of incompressible flows.
\newblock {\em Journal of Computational Physics}, 302:300--328, 2015.

\bibitem{Dong2018}
S.~Dong.
\newblock Multiphase flows of {N} immiscible incompressible fluids: a
  reduction-consistent and thermodynamically-consistent formulation and
  associated algorithm.
\newblock {\em Journal of Computational Physics}, 361:1--49, 2018.

\bibitem{DongL2021}
S.~Dong and Z.~Li.
\newblock Local extreme learning machines and domain decomposition for solving
  linear and nonlinear partial differential equations.
\newblock {\em Computer Methods in Applied Mechanics and Engineering},
  387:114129, 2021.
\newblock (also arXiv:2012.02895).

\bibitem{DongL2021bip}
S.~Dong and Z.~Li.
\newblock A modified batch intrinsic plascity method for pre-training the
  random coefficients of extreme learning machines.
\newblock {\em Journal of Computational Physics}, 445:110585, 2021.
\newblock (also arXiv:2103.08042).

\bibitem{DongN2020}
S.~Dong and N.~Ni.
\newblock A method for representing periodic functions and enforcing exactly
  periodic boundary conditions with deep neural networks.
\newblock {\em Journal of Computational Physics}, 435:110242, 2021.

\bibitem{DongS2012}
S.~Dong and J.~Shen.
\newblock A time-stepping scheme involving constant coefficient matrices for
  phase field simulations of two-phase incompressible flows with large density
  ratios.
\newblock {\em Journal of Computational Physics}, 231:5788--5804, 2012.

\bibitem{DongY2021}
S.~Dong and J.~Yang.
\newblock On computing the hyperparameter of extreme learning machines:
  algorithm and application to computational {PDE}s and comparison with
  classical and high-order finite elements.
\newblock {\em arXiv:2110.14121}, 2021.

\bibitem{DwivediS2020}
V.~Dwivedi and B.~Srinivasan.
\newblock Physics informed extreme learning machine (pielm) $-$ a rapid method
  for the numerical solution of partial differential equations.
\newblock {\em Neurocomputing}, 391:96--118, 2020.

\bibitem{EY2018}
W.~E and B.~Yu.
\newblock The deep {R}itz method: a deep learning-based numerical algorithm for
  solving variational problems.
\newblock {\em Communications in Mathematics and Statistics}, 6:1--12, 2018.

\bibitem{ErichsonZMBKA2020}
N.B. Erichson, P.~Zheng, K.~Manohar, J.N.~Kutz S.L.~Brunton, and A.Y. Aravkin.
\newblock Sparse principal component analysis via variable projection.
\newblock {\em SIAM J. Appl. Math.}, 80:977--1002, 2020.

\bibitem{FabianiCRS2021}
G.~Fabiani, F.~Calabro, L.~Russo, and C.~Siettos.
\newblock Numerical solution and bifurcation analysis of nonlinear partial
  differential equations with extreme learning machines.
\newblock {\em Journal of Scientific Computing}, 89:44, 2021.

\bibitem{GanCCC2018}
M.~Gan, C.L.P. Chen, H.-Y. Chen, and L.~Chen.
\newblock On some separated algorithms for separable nonlinear least squares
  problems.
\newblock {\em IEEE Transactions on Cybernetics}, 48:2866--2974, 2018.

\bibitem{GolubP1973}
G.H. Golub and V.~Pereyra.
\newblock The differentiation of pseudo-inverse and nonlinear least squares
  problems whose variables separate.
\newblock {\em SIAM J. Numer. Anal.}, 10:413--432, 1973.

\bibitem{GolubP2003}
G.H. Golub and V.~Pereyra.
\newblock Separable nonlinear least squares: the variable projection method and
  its applications.
\newblock {\em Inverse Problems}, 19:R1--R26, 2003.

\bibitem{GoodfellowBC2016}
I.~Goodfellow, Y.~Bengio, and A.~Courville.
\newblock {\em Deep Learning}.
\newblock The MIT Press, 2016.

\bibitem{HeX2019}
J.~He and J.~Xu.
\newblock {MgNet}: A unified framework for multigrid and convolutional neural
  network.
\newblock {\em Science China Mathematics}, 62:1331--1354, 2019.

\bibitem{HendrycksG2020}
D.~Hendrycks and K.~Gimpel.
\newblock Gaussian error linear units {(GELU)}.
\newblock {\em arXiv:1606.08415}, 2016.

\bibitem{HerringNR2018}
J.L. Herring, J.G. Nagy, and L.~Ruthotto.
\newblock {LAP}: A linearize and project method for solving inverse problems
  with coupled variables.
\newblock {\em Sampling Theory in Signal and Image Processing}, 17:127--151,
  2018.

\bibitem{HuangZS2006}
G.-B. Huang, Q.-Y. Zhu, and C.-K. Siew.
\newblock Extreme learning machine: theory and applications.
\newblock {\em Neurocomputing}, 70:489--501, 2006.

\bibitem{HuangCS2006}
G.B. Huang, L.~Chen, and C.-K. Siew.
\newblock Universal approximation using incremental constructive feedforward
  networks with random hidden nodes.
\newblock {\em IEEE Transactions on Neural Networks}, 17:879--892, 2006.

\bibitem{Karniadakisetal2021}
G.E. Karniadakis, G.~Kevrekidis, L.~Lu, P.~Perdikaris, S.~Wang, and L.~Yang.
\newblock Physics-informed machine learning.
\newblock {\em Nature Reviews Physics}, 3:422--440, 2021.

\bibitem{KarniadakisS2005}
G.E. Karniadakis and S.J. Sherwin.
\newblock {\em Spectral/hp element methods for computational fluid dynamics,
  2nd edn.}
\newblock Oxford University Press, 2005.

\bibitem{Kaufman1975}
L.~Kaufman.
\newblock A variable projection method for solving separable nonlinear least
  squares problems.
\newblock {\em BIT}, 15:49--57, 1975.

\bibitem{KaufmanP1978}
L.~Kaufman and V.~Pereyra.
\newblock A method for separable nonlinear least squares problems with
  separable equality constraints.
\newblock {\em SIAM J. Numer. Anal.}, 15:12--20, 1978.

\bibitem{KimL2008}
C.-T. Kim and J.-J. Lee.
\newblock Training two-layered feedforward networks with variable projection
  method.
\newblock {\em IEEE Transactions on Neural Networks}, 19:371--375, 2008.

\bibitem{Krogh1974}
F.T. Krogh.
\newblock Efficient implementation of a variable projection algorithm for
  nonlinear least squares problems.
\newblock {\em Commun. ACM}, 17:167--169, 1974.

\bibitem{LiangJHY2021}
S.~Liang, S.W. Jiang, J.~Harlim, and H.~Yang.
\newblock Solving {PDE}s on unknown manifolds with machine learning.
\newblock {\em arXiv:2106.06682}, 2021.

\bibitem{LinYD2019}
L.~Lin, Z.~Yang, and S.~Dong.
\newblock Numerical approximation of incompressible {N}avier-{S}tokes equations
  based on an auxiliary energy variable.
\newblock {\em Journal of Computational Physics}, 388:1--22, 2019.

\bibitem{LuMMK2021}
L.~Lu, X.~Meng, Z.~Mao, and G.E. Karniadakis.
\newblock Deep{XDE}: A deep learning library for solving differential
  equations.
\newblock {\em SIAM Review}, 63:208--228, 2021.

\bibitem{LuoY2020}
T.~Luo and H.~Yang.
\newblock Two-layer neural networks for partial differential equations:
  optimization and generlization theory.
\newblock {\em arXiv:2006.15733}, 2020.

\bibitem{MaoJK2020}
Z.~Mao, A.D. Jagtap, and G.E. Karniadakis.
\newblock Physics-informed neural networks for high-speed flows.
\newblock {\em Computer Methods in Applied Mechanics and Engineering},
  360:112789, 2020.

\bibitem{McLooneBI1998}
S.~McLoone, M.D. Brown, and G.~Irwin.
\newblock A hybrid linear/nonlinear training algorithm for feedforward neural
  networks.
\newblock {\em IEEE Transactions on Neural Networks}, 9:669--684, 1998.

\bibitem{MullenS2009}
K.M. Mullen and I.H.M. van Stokkum.
\newblock The variable projection algorithm in time-resolved spectroscopy,
  microscopy and mass spectrometry applications.
\newblock {\em Numer. Algor.}, 51:319--340, 2009.

\bibitem{NewmanCCR2021}
E.~Newman, J.~Chung, M.~Chung, and L.~Ruthotto.
\newblock Slim{T}rain -- a stochastic approximation method for training
  separable deep neural networks.
\newblock {\em arXiv:2109.14002}, 2021.

\bibitem{NewmanRHW2020}
E.~Newman, L.~Ruthotto, J.~Hart, and B.~van Bloemen~Waanders.
\newblock Train like a ({V}ar){P}ro: Efficient training of neural networks with
  variable projection.
\newblock {\em arXiv:2007.13171}, 2020.

\bibitem{NocedalW1999}
J.~Nocedal and S.~Wright.
\newblock {\em Numerical Optimization}.
\newblock Springer, 1999.

\bibitem{OlearyR2013}
D.P. O'Leary and B.W. Rust.
\newblock Variable projection for nonlinear least squares problems.
\newblock {\em Comput. Optim. Appl.}, 54:579--593, 2013.

\bibitem{Osborne2007}
M.R. Osborne.
\newblock Separable least squares, variable projection, and the gauss-newton
  algorithm.
\newblock {\em Electronic Transactions on Numerical Analysis}, 28:1--15, 2007.

\bibitem{PereyraSW2006}
V.~Pereyra, G.~Scherer, and F.~Wong.
\newblock Variable projections neural network training.
\newblock {\em Mathematics and Computers in Simulation}, 73:231--243, 2006.

\bibitem{RaissiPK2019}
M.~Raissi, P.~Perdikaris, and G.E. Karniadakis.
\newblock Physics-informed neural networks: a deep learning framework for
  solving forward and inverse problems involving nonlinear partial differential
  equations.
\newblock {\em Journal of Computational Physics}, 378:686--707, 2019.

\bibitem{RuanoJF1991}
A.E.B. Ruano, D.J. Jones, and P.J. Fleming.
\newblock A new formulation of the learning problem of a neural network
  controller.
\newblock {\em Proc. 30th IEEE Conf. Decis. Control, Brighton, UK}, pages
  865--866, 1991.

\bibitem{RuheW1980}
A.~Ruhe and P.A. Wedin.
\newblock Algorithms for separable nonlinear least squares problems.
\newblock {\em SIAM Review}, 22:318--337, 1980.

\bibitem{Samaniegoetal2020}
E.~Samanaiego, C.~Anitescu, S.~Goswami, V.M. Nguyen-Thanh, H.~Guo, K.~Hamdia,
  X.~Zhuang, and T.~Rabczuk.
\newblock An energy approach to the solution of partial differential equations
  in computational mechanics via machine learning: concepts, implementation and
  applications.
\newblock {\em Computer Methods in Applied Mechanics and Engineering},
  362:112790, 2020.

\bibitem{ShearerG2013}
P.~Shearer and A.C. Gilbert.
\newblock A generalization of variable elimination for separable inverse
  problems beyond least squares.
\newblock {\em Inverse Problems}, 29:045003, 2013.

\bibitem{SimaH2007}
D.M. Sima and S.~Van Huffel.
\newblock Separable nonlinear least squares fitting with linear bound
  constraints and its application in magnetic resonance spectroscopy data
  quantification.
\newblock {\em Journal of Computational and Applied Mathematics}, 203:264--278,
  2007.

\bibitem{SirignanoS2018}
J.~Sirignano and K.~Spoliopoulos.
\newblock {DGM}: A deep learning algorithm for solving partial differential
  equations.
\newblock {\em Journal of Computational Physics}, 375:1339--1364, 2018.

\bibitem{SjobergV1997}
J.~Sjoberg and M.~Viberg.
\newblock Separable nonlinear least squares minimization - possible
  improvements for neural net fitting.
\newblock {\em Neural Networks for Signal Processing VII. Proceedings of IEEE
  Signal Processing Workshop}, 1997.

\bibitem{SongXHZ2020}
X.~Song, W.~Xu, K.~Hayami, and N.~Zheng.
\newblock Secant variable projection method for solving nonnegative separable
  least squares problems.
\newblock {\em Numerical Algorithms}, 85:737--761, 2020.

\bibitem{Strauss1978}
W.~Strauss.
\newblock Numerical solution of nonlinear klein-gordon equation.
\newblock {\em Journal of Computational Physics}, 28:271--278, 1978.

\bibitem{SzaboB1991}
B.~Szabo and I.~Babushka.
\newblock {\em Finite Element Analysis}.
\newblock John Wiley \& Sons, Inc., 1991.

\bibitem{LeeuwenA2021}
T.~van Leeuwen and A.Y. Aravkin.
\newblock Variable projection for nonsmooth problems.
\newblock {\em SIAM J. Sci. Comput.}, 43:S249--S268, 2021.

\bibitem{WangYP2020}
S.~Wang, X.~Yu, and P.~Perdikaris.
\newblock When and why {PINN}s fail to train: a neural tangent kernel
  perspective.
\newblock {\em arXiv:2007.14527}, 2020.

\bibitem{WeiglB1993}
K.~Weigl and M.~Berthod.
\newblock Neural networks as dynamical bases in function space.
\newblock {\em Report No 2124, INRIA, Sophis-Antipolis, France}, 1993.
\newblock URL: https://hal.inria.fr/inria-00074548/document.

\bibitem{WeiglB1994}
K.~Weigl and M.~Berthod.
\newblock Projection learning: alternative approach to the computation of the
  projection.
\newblock {\em Proc. European Symp. on Artificial Neural Networks, Brussels,
  Belgium}, pages 19--24, 1994.

\bibitem{WeiglGB1993}
K.~Weigl, G.~Giraudon, and M.~Berthod.
\newblock Application of projection learning to the detection of urban areas in
  {SPOT} satellite images.
\newblock {\em Report No 2143, INRIA, Sophia-Antipolis, France}, 1993.
\newblock URL: https://hal.inria.fr/inria-00074529.

\bibitem{YangD2019}
Z.~Yang and S.~Dong.
\newblock An unconditionally energy-stable scheme based on an implicit
  auxiliary energy variable for incompressible two-phase flows with different
  densities involving only precomputable coefficient matrices.
\newblock {\em Journal of Computational Physics}, 393:229--257, 2019.

\bibitem{YangD2020}
Z.~Yang and S.~Dong.
\newblock A roadmap for discretely energy-stable schemes for dissipative
  systems based on a generalized auxiliary variable with guaranteed positivity.
\newblock {\em Journal of Computational Physics}, 404:109121, 2020.
\newblock (also arXiv:1904.00141).

\bibitem{ZangBYZ2020}
Y.~Zang, G.~Bao, X.~Ye, and H.~Zhou.
\newblock Weak adversarial networks for high-dimensional partial differential
  equations.
\newblock {\em Journal of Computational Physics}, 411:109409, 2020.

\bibitem{ZhengD2011}
X.~Zheng and S.~Dong.
\newblock An eigen-based high-order expansion basis for structured spectral
  elements.
\newblock {\em Journal of Computational Physics}, 230:8573--8602, 2011.

\end{thebibliography}

\end{document}